\documentclass[10pt]{amsart}

\usepackage{lineno}

\usepackage{hyperref}
\usepackage{comment}

\usepackage{dsfont}
\usepackage{amsmath}
\usepackage{amssymb}
\usepackage{graphicx}

\usepackage{amssymb}

\usepackage[all,cmtip]{xy}

\usepackage{amsmath}

\usepackage{amsfonts}

\usepackage{soul}

\usepackage{mathpazo}
\usepackage{color}
\usepackage{yfonts}

\usepackage{paralist}

\usepackage{stmaryrd}

\usepackage{amsxtra}

\def\cU{\mathcal A}

\def\cA{\mathcal A}

\def\cI{\mathcal I}

\newcommand{\norm}[1]{\| #1\|}
\def\b{          \beta}

\def\a{          \alpha}

\def\cA{         I}

\def\cB{          \mathcal B}

\def\cD{          \mathcal D}
\def\cF{          \mathcal F}
\def\cH{          \mathcal H}

\def\tWu{       \widetilde W^u}
\def\cX{   {\cal Z}}
\def\cY{   {\cal Y}}

\def\clb{   \color{black}}

\let\cal\mathcal

\def \R{{\mathbb R}}

\def \Z{{\mathbb Z}}
\def \N{{\mathbb N}}

\def \a{{\alpha}}

\def\tpsi{   \widetilde\psi}
\def\tphi{   \widetilde\phi}

\newcommand{\prf}{{\begin{proof}}}
\newcommand{\epf}{{\end{proof}}}

\newcommand{\cT}{{\mathcal T}}

\newcommand{\PP}{{\mathbb P}}
\newcommand{\cP}{{\mathcal P}}

\newcommand{\cK}{{\mathcal K}}

\newcommand{\C}{{\mathbb C}}

\newcommand{\cQ}{{\mathcal Q}}

\newcommand{\orbit}{{\mathcal O}}

\newcommand{\ary}{\begin{eqnarray}}
\newcommand{\eary}{\end{eqnarray}}

\newcommand{\aryst}{\begin{eqnarray*}}
\newcommand{\earyst}{\end{eqnarray*}}

\newcommand{\enmt}{\begin{enumerate}}
\newcommand{\eenmt}{\end{enumerate}}

\newtheorem{thm}{\bf Theorem}[section]

\newtheorem{lemma}[thm]{\bf Lemma}
\newtheorem{prop}[thm]{\bf Proposition}
\newtheorem{cor}[thm]{\bf Corollary}

\newtheorem{claim}[thm]{\bf Claim}

\theoremstyle{definition}
\newtheorem{defi}[thm]{\bf Definition}
\newtheorem{rema}[thm]{\bf Remark}

\newcommand{\rmp}{\varrho_2}
\newcommand{\rhone}{\varrho_1}
\newcommand{\rhothree}{\varrho_3}

\theoremstyle{definition}

\def\bee{\begin{equation}}
\def\eee{\end{equation}}

\newcommand{\pdvr}[2]
{\dfrac{\partial^{#2} #1}{\partial \theta^{#2_1} \partial r^{#2_2}}}
\newcommand{\pdvrs}[2]
{\partial^{#2} #1 /\partial \theta^{#2_1} \partial r^{#2_2}}
\newcommand{\rdiv}{{\rm div}}

\newcommand{\cL}{\mathcal{L}}

\newcommand{\cM}{\mathcal{M}}

\numberwithin{equation}{section}

\definecolor{blue}{rgb}{0,0,1}
\definecolor{red}{rgb}{1,0,0}
\definecolor{green}{rgb}{0,.7,0}

\begin{document}
	
	\title[exponential mixing of 3D Anosov flows]{Smooth mixing Anosov flows in dimension three are exponential mixing}
	

	\author{Masato Tsujii}
	\address{Masato Tsujii, Department of Mathematics, Kyushu University, Fukuoka, 819-0395}
	\email{tsujii@math.kyushu-u.ac.jp}
		
	\author{Zhiyuan Zhang}
	\address{Zhiyuan Zhang, CNRS, Institut Galil\'ee
		Université Paris 13
		99, avenue Jean-Baptiste Cl\'ement
		93430 - Villetaneuse}
	\email{zhiyuan.zhang@math.univ-paris13.fr}
		
	\maketitle
	
	\tableofcontents
	
%	\linenumbers

\section{Introduction}

A $C^1$ flow $g = \{  g^t: M \to M \}_{t \in \R}$ on a compact Riemannian manifold $M$ is {\it Anosov} if there is a continuous flow-invariant splitting
$TM = E^s \oplus N \oplus E^u$ such that $N$ integrates to the flow-line foliation, and  there exist $C, \kappa > 0$ such that for all $t > 0$
\aryst
\norm{Dg^t|_{E^s}} < Ce^{-\kappa t}, \quad 
\norm{Dg^{-t}|_{E^u}} < Ce^{-\kappa t}.
\earyst
%Without loss of generality, we always assume that, after certain smooth change of metric, we have $C = 1$ in the above inequalities.

Ergodic properties of Anosov flows has been studied extensively in the history of dynamical system theory, one because they exhibit a type of chaotic behavior of orbits.  
After the pioneering works of Hadamard \cite{Hadamard} and Hopf \cite{Hopf}, Anosov proved in his thesis \cite{An} that any volume-preserving Anosov flows are ergodic and then Anosov, Sinai \cite{Sinai, AnosovSinai} proved that they are mixing unless the stable and unstable foliations are jointly integrable. 

The Bowen-Ruelle conjecture states that a $C^r$ mixing Anosov flow should be exponential mixing with respect to the so-called equilibrium measures with H\"older potentials.\footnote{Originally, Bowen and Ruelle asked in \cite{BR} whether exponential mixing holds for topologically mixing Axiom A flows. However it turns out that a suspension flow over a full shift with locally constant roof function is never exponential mixing, see \cite{Rue}. The conjecture for Anosov flows remains open.} 
This conjecture seems natural because an Anosov flow strongly mixes in the directions transversal to the flow lines and therefore 
it is difficult to imagine any reason that prevents exponential mixing once the flow started to mix. However, by technical difficulties, there had not been much progress until the epoch-making paper of Dolgopyat \cite{Dol}. He established an argument, which is called Dolgopyat's argument nowadays, to derive exponential mixing (or rapid mixing) of hyperbolic flows from some estimates on joint non-integrability of the stable and unstable foliations. In particular, he proved that an  Anosov flow is exponential mixing if the stable and unstable foliations are $C^1$ and not jointly integrable, and the equilibrium measure satisfies certain doubling property. For instance, a geodesic flow on a surface with negative curvature satisfies these conditions.
After Dolgopyat's work, there is much literature on related researches. 
See \cite{ABV,L,GLP, PS, BDL, FMT} for instance.

Actually the results of Dolgopyat mentioned above were still not very close to the resolution of the Bowen-Ruelle conjecture, because the $C^1$ assumption on invariant foliations is very strong and does not hold usually. 
In order to approach the conjecture, it has remained to develop detailed geometric analysis on non-integrability between the stable and unstable foliations. This is indeed the point that we focus on in this paper. 
A few years ago, the first author proved in \cite{Tsu} that generic $3$-dimensional volume-preserving Anosov flows are exponential mixing with respect to the volume measure. The key idea there was to introduce \lq\lq tempelate functions\rq\rq which describe variations of stable subspaces $E^s$ along unstable manifolds in some intrinsic manner and measure the non-integrability of the stable and unstable foliations. In this paper, we develop the idea much further and, together with other technical developments, we provide a positive answer to the Bowen-Ruelle conjecture in three dimensional cases, though we are still restricted to the case $r=\infty$ \footnote{ Actually we only need to assume that the flow is $C^r$ for some $r$ sufficiently large depending only on the ratio of the exponents in \eqref{eq rate stable unstable}.  \clb}.

\subsection{Main results}

Throughout this paper, we fix a $3$-dimensional compact Riemannian manifold $M$.
	Let $g$ be a $C^\infty$ transitive Anosov flow on $M$.
Given $\theta > 0$ and $F \in C^\theta(M)$, 
we denote by $\nu_{g,F}$ the unique {\it equilibrium measure} for $g$ with potential $F$ (the precise definition is given in Definition \ref{def gibss}). We say that $g$ is {\it exponentially mixing} with respect to $\nu_{g, F}$ and all H\"older test functions if 
for any $\beta > 0$,  there is $\kappa = \kappa(g, \beta) > 0$ such that for any $A, B \in C^\beta(M)$ there exists $C = C(g, A, B,F) > 0$ such that
\aryst
|\int A \circ g^t B d\nu_{g,F} - 
\int A d\nu_{g,F} 
\int  B d\nu_{g,F}| < C e^{-t\kappa}, \quad \forall t > 0.
\earyst
It is straightforward to see that $g$ is exponentially mixing  with respect to $\nu_{g, F}$ and all H\"older test functions as long as the above estimate holds for {\it some }$\beta > 0$.

The main result of this paper is the following.

\begin{thm}\label{main thm}
	Let $g$ be a $C^\infty$ transitive Anosov flow on a $3$-dimensional compact manifold $M$. Then the following are equivalent:
	\enmt
	\item $g$ is topologically mixing;
	\item $g$ is exponentially mixing with respect to any equilibrium measure with H\"older potential and all H\"older test functions. 
 \eenmt
\end{thm}
It is clear that in the above theorem, (2) implies (1). In the rest of the paper we will show that (1) implies (2) as well.

Theorem \ref{main thm} is deduced from Dolgopyat's estimate (Proposition \ref{prop decay near axis}). Another consequence of Proposition \ref{prop decay near axis} is the error term estimate in the {\em prime orbit theorem} for $3D$ Anosov flows. 
For any periodic orbit $\orbit$ of $g$, we denote by $l(\orbit)$ the prime period of $\orbit$. We denote
\aryst
\pi(T) = |\{ \orbit \mid l(\orbit) \leq T  \}|, \quad \forall T > 0.
\earyst
The usual prime orbit theorem, due to Parry and Pollicott \cite{PP}, gives a more precise estimate in the case of topologically weakly mixing hyperbolic flows:
\aryst
\pi(T) = (1+o(1))li(e^{h_{top}T}) \ \mbox{ as } \ T \to +\infty
\earyst
where 
\aryst
li(y) = \int_2^y \frac{1}{\log u} du \sim \frac{y}{ \log y}  \ \mbox{ as } \ y \to +\infty. 
\earyst
Earlier Margulis \cite{Mar} obtained this result in the context of Anosov flows.
For geodesic flows on a surface with negative (variable) curvature, Pollicott and Sharp proved in \cite{PS} that the error term $o(1)$ in the above formula is actually exponentially small. More recently, this result is extended to higher dimensional contact Anosov flows by Giulietti, Liverani and Pollicott \cite{GLP} and Stoyanov \cite{Sto} under some additional conditions. The first author  obtained in \cite{Tsu4} more explicit estimate on the error term for expanding semiflows under a generic condition. In this paper we show the following.
\begin{thm}\label{thm primeorbit}
	Let $g$ be a $C^\infty$ topologically mixing Anosov flow on a $3$-dimensional compact manifold. Denote $h = h_{top}(g^1)$. Then there exists $0 < c < h$ such that
	\aryst
	\pi(T) = li(e^{h T}) + O(e^{cT})
	\earyst
\end{thm}

Both the proofs of Theorem \ref{main thm} and Theorem \ref{thm primeorbit} will be given at the end of Section \ref{sec markov}.

\subsection{Plan of the proof}
Below we explain the idea behind the proof and how we will proceed in the following sections. Note that the explanation is necessarily rather  schematic and the rigorous argument in the later sections will be slightly different. 

Let $X$ be the generator of the transfer operators that describes the natural action of the flow on the space functions. The main point of Dologpyat's argument is that once we prove uniform boundedness of the resolvent 
\begin{equation}\label{eq:resolvent}
\mathcal{R}(s)=(s-X)^{-1}=\int_0^\infty e^{-st}e^{tX} dt
\end{equation}
for $s=a+bi$ with $|a|$ sufficiently small and $|b|$ sufficiently large, we get exponential mixing of the flow\footnote{In order to realize this idea directly, we would need to consider some functions spaces called anisotropic Sobolev (or Banach) spaces. But we will take an alternative way as in \cite{Dol} that uses reduction by Markov partition which we will describe in Section \ref{sec markov}.}. 
We will therefore estimate the integral in \eqref{eq:resolvent}. 
At this moment, we emphasize that we will consider about the parameter $b$ in parallel and that uniformity of constants and estimates on the parameter $b$ is extremely important (note also that we will write $\epsilon = 1/|b|$).

In order to prove uniform boundedness of the resolvent $\mathcal{R}(s)$, we study the temporal distance function defined as follows. Consider a point $x\in M$ and its local (strong) stable and unstable manifold $W^s_{g, loc}(x)$ and $W^u_{g, loc}(x)$. For points $z\in W^s_{g, loc}(x)$ and $y\in W^u_{g, loc}(x)$, there  a point $q(z,y)$ and a real number $\tau_x(z,y)$ so that $q(z,y)\in W^s_{g, loc}(x)$ and $f^{\tau_x(z,y)}(q(z,y))\in W^{u}_{g, loc}(y)$. The function $\tau_x(z,y)$ is called the temporal distance function.
The key fact to prove is that the function $\exp(ib\tau_x(s,t))$ oscillates when $(s,t)$ varies and the oscillation is somewhat uniform with respect to point $x$ and also for the parameter $b$ when $|b|$ is sufficiently large. 

In Section 2, we provide some basic definitions for Anosov flows and then, 
in Section 3, we explain a precise framework for the proof of the main estimate, Proposition \ref{prop decay near axis}.
At the end of Section 3, we present two propositions, Proposition \ref{prop main 1} and \ref{prop main 2}.
The former, Proposition \ref{prop main 1}, claims roughly that joint non-integrability of the stable and unstable foliations implies some uniform estimates on the oscillation of the term $\exp(ib\tau_x(s,t))$. 
But note that, as we will explain below, the uniformity only holds on some subset $\Omega\subset M$ and is viewed in an appropriate scale that depends on point. The latter, Proposition \ref{prop main 2},  
claims that, once we obtain the estimates in the conclusion of Proposition \ref{prop main 1}, we obtain exponential mixing of the flow. The main theorem is an immediate consequence of these propositions. Section \ref{sec Temporal function and non-integrability} to \ref{sec checking the conditions} will be devoted to the proof of Proposition \ref{prop main 1}. The rest will be devoted to the proof of Proposition \ref{prop main 2}. 

In Section \ref{sec Temporal function and non-integrability}, we give the main estimates on temporal distance function $\tau_x(z,y)$. 
To this end, we introduce a convenient flow box coordinates around each point $p\in M$ in Lemma \ref{lem normalcoordinatesystem} and Lemma \ref{lem def of gammax}. These are a kind of ``normal coordinates" that are associated to hyperbolic structure of the flow in some canonical manner. 
Then we will see how the stable and unstable manifolds are transformed by the flow in such coordinates. Not surprisingly,  the transformations are well approximated by the Taylor expansions of some high order, say order $K$, at the points along the orbit of $x$. 
Indeed we obtain a nice approximation of the temporal distance function modulo polynomials of order $K$. 
This approximation, summarized in Corollary \ref{prop template approximation}, is given in terms of the ``template function'' along stable and unstable manifolds, introduced in \cite{Tsu}. 

In Section \ref{sec proof of thm genericity of uni}, we show that non-integrability of the stable and unstable foliation implies some  oscillation of the term $\exp(ib\tau_x(s,t))$ uniform  with respect to the point $x$ and the parameter $b$. In the main step of the proof, we  use Journ\'e's regularity lemma \cite{Jou}. From this lemma,  if the stable subspace $E^s$ does not depend on points $W^u_{g,loc}(x)$ smoothly, then the stable and unstable foliations are not jointly integrable. 
Then we show uniformity of such irregularity using mixing property of the flow. 
When $E^s$ is smooth, we can use the mixing hypothesis to show the uniformity of non-integrability. 

In Section \ref{sec construct partition}, we will introduce a scale depending on the global parameter $\epsilon=1/|b|$, at which we observes the oscillation of the term $\exp(ib\tau_x(s,t))$. Though the estimates in the previous section are basically sufficient, we face another technical problems caused by the fact that the oscillation may be sometimes too fast\footnote{In estimating the oscillating integral, this is usually not a problem. But since we partly consider approximation of temporal distance function in $C^0$ norm, it appears to be a problem.}.  
To avoid this problem, we set up an appropriate scale depending on points (we call it the scale function) and on the parameter $\epsilon=1/|b|$. The definition of the scale function, denote by $\Lambda^\epsilon$, is given in Definition \ref{def hepsilon}. We will then prove that the scale is \emph{tame} in the sense that the oscillation of $\exp(ib\tau_x(s,t))$ looks tame at scale $(\Lambda^\epsilon)^{-1}$. 
Also we can check that the scale is \emph{stable} in the sense that it does not change too much along orbits so that the flow remains  expanding (in the unstable direction) with respect to it. 

In Section \ref{sec Uniform non-integrability on uniform set} and \ref{sec recurrence}, we discuss about the properties of $\Lambda^\epsilon$. A difficulty with the scale function $\Lambda^\epsilon$ is that its dependence on point is not obvious and will not be continuous. In Definition \ref{def uniform set}, we introduce a subset $\Omega$ (more precisely subsets $\Omega(n,\kappa)$ and $\Omega(\epsilon,n,\kappa)$) on which the derivative of the flow satisfies some mild condition on its growth along forward orbit. 
Then, in Lemma \ref{lem oscillation lower bound and smoothness upper bound}, we show that we observe sufficient oscillation of the temporal distance function if the base point $x$ is sufficiently close to the subset $\Omega$. Also, in Lemma \ref{lem verifyadapted}, we show that the scale $\Lambda^\epsilon$ is \emph{adapted} to the subset $\Omega$. This claim is a rather technical one and its meaning may not be very clear at this point. But this will play an important roll in the argument in the proof of Proposition \ref{prop main 2}. Then we will prove in Lemma \ref{lem fast recurrence for expanding set} that the orbits of the flow return often to $\Omega$. 
Finally, in Section \ref{sec checking the conditions}, we prove Proposition \ref{prop main 1} summarizing arguments in the previous sections.

We start the proof of Proposition \ref{prop main 2} from Section \ref{sec proofofpropmaincriterion}. The proof follows the argument in \cite{Dol} but with some modification. 
The point of the modification is that we have to observe the transfer operator at the scale $(\Lambda^\epsilon)^{-1}$. 
On the one hand the scale $\Lambda^\epsilon$ (and related constructions) depends on the parameter $\epsilon=1/|b|$ and, on the other hand, we need to get a uniform bound with respect to the parameter $b$ for the decay of the operator $e^{-st}e^{tX}$ in \eqref{eq:resolvent}. So we have to keep track of dependence of constants on the parameter $b$ carefully. However, once we keep this in mind, the proof is not too difficult because most of the necessary estimates have already been given in the conclusion of Proposition \ref{prop main 1}. The main line of the proof is presented in Proposition \ref{prop L2 iteration}. Then the key induction step in the proof of Proposition \ref{prop L2 iteration} is presented in Proposition \ref{lem dim cancellation}. 
In Section \ref{sec partition}, we introduce a refinement of the Markov partition associated to  $\Lambda^\epsilon$.  In Section \ref{sec: bounds for smoothness}, we prove some simple decay estimates on H\"older norm.  Then, finally, we prove Proposition \ref{lem dim cancellation} in Section \ref{sec proof lem dim cancellation}, completing the proof of the main theorem. At the end, we provide appendices where we collect some elementary estimates that are used in the text. 

\subsection*{Acknowledgement}
The authors would like to thank Alex Eskin and Rafael Potrie for long and stimulating discussions and helpful inputs on the non-integrability condition and many more.
We thank Viviane Baladi, Clark Butler, Dmitry Dolgopyat and Mark Pollicott for discussions on the doubling property of equilibrium measures. 
We also thank useful discussions with  Artur Avila,
Julien Barral, Fran\c{c}ois B\'eguin, Bassam Fayad, Sylvain Crovisier, Simion Filip, Federico Rodriguez-Hertz,  Martin Leguil, Yuri Lima, Carlos Matheus, Amir Mohammadi, Yi Pan, Peter Sarnak, Khadim War, Amie Wilkinson, Jinxin Xue and Lai-Sang Young.

\section{Basic properties of Anosov flow}\label{sec basicpropertiesofanosovflow}

For $1 \leq r \leq \infty$, let $\cU^r$ denote the set of $C^r$ transitive  Anosov flows on  $M$.
For each  $g \in \cU^r$, we denote the generating vector field of $g$ by $V_g$. We denote by the stable subspace, resp. unstable subspace, of $g$ by $E^s$, resp. $E^u$. It is known that $E^s$ and $E^u$ uniquely integrate to the stable foliation $W^s_g$ and unstable foliation $W^u_g$ respectively. 
We denote by $W_g^{cu}$, resp. $W_g^{cs}$, the center-unstable, resp. center-stable foliation.

Given $g \in \cU^1$, we may suppose, after changing the Riemannian metric if necessary, that there exist $\chi_s, \chi_u > 0$, $\bar\chi_u > \chi_u$, $\bar\chi_s > \chi_s$  such that for any $t > 0$,
\ary \label{eq rate stable unstable}
e^{-\bar\chi_s t} < \norm{Dg^t|_{E^s}} < e^{-\chi_s t},  \quad 
e^{-\bar\chi_u t} < \norm{Dg^{-t}|_{E^u}} < e^{-\chi_u t}.
\eary
We set
\ary \label{eq chi *}
\chi_0 = \min(\chi_u, \chi_s), \quad
\chi_*  = \max(\bar\chi_u, \bar\chi_s).
\eary
For any integer $n \geq 0$, for any $x \in M$, we denote
\aryst
\Lambda_n(x) = \norm{Dg^n|_{E^u(x)}}.
\earyst

Now assume that  $g \in \cU^r$ for some $r \geq 2$. 
Let $\upsilon_* > 0$ be a small constant which will be determined later.
It is known that there is a collection of {\it non-stationary normal coordinates} $\{ \Phi^u_x: \R \to W^u_g(x) \}_{x \in M}$ such that for each $x \in M$:
\enmt
\item $\Phi^u_x$ is a $C^r$ diffeomorphism, $\Phi^u_x(0) = x$ and $\norm{D\Phi^u_x(0)} = \upsilon_*$; 
\item $(\Phi^u_{g^t(x)})^{-1} g^t \Phi^u_x(s) =  \pm  \norm{Dg^t|_{E^u(x)}}s$ for all $s,t \in \R$.
\eenmt
We refer the readers to \cite{KK} for a proof of the existence of such coordinate systems.
For any $x \in M$, we choose $\lambda_x$ and $\mu_x$ such that for all $s \in \R$,
\ary
(\Phi^u_{g^1(x)})^{-1} g^1 \Phi^u_x(s) &=&  \lambda_x^{-1}s, \label{def of lambdax} \\
(\Phi^s_{g^1(x)})^{-1} g^1 \Phi^s_x(s) &=&  \mu_x^{-1}s. \label{def of mux}
\eary
By our choices, we have
\ary \label{eq boundsforrates}
|\lambda_x| < e^{-\chi^u} < 1 < e^{\chi^s} < |\mu_x|, \quad \forall x \in M.
\eary

For each $x \in M$ and $\varrho > 0$, we set $W^u_{g}(x, \varrho) = \Phi^u_x((-\varrho,  \varrho))$.  More generally,  for $x \in M$ and an interval $J \subset \R$, we denote
\aryst
W^u_g(x, J) = \Phi^u_x(J).
\earyst
We set $W^u_{g, loc}(x) = W^u_{g}(x, 50\varrho_*)$ for all $x \in M$ where $\varrho_* > 1$ is a large constant so that several intersections in the rest of the paper are well-defined.  
For instance, we let $\varrho_*$ be sufficiently large depending only on $g$, and let $\upsilon_*$ be sufficiently small depending only on $g$ and $\varrho_*$, so that the intersection between $W^{s}_{g, loc}(y)$ and $W^{cu}_{g, loc}(z)$ for any $y,z$ within distance $20\varrho_*\upsilon_*$ uniquely exists.
We define $\Phi^s_x$, $W^s_{g, loc}(x)$, $W^s_{g}(x, J)$ in a similar way. 

We let $\upsilon_*$ be sufficiently small so that for any $x \in M$, for $* = s, u$,
for any $y \in W^*_{g, loc}(x)$, we have
\ary \label{eq local metric comparable}
 \frac{1}{2} < | D[(\Phi^*_y)^{-1} \Phi^*_x](s) | < 2, \quad \forall s \in (-\varrho_*, \varrho_*).
\eary

\begin{defi}
Each element of $\cU^2$ falls into one of the following three non-exclusive classes:
\enmt 
\item[${\it I}.$] Neither $E^u$ nor $E^s$ for $g$ is $C^{1+\delta}$ for any $\delta > 0$; 
\item[${\it II}.$] $E^u$ for $g$ is $C^{1+\delta}$ for some $\delta > 0$;
\item[${\it III}.$] $E^s$ for $g$ is $C^{1+\delta}$  for some $\delta > 0$.
\eenmt
\end{defi}

\begin{rema} \label{rema only class i-ii}
Since the validity of Theorem \ref{main thm} for all $g$ in Class {\it II} implies the validity of Theorem \ref{main thm} for all $g$ in Class {\it III} (by considering the reversed flow), we will only consider the case where $g$ is in Class {\it I} or {\it II} from now on. 
\end{rema}

\begin{rema} 
	We notice that exponential mixing for Anosov flow in Class {\it III}  with respect to the SRB measure is already obtained in \cite{BW}. There the authors only need to assume that the flow is $C^{1+\delta}$. However they have used integration-by-part which only works when $\nu_U$ is absolutely continuous. As we need to consider general equilibrium measures, we give an independent proof. \clb
\end{rema}

The following theorem due to Plante \cite{P} gives a useful characterization of topologically mixing $3D$ Anosov flow.
\begin{thm}[Theorem 3.1 in \cite{P}]\label{thm mixing equals nonintegrable}
	A  transitive $3D$ Anosov flow $g$ is topologically mixing if and only if $E^s$ and $E^u$ are not jointly integrable. 
\end{thm}
Here we say that $E^s$ and $E^u$ are {\em jointly integrable} if there is a $C^1$ foliation whose tangent bundle is $E^u \oplus E^s$. We refer the readers to \cite[page 734-735]{P} for relevant notions.

	For the rest of this paper we make the following notational conventions: 
	\enmt
	\item[$\bullet$] 
	given an Anosov flow $g$,
	we let $C' > 1$ denote a generic constant which depends on $g$, and possibly also on a given potential function $F$ (see Definition \ref{def gibss}), but is independent of anything else.
	We assume that $C'$ may vary from line to line. We say that two positive constants $A,B$ satisfy $A \sim B$ if $A/B \in (C'^{-1}, C')$. We denote by  $O(A)$ a constant of absolute value bounded from above by $C'|A|$.
    We let $\delta \in (0,1)$ be a small constant,  which may vary from line to line, depending on $g$, and possibly also on a given potential $F$, but is independent of anything else.
	\item[$\bullet$] 
	for any function $f: J \to \R$ defined on a subset $J \subset \R^n$ for some $n \geq 1$, we denote
	\aryst
	Osc_{s \in J} f = \sup_{s \in J} f(s) - \inf_{s \in J} f(s),  \quad
	\norm{f}_J = \sup_{s \in J} |f(s)|.
	\earyst
	When the domain $J$ is clear from the context, we may abbreviate $\|\cdot \|_J$ as $\|\cdot\|$ or $\|\cdot \|_{C^0}$.
	\item[$\bullet$] 	for any integer $d \geq 1$, we denote by $Poly^d$ the set of real polynomials  with degree less or equal to $d$. 
	We write $Poly^d(s)$ when we indicate that the independent variable is denoted by $s$.
	For every $C > 0$
	we denote by $Poly^{d}_C$ the collection of polynomials in $Poly^{d}$ with all coefficients of absolute values at most $C$.
	
	Similarly, for any integers $d_1, d_2 > 0$, we denote by  $Poly^{d_1, d_2}$  the set of real polynomials of the form 
	\aryst
	(s,t) \mapsto \sum_{i=0}^{d_1} \sum_{j = 0}^{d_2} a_{i,j} s^i t^j. 
	\earyst
	For every $C > 0$, we denote by $Poly^{d_1, d_2}_C$ the collection of polynomials in $Poly^{d_1, d_2}$ with all coefficients of absolute values at most $C$.
	\eenmt

\section{Markov partition and complex RPF operator} \label{sec markov}

In this section, we recall some properties of the Markov partition and the so-called complex Ruelle-Perron-Frobenius operators for Anosov flows. We follow closely the presentations in \cite{Dol}.

Throughout this section, we let $g$ denote an Anosov flow in $\cU^2$,
and let $\theta \in (0,1)$ denote a constant to be determined depending only on $g$ and a potential $F$ introduced in Section \ref{subsec eqmes}.

\subsection{Markov partition, roof function and H\"older space}
By \cite{Bow, Rat}, any Anosov flow $g$ admits a Markov partition $\Pi = \bigcup_{\a \in \cA} \Pi_\a$ where $I$ is a finite set (hereafter we use the notations in \cite{Dol}).
For each $\a \in \cA$,  $\Pi_\a$ is a parallelogram defined as follows: there is  a local unstable manifold $U_\a$, and a local stable manifold $S_\a$, both of diameters bounded by $\varrho_* \upsilon_*$, such that
\aryst
\Pi_\alpha = \{ [x,y] \mid x \in U_\a, y \in S_\a \}
\earyst
where $[x,y]$ is the unique intersection of $W^s_{g, loc}(x)$ and $W^{cu}_{g, loc}(y)$.  
We define $x_\a$ to be the unique intersection between $U_\a$ and $S_\a$. 
For different $\a, \b \in \cA$, $\Pi_\a$ and $\Pi_\b$ are disjoint. 
We require that $\Pi$ satisfies the following Markov property.
Let $\hat\sigma$ be the Poincar\'e map (of the flow $g$) on $\Pi$, and  let $\tau$ be the return time function. That is, we have for any $x \in \Pi$ that
\ary
\tau(x) &=& \inf \{ t > 0 \mid g^t(x) \in \Pi \}, \\
\hat{\sigma}(x) &=& g^{\tau(x)}(x).
\eary
Then for each $\a \in \cA$, $\hat\sigma(S_\a)$ is contained in $\Pi_{\b_1}$ for some $\b_1 \in \cA$; and  $\hat\sigma^{-1}(U_\a)$ is contained in $\Pi_{\b_2}$ for some $\b_2 \in \cA$.

We denote 
\aryst
U = \bigsqcup_{\a \in \cA} U_\a, \quad S = \bigsqcup_{\a \in \cA} S_\a.
\earyst
We denote by $\pi_U$ the projection from $\Pi$ to $U$: for each $\a \in \cA$, for each $x \in U_\a$ and $y \in S_\a$, we set
\aryst
\pi_U([x,y]) = x.
\earyst

We define $\sigma: U \to U$ to be the $U$ coordinate of $\hat \sigma$.\footnote{The map $\sigma$ is well-defined as the stable manifold foliation is invariant under the flow.} We set
\ary \label{eq tau *}
 \tau_0 = \inf \tau, \quad \tau_*  = \sup \tau.
\eary
Without loss of generality, we may let the diameter of each $\Pi_\a$  be sufficiently small, and let $\rmp > 0$ be a small constant depending on $g$ and $\Pi$ 
 so that for any $\a \in \cA$ and any $x \in U_\a$, 
\ary \label{eq wsgx1 is long}
W^s_g(x, \rmp) \subset \Pi_\a \cap W^s_{g, loc}(x) \subset W^s_g(x, 1/2).
\eary

For each integer $n > 0$, each $\alpha \in I$ and each $x \in U_\alpha$, we denote by $\sigma^{-n}(\alpha)$ the set of inverse branches of $\sigma^n$ restricted to $U_\alpha$, and write $\sigma^{-n}_x = \sigma^{-n}(\alpha)$. 
In this case, the domain of each $v \in \sigma^{-n}(\a)$, denoted by $Dom(v)$, is $U_\alpha$.
The usual caveat applies: we consider an inverse branch $v \in \sigma^{-n}(\a)$ to be continuous on $U_\a$, but we only require that the equality $\sigma^n v = {\rm Id}$ holds in the interior of $U_\a$. 
Finally, we denote
\aryst
\sigma^{-n} = \bigsqcup_{\alpha \in I} \sigma^{-n}(\alpha).
\earyst

By letting $\theta > 0$ be a sufficiently small constant, we may assume that $W^{cs}_g$ is a $C^{1+\theta}$-foliation (see \cite{Ha}). As a result, the center-stable holonomy between any two nearby unstable leaves is $C^{1+\theta}$.
Consequently, for every $\a \in \cA$, every inverse branch $v \in \sigma^{-1}(\a)$ is a $C^1$ contracting map from $U_\a$ to some $U_\b$. 
We have
\ary \label{contraction}
e^{-n\chi_*} d_{W^u_g}(x, y)/C' < d_{W^u_g}(v(x), v(y)) < C' e^{-n \chi_0} d_{W^u_g}(x, y)
\eary
for any integer $n \geq 1$, any $\a \in \cA$, any inverse branch $v \in  \sigma^{-n}(\a)$ and any $x,y \in U_\a$.

We define the uniform norm $\norm{\cdot}_{C^0}$ and the space $C^0(U, \C)$, or simply abbreviated as $C^0(U)$, of continuous complex-valued functions on $U$ in the usual way. We denote by $C^0(U, \R)$ the subset of real-valued functions in $C^0(U, \C)$.
We define  $C^\theta(U, \R)$, resp. $C^\theta(U, \C)$ or simply $C^\theta(U)$, to be the collection of  $\varphi \in C^0(U, \R)$, resp. $C^\theta(U, \C)$, such that 
\aryst 
|\varphi|_{\theta} = \sup_{x \neq y \in U_\a, \a \in I} \frac{|\varphi(x) - \varphi(y)|}{d_{W^u_g}(x,y)^\theta} < \infty.
\earyst
We set 
\aryst
\norm{\varphi}_\theta = \norm{\varphi}_{C^0} + |\varphi|_{\theta}.
\earyst
With a slight abuse of notation, for functions in $C^\theta(J, \C)$ where $J \subset \R$ we also use $\norm{\cdot}_\theta$, resp. $|\cdot|_\theta$ or $|\cdot|_{\theta, J}$, to denote the $\theta$-H\"older norm, resp. semi-norm. 

For any $\varphi \in C^\theta(U, \C)$,
for any integer $n \geq 0$ we denote
\aryst
\varphi_{n}(x) = \sum_{i=0}^{n-1} \varphi \circ	 \sigma^i(x), \quad \forall x \in U.
\earyst

\subsection{Return time and temporal function}\label{lab Temporal function and templates}

\begin{defi} \label{def temporal}
	For any $x \in M$ and $y \in W^s_{g}(x, 5\varrho_*)$, the {\em temporal function} $\Psi_{x,y}: W^u_{g}(x, 5\varrho_*) \to \R$ is defined as follows. For any $z \in W^u_{g}(x, \varrho_*)$, we let $z'$ be the unique intersection of $W^{cu}_{g, loc}(y)$ and $W^s_{g, loc}(z)$; let $z''$ be the unique intersection of $W^{cs}_{g, loc}(z)$ and $W^u_{g, loc}(y)$; and define
	$\Psi_{x,y}(z)$ by equation
	\aryst
	z' = g^{\Psi_{x,y}(z)}(z'').
	\earyst 
\end{defi}
It is known that temporal functions are H\"older continuous, with H\"older exponents and H\"older norms uniformly bounded  depending only on $g$.

The following lemma is an immediate consequence of Definition \ref{def temporal} and the flow invariance of the stable/unstable manifolds. We omit its proof.

\begin{lemma} \label{lem shift invariance of temporal function}
	Let $x \in M$, $y_1, y_2 \in W^s_{g}(x, \varrho_*)$, $z_1,z_2 \in W^u_{g}(x, \varrho_*)$. Then the following are true:
	\enmt
	\item for any $t > 0$ with $g^t(z_1) \in W^u_g(g^t(x), 5\varrho_*)$, we have
	\aryst
	\Psi_{x, y_1}(z_1) = \Psi_{g^t(x), g^t(y_1)}(g^t(z_1));
	\earyst
	\item
	let $y'_1$ be the unique intersection of $W^s_{g, loc}(z_1)$ and $W^{cu}_{g, loc}(y_1)$, then
	\aryst
	\Psi_{x,y_1}(z_2) - \Psi_{x,y_1}(z_1) = \Psi_{z_1, y'_1}(z_2);
	\earyst
	\item
	let $z'_1$ be the unique intersection of $W^u_{g, loc}(y_1)$ and $W^{cs}_{g, loc}(z_1)$, then
	\aryst
	\Psi_{x,y_2}(z_1) - \Psi_{x,y_1}(z_1) = \Psi_{y_1, y_2}(z'_1).
	\earyst
	\eenmt
\end{lemma} 
The return time function $\tau$ is closely related to the temporal functions. 
Notice that $\tau([x,y])$ is independent of $y$ since the stable manifold foliation is invariant under the flow. Hence $\tau$ can be treated as a function on $U$.
We have the following.
\begin{lemma} \label{eq b tau k}
	For any $\a \in \cA$, any $k \geq 1$, any $w \in \sigma^{-k}(\a)$,  there exists a unique $x_{w} \in S_\a$ such that 
	$g^{\tau_k(w(x_\a))}(w(x_\a)) = x_w$. 
	Consequently
	\ary 
	\tau_k \circ  w(z) -  \tau_k \circ  w(x_\a)  =  \Psi_{x_\a, x_{w}}(z), \quad \forall z \in U_\a. 
	\eary
\end{lemma}

\begin{cor} \label{cor eq b tau k}
	For any $\a \in \cA$, any $x \in U_\a$, any $k \geq 1$, any $w \in \sigma^{-k}(\a)$, the unique intersection $x^w$ of $W^{s}_{g, loc}(x)$ and $W^{cu}_{g, loc}(x_w)$ satisfies 
	\aryst
	\tau_k \circ w(z) - \tau_k \circ w(x) = \Psi_{x, x^w}(z), \quad \forall z \in U_\a.
	\earyst
	Moreover by letting $\theta$ be sufficiently small depending only on $g$ we have
	\ary \label{eq tau n circ v f n circ v}
	\norm{\tau_k \circ w }_{\theta} < C'.
	\eary	
\end{cor}
\begin{proof}
	By Lemma \ref{eq b tau k}, we have
	\aryst
	\tau_k \circ w(z) - \tau_k \circ w(x) = \Psi_{x_\a, x_w}(z) - \Psi_{x_\a, x_w}(x), \quad \forall x, z \in U_\a.
	\earyst	
	By Lemma \ref{lem shift invariance of temporal function}(2), we have
	\aryst
	\Psi_{x_\a, x_w}(z) - \Psi_{x_\a, x_w}(x) = \Psi_{x, x^w}(z).
	\earyst
	We can then deduce \eqref{eq tau n circ v f n circ v} from the H\"older bounds of the temporal functions.
\end{proof}

\subsection{Equilibrium measure} \label{subsec eqmes}

\begin{defi} \label{def gibss}
	Let $F$ be a real-valued H\"older function on $M$. The {\em pressure} of $F$ is defined by 
	\aryst
	Pr_g(F) = \sup_{\nu} [\int F  d\nu + h_\nu(g^1)].
	\earyst
	Here $\nu$ ranges over all the $g$-invariant probability measures and $h_\nu$ denotes the measure entropy with respect to $\nu$. It is known that for any $g \in \cU^2$, there is a unique $g$-invariant measure, denoted by $\nu_{g, F}$, which realizes the above supremum.
	We say that $\nu_{g, F}$ is the {\em equilibrium measure} for $g$ with {\em potential} $F$.
\end{defi}

We fix a H\"older potential $F$. We let
$\theta$ be sufficiently close to $0$ so that $F \in C^\theta(M)$.

Let $\nu_\Pi$ be the probability measure on $\Pi$ that is induced from $\nu_{g, F}$ via the flow (here we suppress the dependence of $\nu_\Pi$ on $g, F$). 
 In the following, we let $\nu_U$ denote the projection of $\nu_\Pi$ to $U$ via $\pi_U$. It is known that $\sigma$ preserves $\nu_U$. We denote by  $L^1(U, d\nu_U)$, resp. $L^2(U, d\nu_U)$, the space of  complex-valued functions on $U$ which are integrable, resp. square-integrable, with respect to $\nu_U$. 

To state various proofs in later sections in a succinct way, it is convenient to introduce the probabilistic formalism as follows. We view $(U, \cB_U, \nu_U)$ as a probability space where $\cB_U$ denotes the Borel $\sigma$-algebra of $U$. Given a measurable space $(Z, \cB_Z)$, a measurable map $\varphi: U \to Z$ is viewed as a $Z$-{\it valued random variable}.
For any function $\varphi \in L^1(U, d\nu_U)$, we denote by $\mathbb E(\varphi)$ the {\it expectation} of $\varphi$. In other words, we have
\aryst
\mathbb E(\varphi) = \int_U \varphi(x) d\nu_U(x).
\earyst
We denote the $\nu_U$-probability of a $\cB_U$-measurable set $A$ by $\PP(A)$, e.g., $\PP(A) = \mathbb E({\rm 1}_{A})$. 
The {\it conditional expectation} $\mathbb E(\varphi|A)$ is defined by $\mathbb E(\varphi|A) = \PP(A)^{-1}\mathbb E(\varphi {\rm 1}_A)$ whenever  
$\PP(A) > 0$ and $\varphi \in L^1(U, d\nu_U)$. 
The {\em conditional probability} $\PP(B \mid A) = \mathbb E(1_B \mid A)$ whenever  $\PP(A) > 0$.
In general, for any $\sigma$-subalgebra $\cB_0 \subset \cB_U$, we may define $\mathbb E(\varphi \mid \cB_0)$ and $\mathbb E(B \mid \cB_0)$ at $\nu_U$-almost every point.

%\clr We define function $\tilde f \in C^\theta(\Pi, \R)$ by
%\aryst
%\tilde f(x) = \int_0^{\tau(x)} F(x,s) ds - Pr_g(F)\tau(x).
%\earyst
%By letting $\theta$ be smaller  if necessary while keeping $F$ fixed, 
%we may choose functions $h \in C^\theta(\Pi, \R)$ and $\hat f \in C^\theta(U, \R)$ such that
%\aryst
%\norm{\hat f}_\theta &\leq& C', \\
%\tilde f(x) &=& \hat f(x) + h(x) - h(\sigma(x)), \quad \forall x \in \Pi.
%\earyst
%\clb
Moreover, by letting $\theta$ be smaller  if necessary while keeping $F$ fixed, and by \cite[Proposition 1]{Dol} as well as the remark below it, there exists $\hat f \in C^\theta(U, \R)$ with $\norm{\hat f}_\theta \leq C'$
and
$\nu_U$ satisfies the following the {\it Gibbs property}:
give an integer $n \geq 1$ and a $U$-valued random variable \clb $Y$ on $(U, {\cal B}_U)$ with distribution $\nu_U$, we set $Z = \sigma^n(Y)$. Then for any $y, z \in U$, the conditional probability $\PP(Y = y \mid Z = z)$ exists and satisfies
	\aryst
\PP(Y = y \mid Z = z) = \begin{cases}
	0, & z \neq \sigma^{n}(y), \\
	e^{\hat f_{n}(y)}, & z = \sigma^{n}(y).
\end{cases}
\earyst
Since $U$ is of dimension $1$, it is direct to deduce from the above that $\nu_U$ satisfies the {\it doubling property} (or {\it Federer property}): there exists  $c_{g,F} > 0$ such that for any $x \in U$ and any $\varrho  > 0$, we have
\aryst
\nu_U(W^u_g(x, \varrho/2) \cap U)  > c_{g,F} \nu_U(W^u_g(x, \varrho) \cap U).
\earyst

We introduce the following notion which divides the proof into two cases.
\begin{defi}
	We say that $\nu_{g,F}$ is {\em non-expanding} if
	\aryst
	\int \rdiv V_g d\nu_{g, F} \leq 0.
	\earyst
	
 We say that $g \in \cU^2$ is in {\em Class $I_F$} if $g$ is in Class I and $\nu_{g, F}$ is non-expanding.
\end{defi}
\begin{rema} \label{rem only class i-1}
Notice that for any $F \in C^\theta(M, \R)$, if $g$ is in Class $I$ but not in Class $I_F$, then the reversed flow $g^{-1}$ is in Class $I_F$.  Thus the validity of Theorem \ref{main thm} for all $g$ in Class $I_F$ implies the validity of Theorem \ref{main thm} for all $g$ in Class I. Then by Remark  \ref{rema only class i-ii}, it suffices to prove  Theorem \ref{main thm} for all $g \in \cU^\infty$ in Class $I_F$ and {\it II}.
\end{rema}

We have the following fractional moment estimate.
\begin{lemma}\label{lem subexponential growth on average}
	Let $g \in \cU^2$ and let $F \in C^\theta(M, \R)$. If $\nu_{g, F}$ is  non-expanding,
	then for any $\kappa_0 > 0$, there exist an integer $n_0 > 0$ and some $\gamma_0 > 0$ such that for any $n > n_0$  we have
	\aryst
	\int \det(D g^{\tau_n})^{\gamma_0} d\nu_{U} &<& e^{n\kappa_0 \gamma_0}.
	\earyst	
\end{lemma}
\begin{proof}
	We defer the proof to Appendix \ref{ap a}.
\end{proof}

\subsection{Complex Ruelle-Perron-Frobenius operator} \label{subsec complex rpf}
By the reduction in \cite{Dol}, it is enough to study the so-called {\it complex Ruelle-Perron-Frobenius operator} $\widehat \cL_{a,b}$ for any  $a + ib \in \C$ with $|a|$ sufficiently small and $|b|$ sufficiently large.

Fix a potential $F \in C^\theta(M)$.
For any complex number $a + ib \in \C$,
we consider the complex Ruelle-Perron-Frobenius operator $\widehat\cL_{a,b} : C^{\theta}(U, \C) \to C^{\theta}(U, \C)$ defined by
\aryst
\widehat\cL_{a,b}\varphi(x) &=& \sum_{\bar x \in \sigma^{-1}(x)} \exp(\hat f(\bar x) + (a + ib) \tau	(\bar x)) \varphi(\bar x),
\earyst
where $\hat f$ is related to $F$ as in Section \ref{subsec eqmes}.

From the Gibbs property, we see that 
\aryst
\widehat \cL_{0,0} 1 = 1.
\earyst
From the assumption that $g$ is transitive, $\widehat \cL_{0,0}$ has an isolated simple eigenvalue $1$, and hence its perturbation $\widehat \cL_{a,0}$, for all $a \in \R$ close to $0$, has a unique, isolated simple real eigenvalue $E_{a}$ close to $1$.
Let $\rho_{a}$ be the normalized (i.e., $\int \rho_a d\nu_U = 1$) eigenfunction for $\widehat \cL_{a,0}$ at $1$. It is clear that $\rho_a$ is real-valued.
Moreover by letting $a$ be close to $0$, we can ensure that $\rho_a > 1/2$ everywhere. 

The {\it normalised complex Ruelle-Perron-Frobenius operator}  $\cL_{a,b}$ for $a,b \in \R$ with $|a|$ sufficiently small is defined as follows:
\aryst
\cL_{a,b} \varphi(x)
&=& (E_a \rho_a)^{-1} \widehat \cL_{a,b}(\rho_a \varphi)(x) \\
&=& \sum_{\bar x \in \sigma^{-1}(x)} \exp(f^{(a)}(\bar x) + ib \tau	(\bar x)) \varphi(\bar x)
\earyst
where we set
\aryst \label{def of fa}
f^{(a)} = \hat f + a \tau + \log \rho_a - \log \rho_a \circ \sigma   - \log  E_a.
\earyst
We can see that $\cL_{a,0} 1 = 1$.

\subsection{Smoothing}\label{sec smoothing}

As we will see in Proposition \ref{prop decay near axis} that we are going to estimate ${\cL}_{a,b}^n$ for $n < C' \log |b|$. Hence it is harmless to consider a small perturbation of the operator which only create a negligible error up to time $C' \log |b|$. We indeed consider such a perturbation to make the real part of the coefficients of $\cL_{a,b}$ differentiable functions. This will make the argument in Section \ref{sec proofofpropmaincriterion}-\ref{sec proof lem dim cancellation} simpler.

We follow the smoothing procedure in \cite{Dol}.
Let $\delta_1 \in (0,1)$ be a small constant depending only on $g,F$ which will be determined in due course.
We choose $f_{(b)}, \tau_{(b)} \in C^1(U, \R)$, obtained from $\hat f, \tau$ respectively by making convolution with an appropriate $C^1$ mollifier function with support of size $|b|^{-\delta_1/2}$. It is direct to verify that
\ary 
\ \ \ \ \ \norm{\hat f - f_{(b)}}_{\theta /2 } \leq C' |b|^{-\delta_1 \theta/4}, & \norm{f_{(b)}}_{\theta} < C', & \norm{f_{(b)}}_{C^1} \leq C' |b|^{\delta_1},  \label{eq smoothing1} \\
\ \ \ \ \ \norm{\tau - \tau_{(b)}}_{\theta /2} \leq C' |b|^{-\delta_1 \theta/4}, &  \norm{\tau_{(b)}}_{\theta} < C', & \norm{\tau_{(b)}}_{C^1} \leq C' |b|^{\delta_1}. \label{eq smoothing2}
\eary

For any $a,b,B \in \R$ with $B \neq 0$ we set
\aryst
\cL'_{a,b,B}\varphi(x) = \sum_{\bar x \in \sigma^{-1}(x)} \exp(f_{(B)}(\bar x) + a\tau_{(B)}(\bar x) +ib\tau(\bar x)) \varphi(\bar x).
\earyst

In the rest of this section, we always let $a \in \R$ be a constant with $|a|$ sufficiently small, and let $b \in \R$ be a constant with $|b|$ sufficiently large, both depending only on $g$. Then there is a unique, isolated simple real eigenvalue $E_{a,b}$ of $\cL'_{a,0, b}$ close to $1$. 
As before, we let $\rho_{a,b}$ be the normalized eigenfunction for $\cL'_{a,0,b}$ at $E_{a,b}$. 
Then for all $b$ with $|b|$ sufficiently large, $\rho_{a,b}$ is a positive function with a lower bound independent of both $a$ and $b$, e.g., we may assume $\rho_{a,b} > 1/3$. 
We set
\ary 
f^{(a,b)} &=&  f_{(b)} + a \tau_{(b)} + \log \rho_{a,b} - \log \rho_{a,b} \circ \sigma   -   \log E_{a,b},  \label{eqfabdef} \\
\cM_{a,b} \varphi(x) &=& \sum_{\bar x \in \sigma^{-1}(x)} \exp(f^{(a,b)}(\bar x) ) \varphi(\bar x). \label{eqmabdef}
\eary
Clearly, we have
\ary \label{eq hatftofab}
\norm{\hat  f - f^{(a,b)}}_{C^0} \to 0 \quad \mbox{ as } \quad |a| \to 0, |b| \to \infty.
\eary
Moreover, we have $\cM_{a,b}1 = 1$.

We claim that
for all $b$ with $|b|$ sufficiently large
\ary \label{eq rhoarhoabeaeabclose}
\norm{\rho_a - \rho_{a,b}}_{C^0}, |E_{a} - E_{a,b}| < C'|b|^{-\delta_1 \theta/4}.
\eary
Indeed, by \eqref{eq smoothing1} and \eqref{eq smoothing2},  for any $\varphi \in C^{\theta/2}(U, \C)$ we have
\ary
\label{eq c0different1} \norm{\widehat \cL_{a,0}\varphi - \cL'_{a,0, b}\varphi}_{\theta/2}
&\leq& C' \max( \norm{\hat f - f_{(b)}}_{\theta/2}, \norm{\tau - \tau_{(b)}}_{\theta/2} )\norm{\varphi}_{\theta/2}  \\
\nonumber &\leq& C' |b|^{-\delta_1 \theta/4}\norm{\varphi}_{\theta/2}. \nonumber 
\eary
Then the claim follows from the standard perturbation theory (see \cite{Ka}) by considering $\widehat \cL_{a,0}$ and $\cL'_{a,0, b}$ as bounded linear operators in $C^{\theta/2}(U, \C)$.
Moreover  for any $\varphi \in C^\theta(U, \C)$ we have
\ary
\norm{\cL'_{a,0, b}\varphi}_{\theta} &\leq& C'\norm{\varphi}_{\theta}; 
\eary
and for any $\varphi \in C^1(U, \C)$ and any integer $n \geq 1$ we have
\ary
\norm{(\cL'_{a,0,b})^n\varphi}_{C^1} &\leq&  C' e^{-n\chi_0}\norm{\varphi}_{C^1} + C'|b|^{\delta_1}\norm{\varphi}_{C^0}.
\eary
Thus for all $b$ with $|b|$ sufficiently large we have
\ary \label{eq rhoabholderc1}
\norm{\rho_{a,b}}_\theta  < C', \quad \norm{\rho_{a,b}}_{C^1}  < C'|b|^{\delta_1}.
\eary
\clb

We set
\ary \label{eq deflab}
\widetilde \cL_{a,b} \varphi(x) &=& (E_{a,b}\rho_{a,b})^{-1} \cL'_{a,b,b}( \rho_{a,b} \varphi)(x) \\
&=& \sum_{\bar x \in \sigma^{-1}(x)} \exp(f^{(a,b)}(\bar x) + ib \tau(\bar x)) \varphi(\bar x)  \nonumber
\eary
By definition and \eqref{eq rhoarhoabeaeabclose} we have $f^{(a,b)} \in C^1(U, \R)$ and
\aryst
\norm{f^{(a,b)} - f^{(a)}}_{C^0} < C'|b|^{-\delta_1 \theta/4}.
\earyst
Then similar to \eqref{eq c0different1}, we also have
 for any $\varphi \in C^0(U, \C)$ that
\ary \label{eq c0differenceLandtildeL}
\norm{\cL_{a,b}\varphi - \widetilde\cL_{a,b}\varphi}_{C^0} \leq C' |b|^{-\delta_1 \theta/4}\norm{\varphi}_{C^0}.
\eary
We can also deduce from \eqref{eq smoothing1}, \eqref{eq smoothing2} and \eqref{eq rhoabholderc1} the above the following. We omit the proof.
\begin{lemma} \label{lem f a b}
	For all $a,b \in \R$ with $|a|$ sufficiently small and with $|b|$ sufficiently large, we have
	\aryst
	\|f^{(a,b)}\|_\theta < C', \quad
	\norm{D f^{(a,b)}}_{C^0} < C' |b|^{\delta_1}.
	\earyst
\end{lemma}

\subsection{A criterion for exponential decay}

The rest of the paper is dedicated to the following cruicial step.
For any function $u \in C^\theta(U)$, we set
\aryst
\norm{u}_{\theta,b} = \max(\norm{u}_{C^0}, |b|^{-1}|u|_\theta).
\earyst
\begin{prop}[Dolgoyat's estimate]\label{prop decay near axis}
	Given a function $F \in C^\theta(M, \R)$ for some $\theta > 0$, and $g \in \cU^\infty$ in Class $I_F$ or {\it II} such that $E^s$ and $E^u$ are not jointly integrable, there exist $ \kappa, C  > 0$ such that for all $a \in \R$ with $|a|$ sufficiently small, for all $b \in \R$ with $|b|$ sufficiently large, for all $u \in C^\theta(U)$, for all $n > C \ln |b|$, we have
	\ary \label{eq decay tol2fromtheta}
	\norm{\cL_{a,b}^n u}_{L^2(U, d\nu_U)} < |b|^{-\kappa} \norm{u}_{\theta, b}.
	\eary
\end{prop}

We first give some ideas on the proof of Proposition \ref{prop decay near axis}.
Given a H\"older function $u$, we will introduce a sequence of functions which control the modulus, as well as the regularity, of functions $\cL_{a,b}^n u$, $n \geq 1$ at a collection of {\it appropriate scales}. The main subtlety in our case is in the choice of these scales (in our case, these scales may vary from point to point):

\enmt
\item on one hand, we need to know that at the scales we choose, the temporal functions, after an appropriate normalization, has uniformly bounded H\"older norms;
\item on the other hand, we need to show that at sufficiently many places the normalized temporal function has a definite amount of oscillation.
\eenmt

Item (1) above is used to control the H\"older norm of $\cL_{a,b}^{n+1} u$ in terms of the uniform norm of $\cL_{a,b}^{n} u$. Item (2) can be used to show the decay of the uniform norm of $\cL_{a,b}^{n}u$. Indeed, one can apply Dolgopyat's argument in \cite[Section 7]{Dol} to get certain cancellation using the inductive information about the H\"older norm of $\cL_{a,b}^n u$, combined with the oscillation provided by item (2). Moreover, to get the good $L^2$ bound in Proposition \ref{prop decay near axis}, we need to know that the locus where we get cancellation is sufficiently rich: most of the points will run into this set at a definite frequency.

In the rest of the paper, for a given function $\Lambda: U \to \R_+$, we denote for every $\a \in \cA$ and every $x \in U_\a$ that
\ary \label{def jx}
J^\Lambda_x = \{ s\in (-1,1) \mid \Phi^u_x(\Lambda(x)^{-1}s) \in U_\a \}.
\eary
By definition, $J^\Lambda_x$ is a sub-interval of $(-1,1)$,  and it contains either $(-1,0]$ or $[0,1)$ when $\inf\Lambda$ is sufficiently large.
To simplify notation, we assign a map $\Phi^\Lambda_x : J^\Lambda_x \to U$ to each $x \in U$ by setting
\ary \label{def wx}
\Phi^\Lambda_x(z) = \Phi^u_x(\Lambda(x)^{-1}z), \quad \forall z \in J^\Lambda_x.
\eary

We introduce the following definitions relative to a given $g \in \cU^2$.
We first introduce a function that describes the scales we work on. Given $a,b$ as in Proposition \ref{prop decay near axis}, we will suppose $\epsilon = |b|^{-1}$ so that $\epsilon>0$ is a small number and we will consider the action of the operator ${\cal L} = {\cal L}_{a,b}$ at \lq\lq scale\rq\rq  $(\Lambda^\epsilon)^{-1}$.

\begin{defi}\label{def stableandtame}
	We say that a sequence of functions $\{ \Lambda^\epsilon: U \to \R_+ \}_{\epsilon > 0}$ in $L^\infty(U)$ is 
	\enmt
	\item[$\bullet$]
	{\bf stable} if there exist $n, \kappa > 0$ such that for all sufficiently small $\epsilon$, we have
	\ary
	\Lambda^\epsilon(x)  > \epsilon^{-\kappa}, \quad \forall x \in U,
	\eary
	and for any integer $m \geq n$, for any $v \in \sigma^{-m}$,
	\ary \label{eq stable function}
%		\norm{Dg^{\tau_{m}(x)}|_{E^u(x)}}^{-1} \Lambda^\epsilon(\sigma^m(x))^{-1}
%	&<& e^{-m \kappa} \Lambda^\epsilon(x)^{-1}, \quad 	\forall x \in U.
	\norm{Dg^{\tau_{m}(v(x))}|_{E^u(x)}}^{-1} \Lambda^\epsilon(x)^{-1}
	&<& e^{-m \kappa} \Lambda^\epsilon(v(x))^{-1}, \quad 	\forall x \in Dom(v).
	\eary
	In this case, we also say that $\{ \Lambda^\epsilon\}_{\epsilon > 0}$ is $(n, \kappa)$-stable.

	\item[$\bullet$]  $n$-{\bf adapted to} some subset $\Omega \subset U$ for some integer $n \geq 1$ if there is a constant $C > 0$ such that for all sufficiently small $\epsilon$, for any $x \in \Omega$, for any $v \in \sigma^{-n}_x$,  
	for any $y \in U$ such that $y \in W^u_g(v(x), 4\Lambda(y)^{-1})$, we have
	\aryst
	\Lambda(x) < C \Lambda(y).
	\earyst 
	In this case, we also say that $\{ \Lambda^\epsilon\}_{\epsilon > 0}$ is $(n, C)$-adapted to $\Omega$.
	\eenmt
\end{defi}

We introduce the following notion to describe the subset on which we can show cancellations. 

\begin{defi}\label{def recu}
	Given a $\sigma$-invariant measure $\nu$ on $U$ and an integer $n \geq 1$,
	we say that a subset $\Omega \subset U$ is $n$-{\bf recurrent} with respect to $\nu$ if there exist $C, \kappa > 0$ such that  for any integer $m > C$ we have
	\aryst
	\nu( \{  x \in U \mid  |\{ 1 \leq j \leq m \mid \sigma^{jn}(x) \in \Omega \}| < \kappa m  \} ) < e^{-m \kappa}.
	\earyst
	In this case, we also say that $\Omega$ is $(n, C, \kappa)$-recurrent.
\end{defi}

To illustrate the meaning of the above definitions, let us denote by $J_x$ an segment centered at $x$ with radius $\Lambda^\epsilon(x)^{-1}$ for every $x \in U$.
The notion of stable says that as $\epsilon$ decreases, these segments shrink; and for any inverse branch $v \in \sigma^{-n}_x$, $v(I_x)$ is exponentially small compared to $I_{v(x)}$. On the other hand, the notion of $n$-adaptedness implies that $v(I_x)$ is not too small compared to $I_{v(x)}$ if $x$ belongs to $\Omega$. 
The meaning of $n$-recurrence should be rather intuitive.  

Finally, we introduce the non-integrability condition which is responsible for the exponential mixing.

\begin{defi} \label{def tameanduni}
	We say that a sequence of functions $\{ \Lambda^\epsilon: U \to \R_+ \}_{\epsilon > 0}$ is
	{\bf tame} if there exist $C,  \kappa > 0$ such that for all sufficiently $\epsilon > 0$,  
	for every $x \in U$, for every $y \in (-1,1)$,
	there exists $R \in C^\theta(J^{\Lambda^\epsilon}_x,\R )$ such that
	\aryst
	\|R\|_\theta &\leq& C|y|^\kappa, \\
	|\epsilon^{-1}\Psi_{x,\Phi^s_x(y)}(\Phi^{\Lambda^\epsilon}_x(s))  -  R(s)| &<& \epsilon^{\kappa}, \quad \forall s \in J^{\Lambda^\epsilon}_x.
	\earyst
	In this case, we also say that $\{ \Lambda^\epsilon\}_{\epsilon > 0}$ is $(C, \kappa)$-tame.
	
	Given $C > 0$ and a subset $\Omega \subset U$, we say that $C$-{\rm UNI} (short for {\bf uniform non-integrability}) holds on $\Omega$ at scales $\{ \Lambda^\epsilon: U \to \R_+ \}_{\epsilon > 0}$ if there exists $\kappa > 0$ such that for 
	every sufficiently small $\epsilon > 0$,
	 for every $x \in U$ with $W^u_g(x, C\Lambda^\epsilon(x)^{-1}) \cap \Omega \neq \emptyset$, there exists $\bar y \in (-\varrho_2, \varrho_2)$ 
	such that for any $y \in (\bar y - \kappa, \bar y + \kappa)$, for any $\omega \in \R/2\pi \Z$, for $J_0 = [0,1)$ or $(-1,0]$, there is a sub-interval $J_1 \subset J_0$ with $|J_1| > \kappa$ such that
	\aryst
	\inf_{s \in J_1} \|\epsilon^{-1}\Psi_{x,y}( \Phi^{\Lambda^\epsilon}_x(s))  - \omega\|_{\R / 2\pi \Z} > \kappa.
	\earyst
		In this case, we say that $(C, \kappa)$-{\rm UNI} holds on $\Omega$ at scales $\{ \Lambda^\epsilon\}_{\epsilon > 0}$.
\end{defi}

Definition \ref{def tameanduni} is about the \lq\lq normalized\rq\rq temporal function $\epsilon^{-1}\Psi_{x,y}( \Phi^{\Lambda^\epsilon}_x(s))$. 
The notion of tame says that these functions have controlled regularity. The notion of {\rm UNI} says that these functions have lots of oscillation near $\Omega$.

The main result of this section is the following. We show that Dolgopyat's estimate (Proposition \ref{prop decay near axis}) follows from some geometric condition on the temporal functions.  We divide the proof into the following two propositions.

\begin{prop}\label{prop main 1}
	Given a potential function $F \in C^\theta(M, \R)$ for some $\theta > 0$, an Anosov flow $g \in \cU^\infty$ in Class $I_F$ or {\it II} such that $E^s$ and $E^u$ are not jointly integrable.
	Then for any $C_1 > 1$,
	for any sufficiently large integer $n_1 > 0$,  there exist 
	\enmt
	\item a subset $\Omega \subset U$ which is $n_1$-recurrent  with respect to $\nu_U$; 
	\item a stable, tame sequence of functions $\{ \Lambda^\epsilon: U \to \R_+ \}_{\epsilon > 0}$ that is $n_1$-adapted to $\Omega$
	\eenmt
	such that $C_1$-{\rm UNI} holds on $\Omega$ at scales $\{ \Lambda^\epsilon \}_{\epsilon > 0}$.
\end{prop}

\begin{prop}\label{prop main 2}
		Given a potential function $F \in C^\theta(M, \R)$ for some $\theta > 0$, an Anosov flow $g \in \cU^\infty$, there exists $C_1 > 1$ such that if the conclusion of Proposition \ref{prop main 1} is satisfied for $C_1$ and all sufficiently large $n_1$,
	then there exist $\kappa, C  > 0$ such that for all $a \in \R$ with $|a|$ sufficiently small, for all $b \in \R$ with $|b|$ sufficiently large, for all $u \in C^\theta(U)$, for all $n > C \ln |b|$, we have
	\ary \label{eq decay tol2fromtheta}
	\norm{\cL_{a,b}^n u}_{L^2(U, d\nu_U)} < |b|^{-\kappa} \norm{u}_{\theta, b}.
	\eary
\end{prop}

Proposition \ref{prop main 1} is proved in Section \ref{sec construct partition} to \ref{sec checking the conditions}, and
 Proposition \ref{prop main 2} is proved in Section \ref{sec proofofpropmaincriterion} to \ref{sec proof of lemma prop exp rec in vol exp}.
 
 \begin{proof}[Proof of Proposition \ref{prop decay near axis} and Theorem \ref{main thm}]
 	Proposition \ref{prop decay near axis} follows immediately from Proposition \ref{prop main 1} and Proposition \ref{prop main 2}.
 	It remains to give the proof of Theorem \ref{main thm}. 	By Remark \ref{rem only class i-1}, it suffices to prove Theorem \ref{main thm} for $g$ in Class $I_F$ or {\it II}.
 	By Theorem \ref{thm mixing equals nonintegrable}, we may apply Proposition \ref{prop decay near axis} to any topologically mixing $g \in \cU^\infty$ in Class $I_F$ or {\it II}.
 	By \cite{Dol}, any Anosov flow $g$ satisfying the conclusion of
 	Proposition \ref{prop decay near axis} is exponentially mixing with respect to $\nu_{g, F}$ and any H\"older test functions.
 \end{proof}

\begin{proof}[Proof of Theorem \ref{thm primeorbit}]
	By reversing the flow if necessary, we may assume without loss of generality that $g$ is in Class $I_0$ (with potential $F=0$) or {\em II}.
	Then we can follow the argument in \cite{PS} to deduce Theorem \ref{thm primeorbit}. The only modifications are: 1. we replace the norm $\|\cdot\|_{1,t}$ in \cite{PS} by norm $\|\cdot\|_{\theta, t}$; 2. we replace \cite[Proposition 4]{PS} by Proposition \ref{prop decay near axis}.
\end{proof}

\section{Temporal function and non-integrability}\label{sec Temporal function and non-integrability}

In this section, we give the main technical estimates of this paper. This will be used in Section \ref{sec proof of thm genericity of uni} to \ref{sec Uniform non-integrability on uniform set} to show Proposition \ref{prop main 1}.
Our main observation is that
temporal function can be approximated by certain finite dimensional subspace of functions which we call \lq\lq templates\rq\rq. This idea first appeared in \cite{Tsu} on  three dimensional volume preserving Anosov flows.

\subsection{Construction of charts}\label{sec constructionofcharts}
\newcommand{\rUnit}{1}
\newcommand{\rIn}{2}
\newcommand{\rOut}{\varrho_0}

In this subsection, we let $r \geq 2$ and let $g \in \cU^r$. Recall that  $\varrho_*$ and $\upsilon_*$ are given in the definition of non-stationary normal coordinates in Section \ref{sec basicpropertiesofanosovflow}. 
We will choose constants $\rOut, \rhone$ such that $\varrho_* > \rOut >  \rhone > 1$. we may assume that $\rhone$ is sufficiently large depending only on $g$, at the cost of letting $\varrho_*$ be sufficiently large and let $\upsilon_*$ be sufficiently small.

We will define below a family of charts, under which we have, among other things, that the center-unstable foliation $W^{cu}_g$ is \lq\lq straightened\rq\rq. 
\begin{defi}\label{def normalcharts}
A family of $C^{r}$ embeddings $\{ \iota_x: (-\rOut,\rOut)^3 \to M \}_{x \in M}$ is called a {\em normal coordinate system} if the $C^r$ norm of $\iota_x$ are  bounded uniformly  over all $x \in M$; and for every $x \in M$, for every $(z,y, t) \in (-\rOut,\rOut)^3$, we have 
 \aryst
 \iota_x(0,0,0) = x, \ \iota_x(z,0,0) = \Phi^u_x(z), \  \iota_x(0,y,0) = \Phi^s_x(y), \ \iota_x(z,y,t) = g^t(\iota_x(z,y,0)). 
 \earyst
Moreover, there is an integer $K > 0$ such that for every $x \in M$, the map
\aryst
g_x = \iota_{g^1(x)}^{-1} g^1 \iota_x : (-2\rhone, 2\rhone)^3 \to (-\varrho_0, \varrho_0)^3
\earyst
is a well-defined $C^r$ embedding of the form
\ary  \label{eq gx formula}
\ \ \ \ \ \ \ \ \ \ \ \ \ \ \ \ \ g_x(z,y,t) = (f_x(z,y), t + \psi_x(z,y))
\eary
where $f_x(z,y) = (f_{x,1}(z,y), f_{x,2}(z,y))$ is a $C^r$ embedding mapping $(-2\rhone, 2\rhone)^2$ into $(-\rOut,\rOut)^2$, and $\psi_x$ is a $C^r$ function on $(-2\rhone, 2\rhone)^2$ such that
\ary  \label{eq straignten 1}
\partial_y f_{x,2}(\cdot,0) \equiv \mu_x^{-1}, \quad \partial_z f_{x,1}(0,\cdot) \equiv \lambda_x^{-1}, \\ \label{eq straignten 2}
\partial_y \psi_x(\cdot,0)  \in Poly^K, \quad \partial_z\psi_x(0,\cdot) \in Poly^K.
\eary
\end{defi}

Now we proceed to the construction of a normal coordinate system. We first choose a family of charts  $\{ \bar\iota_x: (-\rOut,\rOut)^3 \to M \}_{x \in M}$ such that all the properties of a normal coordinate system except for \eqref{eq straignten 1} and \eqref{eq straignten 2} are satisfied.
%We will then modify these charts so that \eqref{eq straignten 1} and \eqref{eq straignten 2} are satisfied as well.
We define for each $x \in M$ that
$$\check g_x = (\bar\iota_{g^1(x)})^{-1} g^1 \bar\iota_x.$$
The following lemma will be proved in Appendix \ref{appendix chart normalization}.

\begin{lemma} \label{lem normalcoordinatesystem}
For any integers $r > K > \frac{\chi_*}{\chi_0}$, there is a family of $C^r$ embeddings  $\{ \check h_x: (-2\rhone, 2\rhone)^3 \to (-\varrho_0, \varrho_0)^3 \}_{x \in M}$  such that the following is true. 
Take $x \in M$ and set 
\ary
 \bar g_x =  \check h_{g^1(x)}^{-1} \check g_x \check h_x.
\eary   
Write
\aryst
 \bar g_x(z,y,t) = ( \bar f_{x,1}(z,y),  \bar f_{x,2}(z,y), t +  \bar \psi_{x}(z,y)).
\earyst
Then we have
\aryst  
\partial_y  \bar f_{x,2}(\cdot,0) \equiv \mu_x^{-1}, \quad \partial_z  \bar f_{x,1}(0,\cdot) \equiv \lambda_x^{-1}, \\ 
\partial_y \bar \psi_x(\cdot,0)  \in Poly^K, \quad \partial_z \bar \psi_x(0,\cdot) \in Poly^K.
\earyst
\end{lemma}

Let $\check h_x$ be given by Lemma \ref{lem normalcoordinatesystem}. We may define our charts by 
$$\iota_x = \bar \iota_x \check h_x.$$ 
Then all properties of a normal coordinate system, including  \eqref{eq straignten 1} and \eqref{eq straignten 2}, are satisfied by $\{ \iota_x \}_{x \in M}$. This completes the constructions of the charts.
In the rest of the paper, we let $g \in \cU^r$ where $r$ is a large integer satisfying the condition of Lemma \ref{lem normalcoordinatesystem}, and fix a normal coordinate system $\{\iota_x\}_{x \in M}$. We let $K$ be a large integer to be determined depending only on $g$ so that, among other things, \eqref{eq straignten 1} and \eqref{eq straignten 2} hold. 
%largeness of K--1%

For any $x \in M$, for any $(z,y,t) \in (-\rOut, \rOut)^3$, we set
\aryst
\tWu_x(z,y,t) = \iota_x^{-1} W^u_{g, loc}(\iota_x(z,y,t)).
\earyst
The intersection of $\iota_x^{-1}(W^{cu}_g)$ with $(-\varrho_0, \varrho_0)^2 \times \{0\}$ gives a lamination of a neighborhood of $(0,0)$ in $(-\varrho_0, \varrho_0)^2 \times \{0\}$ by $C^r$ curves in $W^{cu}_g$.
In the following lemma, whose proof is deferred to Appendix \ref{ap c}, we introduce a special parametrisation $\gamma_x$ of this lamination. %, which we call the \lq\lq u-linearization coordinate\rq\rq. 
\begin{lemma}\label{lem def of gammax}
	By letting $\upsilon_*$ in Section \ref{sec basicpropertiesofanosovflow} be sufficiently small, the following is true.
	There is a collection of maps $\{ \gamma_x = (\cX_x, \cY_x): (-\rhone,\rhone)^2 \to (-2\rhone, 2\rhone)^2 \}_{x \in M}$ and functions $\{ \phi_x: (-\rhone, \rhone)^2 \to \R \}_{x \in M}$ satisfying that
	\aryst
	\gamma_x(0,y) &=& (0,y), \quad \forall y \in (-\rhone,\rhone), \\
	\gamma_x(z,0) &=& (z,0), \quad \forall z \in (-\rhone, \rhone), \\
	\tWu_x(0, y, 0)
	&\supset& \{   ( \gamma_x(z,y), \phi_x(z,y) )   \mid z \in (-\rhone, \rhone) \}, \quad \forall y \in (-\rhone,\rhone), \\
	f_x(\gamma_x(z, y)) &=& \gamma_{g^1(x)}( \lambda_x^{-1}z,\mu_x^{-1}y), \quad \forall y \in (-\rhone, \rhone),  z \in (- |\lambda_x|\rhone, |\lambda_x|\rhone).
	\earyst
	Moreover we may suppose that for each $x \in M$,
	\enmt
	\item for each $y \in (-\rhone,\rhone)$, $\gamma_x(z,y)$ is $C^r$ in $z$, and its $C^r$ norm is bounded uniformly over all $x,y$;
	\item for every $z \in (-\rhone, \rhone)$, we have 
	\ary \label{eq straighten 3}
	\cY_x(z,y) &=& y  + O(|y|^{1+\delta}),
	\eary
	and for every integer $1 \leq l \leq r-1$ and every $y \in (-\rhone, \rhone)$, we have
	\ary \label{eq partiallcyxbound}
	\partial^l_z \cY_x(0,y) = O(|y|);
	\eary
	\item  for every $z, y \in (-\rhone, \rhone)$ we have
	\ary
	\label{eq straighten 2}
	\partial_z \cX_x(0,y) &=& 1, \\
	 \label{eq straighten 4}
	\cX_x(z,y) &=& z  + O(|y|^{\delta}).
	\eary
	\eenmt
\end{lemma}

   We define function $\widetilde \phi^s_x : (-\rhone, \rhone)^2 \to \R$ by the following relation:
   \ary \label{eq deftildephis}
   W^s_{g,loc}(
   \iota_x(\gamma_x(z,y), \tphi^s_x(z,y) ) ) \cap W^u_{g, loc}(x) \neq \emptyset.
   \eary
   %Since stable holonomy is H\"older, we know that $\widetilde \phi^s_x(z,y)$ is H\"older in $z$ for each fixed $y$.  Since $W^{cs}_g$ is a $C^1$ foliation, we know that $C^1$ in $y$ for each fixed $z$.
   We need the following linear approximation of $\widetilde \phi^s_x$.
   Given $x \in M$, we define $v^u_{g,x}, v^c_{g,x}: (-\rhone,\rhone) \to \R$ by
   \ary \label{eq defvcgvug}
   E^s(\Phi^u_x(z)) = \R (D\iota_x)_{(z,0,0)}(v^u_{g,x}(z), 1, v^c_{g,x}(z)).
   \eary
   It is known that both $v^u_{g,x}, v^c_{g,x}$ are H\"older continuous functions, and the H\"older exponents and H\"older norms are bounded uniformly depending only on $g$.

\begin{lemma} \label{lem linear approximation by stable templates}
	For every $x \in M$ and $z,y \in (-\rhone,\rhone)$, we have
	\aryst
	\widetilde \phi^s_x(z,y) = y v^c_{g,x}(z) + O(|y|^{1+\delta}).
	\earyst
\end{lemma}
\begin{proof}
	We denote by $\iota_x(z', 0, 0)$ the unique intersection between  $W^{cs}_{g,loc}(\iota_x(\gamma_x(z,y), 0))$ and $W^u_{g,loc}(x)$.
	We have
	\aryst
	\widetilde \phi^s_x(z,y) &=& \cY_x(z,y) v^c_{g,x}(z') + O(|y|^{1+\delta}) \\
	&=& y v^c_{g,x}(z') + O(|y|^{1+\delta}).
	\earyst
	By \eqref{eq straighten 4}, we have $|z-z'| = O(|y|^\delta)$. Thus we have
	\aryst
 |v^c_{g,x}(z) - v^c_{g,x}(z')| = O(|y|^\delta).
	\earyst
	This concludes the proof.
\end{proof}
   
\begin{lemma} \label{lem stable templates transformation}
	For any $x \in M$,  there exists $P \in Poly^K$ such that for any $z \in (-|\lambda_x|\varrho_1, |\lambda_x|\varrho_1)$ 
	\aryst
	 \mu_x^{-1}  v^c_{g,g^1(x)}(\lambda_x^{-1}z) =  v^c_{g,x}(z) + P(z).
	\earyst
\end{lemma}
\begin{proof}
	By definition,  for any $z \in (-|\lambda_x|\varrho_1, |\lambda_x|\varrho_1)$ and any $y \in (-\rhone, \rhone)$ we have
	\aryst
	\widetilde \phi^s_{g^1(x)}(\lambda_x^{-1}z, \mu_x^{-1}y) = \widetilde \phi^s_{x}(z,y) + \psi_x(\gamma_x(z,y)).
	\earyst
    By Lemma \ref{lem def of gammax} and Lemma \ref{lem linear approximation by stable templates}, we can differentiate both side of the above equality with respect to $y$ and evaluate at $(z,0)$ to show that
    \aryst
    \mu_x^{-1}  v^c_{g,g^1(x)}(\lambda_x^{-1}z) =  v^c_{g,x}(z) + \partial_y \psi_x(z,0).
    \earyst
	We conclude the proof by \eqref{eq straignten 2}.
\end{proof}

\subsection{Approximation of temporal function}

For any $x \in M$, and any $y,z \in (-\varrho_1, \varrho_1)$, we denote
\aryst
\Psi_{x}(z,y) &=& \Psi_{x, \Phi^s_x(y)}(z')
\earyst
where $z'$ is the unique intersection between $W^u_{g,loc}(x)$ and $W^{cs}_{g,loc}(\iota_x( \gamma_x(z,y), 0))$.

By Lemma \ref{lem shift invariance of temporal function} and Lemma \ref{lem def of gammax}, we have for any $x \in M$, any $z,y \in (-\rhone, \rhone)$ that
\ary \label{lem finvariance}
\Psi_{x}(\lambda_x z,y) = \Psi_{g^1(x)}(z, \mu_{x}^{-1}y).
\eary

The following estimate is easy to verify. We omit its proof.
\begin{lemma} \label{lem psixnorm1}
	There is $\delta_2 > 0$ such that for any $x \in M$ and any $y \in (-\rhone, \rhone)$
	\aryst
	\norm{\Psi_x(\cdot, y)}_{(-\varrho_1, \varrho_1)} < C'|y|^{\delta_2}.
	\earyst
\end{lemma}

The main result of this subsection is the following.
\begin{prop} \label{prop template approximation 0}
  There exist $\delta_0, \eta_0 \in (0,1/2)$, $C_2 > 0$  and a sequence $\{ D_n > 0 \}_{n \geq 1}$ with $\lim_{n \to \infty} D_n = 0$  such that for any $x \in M$, for any integer $n \geq 1$, there exist functions $a_1, \cdots, a_{K}: (-\rhone, \rhone) \to \R$ and $P \in Poly^K$ such that the following are true:
  \enmt
  \item 
  there is $Q \in Poly^K$ such that for any $y \in (-\varrho_1, \varrho_1)$
  \aryst
  a_1(y) = \lambda_0 \cdots \lambda_{n-1}(  -\partial_z  \phi_{x}(0,y)  + Q(y)); 
  \earyst
  \item 
  set $\tilde n = \lfloor (1-\eta_0)n \rfloor$.
   Then for any $2 \leq i \leq K$, any $y \in (-\rhone, \rhone)$, we have
  \ary \label{eq aiupperbound}
  |a_i(y)| \leq  C'|y| \sum_{m=0}^{\tilde n-1} |\mu_0 \cdots \mu_{m-1}|^{-1}|\lambda_m \cdots \lambda_{n-1}|^i;
  \eary
  \item  there is $\varkappa \in \{ \pm 1 \}$ such that 
  for any $z, y \in (-\rhone, \rhone)$, the function 
  \aryst
  \varphi_{x,y,n}(z) = (\mu_0 \cdots \mu_{n-1})^{-1}y ( v^c_{g,g^n(x)}(z) + P(z)) + \sum_{i=1}^{K} a_i(y) z^i
  \earyst
  satisfies
	\aryst
	\| \Psi_{x}(\varkappa \Lambda_n(x)^{-1}\cdot,y) - \varphi_{x,y,n} \|_{(-\varrho_1, \varrho_1)} &<& C_2(|\mu_0 \cdots \mu_{n-1}|^{-1}|y|)^{1+\delta_0} +  |\lambda_{\tilde n} \cdots \lambda_{n-1} |^{K}  \\
	&&  + |\lambda_0 \cdots \lambda_{n-1}|^2 ),  \\
	\| \varphi_{x,y,n} \|_{(-\varrho_1, \varrho_1)} &<& D_n |y|.
	\earyst
\eenmt
\end{prop}

Before presenting the proof of the above proposition, we briefly state the idea behind the proof. Imagine that for a small constant $\epsilon > 0$ we want to get an $\epsilon^{1+\delta}$-approximation of the temporal function $\Psi_{x,y}(z)$ when $z \in W^u_g(x, \epsilon)$ and $y \in W^s_{g,loc}(x)$. We can write $\Psi_{x,y}(z)$ as the difference between two parts:

(i) the flow coordinate in  $\iota_x$ of the intersection between $W^{u}_{g,loc}(y)$ and $W^{cs}_{g,loc}(z)$;

(ii) the flow coordinate in  $\iota_x$ of the intersection between $W^{cu}_{g,loc}(y)$ and $W^{s}_{g,loc}(z)$.

Notice that since the intersection in (i) is $O(\epsilon)$-close to $W^{s}_{g,loc}(x)$ (by the fact that $W^{cs}_g$ is a $C^1$ foliation), we can do some linear approximation of part (i).

To approximate part (ii), we can, for some well-chosen integer $n > 0$, iterate the flow by time $n$ to show that part (ii) can be written as a sum of two parts:

(ii-a) the flow coordinate in $\iota_{g^n(x)}$ of the intersection $W^{cu}_{g,loc}(g^n(y))$ and $W^{s}_{g,loc}(g^n(z))$;

(ii-b) the remaining part contributed by the \lq\lq coordinate change\rq\rq  from $\iota_x$ to $\iota_{g^n(x)}$, or in very rough terms, the cumulative turnings of the coordinate systems along the flow from $x$ to $g^n(x)$. 
%In this step we need to use property \eqref{lem finvariance} of the u-linearization coordinates.

 Term (ii-a) can again be approximated by linearization. To estimate part (ii-b), we notice that for any integer $0 \leq m \leq (1-\eta_0)n$, $g^m(z)$ is $\epsilon^\delta$-close to $W^s_{g,loc}(g^m(x))$; and for any integer $(1-\eta_0)n < m \leq n$, $g^m(y)$ is $\epsilon^{1-\delta}$-close to $W^u_{g,loc}(g^m(x))$, at least when $\eta_0$ is small enough. Then roughly speaking, the first few turnings are along a sequence of local stable leaves, and the rest of the turnings are along a sequence of local unstable leaves. By our choices of the normal coordinate systems, we can approximate both parts by polynomials.

\begin{proof}[Proof of Proposition \ref{prop template approximation 0}]
	Recall that $\{ f_x \}_{x \in M}$ and $\{ \psi_x \}_{x \in M}$ are defined in \eqref{eq gx formula}.
	Given any $x \in M$ and $m \in \Z$ we will denote for simplicity the following
	\ary \label{eq xmfmmumlambam}
	x_m = g^m(x), \quad f_m = (f_{m,1},f_{m,2}) = f_{x_m},   \quad \mu_m = \mu_{x_m}, \quad \lambda_m = \lambda_{x_m}.
	\eary
	Similarly, we define $\psi_m$, $\phi_m$ (see Lemma \ref{lem def of gammax}), etc.
Given an integer $n \geq 1$, $z,y \in (-\varrho_1, \varrho_1)$,
 we set
\aryst
 z_m = \lambda_m \cdots \lambda_{n-1}z, \quad y_m = (\mu_0 \cdots \mu_{m-1})^{-1}y,  \quad \forall 0 \leq m \leq n.
\earyst
By our choices of $\gamma_m, z_m, y_m$ and by Lemma \ref{lem def of gammax},  we have
\ary \label{eq fmgammamgammam+1}
f_{m}(\gamma_{m}(z_m, y_m) ) = \gamma_{m+1}(z_{m+1}, y_{m+1}), \quad \forall 0 \leq m \leq n-1.
\eary

By the invariance of unstable manifolds, we obtain
\aryst
g_{m}(\gamma_m(z_m, y_m), \phi_{m}(z_m, y_m)) = (\gamma_{m+1}(z_{m+1}, y_{m+1}), \phi_{m+1}(z_{m+1}, y_{m+1})).
\earyst
Then we have
\aryst
 \phi_{m+1}(z_{m+1}, y_{m+1}) =  \phi_{m}(z_m, y_m) + \psi_{m}(\gamma_{m}(z_m, y_m) ).
\earyst
We iterate the above equality and obtain
\ary \label{eq phin decompose}
\phi_{n}(z, y_n) = \widetilde\psi_{n}(z, y) + \phi_x(z_0,  y)
\eary
where we set
\aryst
\widetilde \psi_{n}(z, y) = \sum_{m = 0}^{n-1} \psi_m(\gamma_m(z_m, y_m)).
\earyst

Consider the decomposition $\tpsi_{n} =
\tpsi^+_{n} + \tpsi^-_{n}$ where
\ary
\tpsi^+_{n}(z,y) &=& \sum_{m = \tilde n }^{n-1}  \psi_m(\gamma_m(z_m, y_m)), \label{eq tpsiplusn} \\
\tpsi^-_{n}(z,y) &=&   \sum_{m = 0}^{\tilde n -1} \psi_m(\gamma_m(z_m, y_m)). \label{eq tpsiminusn}
\eary

We now collect some basic properties of the functions above.

\begin{lemma} \label{lem tpsi+xn app}
	Let  $\eta_0$ be sufficiently small, and let $K$ be sufficiently large. Then there exists $P_1 \in Poly^K$  such that 
	\aryst
	\tpsi^+_{n}(z,y) &=&  yP_1(z) + O(|y_{ n}|^{1+\delta}), \quad \forall z,y \in (-\varrho_1, \varrho_1).
	\earyst
\end{lemma}
\begin{proof}
	For any $0 \leq  m \leq n-1$,
	\ary \label{eq partializpsim}
	 \psi_m (z,0) = 0, \quad \forall z \in (-\rhone, \rhone).
	\eary

    By \eqref{eq straighten 3} and \eqref{eq straighten 4} in Lemma \ref{lem def of gammax}, we have
	for any $0 \leq  m \leq n$ that	
	\ary \label{sublem 2}
	|\cX_m(z_m,y_m) - z_m| &=& O(|y_m|^\delta), \\
	|\cY_m(z_m, y_m)| &=& O(|y_m|). \label{sublem 2y}
 	\eary

By taking the Taylor expansion of $\psi_m$ at $(z_m,0)$, by \eqref{eq partializpsim}, \eqref{sublem 2}, \eqref{sublem 2y} and Lemma \ref{lem def of gammax}, we have
\aryst
&&\psi_m (\cX_m(z_m, y_m), \cY_m(z_m,y_m)) \\
 &=&  \sum_{\substack{1 \leq i + j \leq K, \\ i, j \geq 0}} \partial_z^i \partial_y^j \psi_m(z_m,0)(\cX_m(z_m,y_m) - z_m)^i \cY_m(z_m, y_m)^j + O(|y_m|^{\delta K})  \nonumber \\
 &=& y_m \partial_y \psi_m (z_m,0) + O(|y_m|^{\delta K} + |y_m|^{1+\delta}) \\
 &=& (\mu_0 \cdots \mu_{m-1})^{-1} y \partial_y \psi_m(\lambda_m \cdots \lambda_{n-1} z, 0) + O(|y_m|^{\delta K} + |y_m|^{1+\delta}).
\earyst
Then by letting $\eta_0$ be sufficiently small and letting $K$ be sufficiently large,  we have
	\ary \label{eq tpsiplusexpression} \ \ \ \ \ \
	\tpsi^+_{n}(z,0) =  y  \sum_{m = \tilde n}^{n-1}  (\mu_0 \cdots \mu_{m-1})^{-1} \partial_y \psi_m(\lambda_m \cdots \lambda_{n-1}z,0)  + O(|y_n|^{1+\delta}). 
	\eary
	We conclude the proof by \eqref{eq straignten 2}.
\end{proof}

\begin{lemma} \label{lem tpsi-xn app}
	For every integer $0 < \ell \leq r-1$, we have
	\aryst
	|\partial^\ell_z\tpsi^-_{n}(z,y)| &<& C'\sum_{m=0}^{\tilde n - 1}|\lambda_m \cdots \lambda_{n-1}|^\ell, \quad \forall z,y \in (-\varrho_1, \varrho_1), \\
	|\partial^\ell_z\tpsi^-_n(0,y)| &<& C' |y|\sum_{m=0}^{\tilde n-1} |\mu_0 \cdots \mu_{m-1}|^{-1}|\lambda_m \cdots \lambda_{n-1}|^\ell, \quad \forall y \in (-\varrho_1, \varrho_1).
	\earyst
    Moreover, let $K > \frac{\chi_*}{\chi_0}$, then there exists $P_2 \in Poly^K$ such that 
	\aryst
	\partial_z\tpsi^-_{n}(0,y)
	&=&   P_2(y), \quad \forall y \in (-\varrho_1, \varrho_1).
	\earyst
\end{lemma}
\begin{proof}
	By definition, we have
	\aryst
	\partial^\ell_z\tpsi^-_{n}(z,y) &=& \sum_{m=0}^{\tilde n -1} \partial^\ell_z(\psi_m (\cX_m(z_m, y_m), \cY_m(z_m, y_m) )).
	\earyst
	To simplify the computation of the right hand side, we introduce the following notation. For any integer $a \geq 1$, an $a$-multi-index ${\bf p} = (p_1, \cdots, p_a)$ is an element of $\Z_+^a$ whose weight is defined by $|{\bf p}| = p_1 + \cdots + p_a$. We let $\emptyset$ denote the unique $0$-multi-index with weight $|\emptyset|=0$. 
	Since $\partial_z z_m = \lambda_m \cdots \lambda_{n-1}$, we have
	\aryst
	\partial^\ell_z(\psi_m (\cX_m(z_m, y_m), \cY_m(z_m, y_m) )) &=& (\lambda_m \cdots \lambda_{n-1})^\ell \sum_{ \substack{ {\bf p}, {\bf q} \\  |{\bf p}| + |{\bf q}| = \ell }  } c_{{\bf p}, {\bf q}} D_{{\bf p}, {\bf q}}(z_m, y_m) 
	\earyst
	where the sum is taken over all multi-indices ${\bf p}, {\bf q}$ satisfying $|{\bf p}| + |{\bf q}| = \ell$; $c_{{\bf p}, {\bf q}}$ is an integer depending only on ${\bf p}, {\bf q}$; and for ${\bf p} = (p_1, \cdots, p_a)$ and ${\bf q} = (q_1, \cdots, q_b)$ we set 	\aryst
	 D_{{\bf p}, {\bf q}}(z_m, y_m) &=&
	\partial^{a}_z \partial^{b}_y\psi_m( \gamma_m(z_m,y_m)) \partial^{p_1}_z\cX_m(z_m,y_m) \cdots \partial^{p_a}_z\cX_m(z_m,y_m) \\
	 &&\cdot \partial^{q_1}_z\cY_m(z_m,y_m) \cdots \partial^{q_b}_z\cY_m(z_m,y_m).
	\earyst
    Then we have
	\aryst
	|\partial^\ell_z(\psi_m( \gamma_m(z_m, y_m)))| < C'(\lambda_m \cdots \lambda_{n-1})^\ell.
	\earyst
	By Lemma \ref{lem def of gammax} \eqref{eq partiallcyxbound}, we have for any $1 \leq \ell \leq r-1$ that
	\aryst
	|\partial^\ell_z \cY_m(0, y_m)| \leq  C' |y_m|.
	\earyst
	Thus for any multi-indices ${\bf p}, {\bf q}$ with $|{\bf p}| + |{\bf q}| = \ell$ and $|{\bf q}| \geq 1$
	\ary \label{eq dboundforb>0}
	|D_{{\bf p}, {\bf q}}(0, y_m)| \leq C' |y_m|.
	\eary
	We next consider the case $|{\bf q}| = 0$. For any integer $0 \leq a \leq r-1$, we have $\partial^a_z\psi_m(0, 0) = 0$. Then we have
	\aryst
	|\partial^a_z\psi_m(\gamma_m(0, y_m))| = 	|\partial^a_z\psi_m(0, y_m)| \leq C'|y_m|,
	\earyst
	and, by Lemma \ref{lem def of gammax}(1), for any ${\bf p}$ with $|{\bf p}| = \ell$ we have
	\ary \label{eq dboundforb=0}
	|D_{{\bf p}, \emptyset}(0, y_m)| \leq C' |y_m|.
	\eary
	By \eqref{eq dboundforb>0} and \eqref{eq dboundforb=0}, we deduce that
	\aryst
	|\partial^\ell_z(\psi_m (\cX_m(z_m, y_m), \cY_m(z_m, y_m) ))|_{z=0}| \leq C' |\lambda_m \cdots \lambda_{n-1}|^\ell|y_m|.
	\earyst
	Consequently, we have
	\aryst
	|\partial^\ell_z\tpsi^-_{n}(0,y)| \leq C' |y|\sum_{m=0}^{\tilde n-1} |\mu_0 \cdots \mu_{m-1}|^{-1}|\lambda_m \cdots \lambda_{n-1}|^\ell.
	\earyst
	
	We claim that for all $0 \leq m \leq n-1$
	\ary \label{eq computepartialzcXm}
	\partial_z [\cX_m(z_m, y_m)]|_{z = 0} = \lambda_m \cdots \lambda_{n-1}.
	\eary
	Indeed, for $m=n-1$ the above is true by \eqref{eq straighten 2} in Lemma \ref{lem def of gammax}; and the case for $m < n-1$ follows from induction, $\partial_y f_{m,1}(0,\cdot) \equiv 0$, $\partial_z f_{m,1}(0, \cdot) \equiv \lambda_m$ and the formula
	\aryst
	\cX_{m+1}(z_{m+1}, y_{m+1}) = f_{m, 1}(\cX_{m}(z_{m}, y_{m}), \cY_{m}(z_{m}, y_{m})).
	\earyst

By \eqref{eq straignten 2} and \eqref{eq computepartialzcXm} we have
	\aryst
	\partial_z\tpsi^-_{n}(0,y) =  \sum_{m = 0}^{\tilde n -1}  \lambda_m \cdots \lambda_{n-1}\partial_z \psi_m(0,y_m) \in Poly^{L}(y).
	\earyst
\end{proof}

By Lemma \ref{lem tpsi-xn app}, we have
\ary \label{eq psi-xn expansion}
\widetilde \psi^{-}_{n}(z,y) &=& \sum_{i = 1}^{K} \frac{z^i}{i!}\partial^\ell_z\tpsi^-_{n}(0,y) + O((\lambda_{\tilde n} \cdots \lambda_{n-1} )^{K} ). 
\eary

By \eqref{lem finvariance}, there is $\varkappa \in \{\pm 1\}$ such that
\ary \label{eq psipushforward}
 \Psi_{x}(\varkappa \Lambda_n(x)^{-1} z, y) = \Psi_{x_n}(z, y_n).
\eary

By defintion,  \eqref{eq phin decompose}, \eqref{eq psi-xn expansion}, \eqref{eq psipushforward} and by Lemma \ref{lem linear approximation by stable templates}, we have
	\aryst
	\Psi_{x_n}(z,y_n) &=&   \tphi^s_{x_n}(z, y_n) - \phi_{n}(z, y_n) \\
	&=&   y_n  v^c_{g,x_n}(z) - \phi_{n}(z,y_n) +  O(|y_n|^{1+\delta}) \\
	&=&   y_n  v^c_{g,x_n}(z) - \phi_x(z_0,  y) - \widetilde\psi_{n}^+(z, y) - \widetilde\psi_{n}^-(z, y) +  O(|y_n|^{1+\delta}).\\
	&=&   y_n  v^c_{g,x_n}(z) - z_0 \partial_z\phi_x(0,  y) - \widetilde\psi_{n}^+(z, y) - \widetilde\psi_{n}^-(z, y) +  O(|y_n|^{1+\delta} + |z_0|^2).
	\earyst
We define for every $y \in (-\varrho_1, \varrho_1)$ that
\aryst
a_i(y) =
\begin{cases}
 - \lambda_0 \cdots \lambda_{n-1} \partial_z \phi_x(0,y) - \partial_z\tpsi^-_{n}(0,y), & i=1, 	\\
 -\frac{1}{i!}\partial^\ell_z\tpsi^-_{n}(0,y), & i=2,\cdots, K.  
\end{cases}
\earyst
We conclude the proof of  Proposition \ref{prop template approximation 0} by combining Lemma \ref{lem tpsi+xn app} and Lemma \ref{lem tpsi-xn app}.
\end{proof}

\subsection{Template function}

In this section, we introduce a collection of intrinsically defined objects called \lq\lq stable templates\rq\rq and \lq\lq unstable templates\rq\rq, and reformulate the main result in the previous section using these new objects.

\def\Tplt{          \Xi}
\def\tplt{          \xi}

\def\cGu{          \mathcal T^u}
\def\cGs{          \mathcal T^s}

For any $x \in M$, we let $\Tplt^u(x)$ denote the set of continuous sections $\tplt: W^u_{g}(x, \varrho_1) \to T^*M$ such that for all $\bar x \in W^u_{g}(x, \varrho_1)$, $\tplt(\bar x) \in T^*_{\bar x}M$ vanishes on $E^u(\bar x)$.
Let $\Tplt^u_0(x)$, resp. $\Tplt^u_1(x)$, be the subset of $\Tplt^u(x)$ that consists of $\tplt$ satisfying $\tplt(\bar x)(V_g(\bar x)) = 0$, resp. $1$, for all $\bar x \in W^u_{g}(x, \varrho_1)$. If $g \in \cU^r$ for some $r \geq 1$, we denote by $\Tplt^{u,k}$, resp. $\Tplt^{u,k}_0, \Tplt^{u,k}_1$, for $0 \leq k \leq r$ the subset of $C^k$ sections in $\Tplt^u$, resp. $\Tplt^{u}_0, \Tplt^u_1$.

For any $x \in M$, any $t > 0$,  $g^t$ induces a linear map $L^t_x: \Tplt^u(x) \to \Tplt^u(g^t(x))$ given by
\aryst
L^t_x \tplt(w) = (Dg^{-t})^*[\tplt(g^{-t}(w))], \quad \forall w \in W^u_{g, loc}(g^t(x)).
\earyst
If $g \in \cU^r$,  $L^t_x$ maps $\Tplt^{u,r-1}(x)$, $\Tplt^{u,r-1}_0(x)$ and $\Tplt^{u,r-1}_1(x)$ into $\Tplt^{u,r-1}(g^t(x))$, $\Tplt^{u,r-1}_0(g^t(x))$ and $\Tplt^{u,r-1}_1(g^t(x))$ respectively.

For each $x \in M$, we denote by $\tplt^{u,s}_x \in \Tplt^u_1(x)$ the unique section such that for every $\bar x \in W^u_{g}(x, \varrho_1)$,
\ary \label{eq def xis}
{\rm Ker} \tplt^{u,s}_x(\bar x) \supset E^s(\bar x).
\eary 

By the flow-invariance, we have
\aryst
L^t_x \tplt^{u,s}_x = \tplt^{u,s}_{g^t(x)}, \quad \forall t > 0, x \in M.
\earyst

\begin{lemma} \label{lem tpltperp}
	Assume that $r \geq 2$ and $g \in \cU^r$.
    Then there is a family of $C^{r-1}$ sections $\{ \tplt^{u,\perp}_x \in \Tplt^{u,r-1}_0(x) \}_{x \in M}$ with uniformly bounded $C^{r-1}$ norms such that 
	\aryst
	\norm{\tplt^{u, \perp}_x|_{E^s(x)}} = 1, \quad
	L^1_x \tplt^{u,\perp}_x =  \pm \mu_x \tplt^{u, \perp}_{g^1(x)}, \quad \forall  x \in M.
	\earyst 
	Moreover, for every $x \in M$,  $\tplt^{u,\perp}_x$ is uniquely determined up to a sign.
\end{lemma}

\begin{proof}
	We defer the proof to Appendix \ref{ap d}.
\end{proof}

An explicit choice of the family of sections in Lemma \ref{lem tpltperp} is given as follows.
	On $\R^3$ we use the coordinate $(z,y,t)$. Then it is natural to denote by $dy, dt : \R^3 \to \R$ the linear maps	
\ary \label{eq defdydt}
dy(a,b,c) = b, \quad dt(a,b,c) = c, \quad \forall (a,b,c) \in \R^3.
\eary
We define $\widetilde \tplt^\perp_x: W^u_g(x, \varrho_1) \to T^*M$ by 
\ary \label{eq tildexiperpxtilde0x}
\widetilde \tplt^\perp_x(\bar x) = (D\iota_x^{-1})^*_{\bar x}dy, \quad \forall \bar x \in  W^u_g(x, \varrho_1).
\eary
We will see in Appendix \ref{ap d} that $\{ \widetilde \xi^{\perp}_x \}_{x \in M}$ satisfies Lemma \ref{lem tpltperp} in place of $\{  \xi^{u, \perp}_x \}_{x \in M}$.

\begin{lemma} \label{lem choose straight sections}
	Assume that $K > \frac{\chi_*}{\chi_0}$ and $g \in \cU^r$ for some $r \geq K+2$. 
	Then there is a family of sections $\{ \tplt^{u,0}_x \in \Tplt^{u,K+1}_1(x) \}_{x \in M}$ such that the following is true.
	For any $x \in M$, any $\tplt \in \Tplt^{u}_1(x)$, we denote by $\varphi^u_{x, \tplt} \in C^0(-\rhone, \rhone)$ the unique function such that
	\aryst
	\tplt(\Phi^u_x(s)) = \tplt^{u,0}_x(\Phi^u_x(s)) + \varphi^u_{x,\tplt}(s) \tplt^{u, \perp}_x(\Phi^u_x(s)),
	\quad \forall s \in (-\rhone, \rhone).
	\earyst
	Then the following are equivalent:
	\enmt
	\item $\varphi^u_{x,\tplt} \in Poly^K$;
	\item there exist $C > 0$ and a sequence of sections $\{ \tplt_n \in \Tplt^{u,K+1}(g^{-n}(x)) \}_{n \geq 1}$ such that for all $n \geq 1$, we have $\norm{D^{K+1}\varphi^u_{g^{-n}(x), \tplt_n}}_{(-\rhone, \rhone)}  < C$ and $L^n_{g^{-n}(x)} \tplt_n = \tplt$.
	\eenmt
\end{lemma}

\begin{proof}
	We defer the proof to Appendix \ref{ap d}.
\end{proof}

An explicit choice of the family of sections in Lemma \ref{lem choose straight sections} is given as follows.
We define $\widetilde \tplt^0_x: W^u_g(x, \varrho_1) \to T^*M$ by 
\ary \label{eq tildexiperpxtilde0x2} 
\widetilde \tplt^0_x(\bar x) = (D\iota_x^{-1})^*_{\bar x}dt, \quad \forall \bar x \in  W^u_g(x, \varrho_1) 
\eary
where $dt$ is given by \eqref{eq defdydt}.
We will see in Appendix \ref{ap d} that $\{ \widetilde \tplt^0_x \}_{x \in M}$ satisfies Lemma \ref{lem choose straight sections} in place of $\{  \tplt^{u,0}_x \}_{x \in M}$.

We have the following corollary of Lemma \ref{lem choose straight sections}.
\begin{cor} \label{cor uniqueness of varphi}
	Let  $\{ \tplt'_x \in \Tplt^{u,K+1}_1(x) \}_{x \in M}$ be a family of sections, and let us denote for any $x \in M$ and any $\xi \in \Xi^u_1(x)$ the unique function $\varphi'_{x, \xi} \in C^0(-\varrho_1, \varrho_1)$ such that
	\aryst
		\tplt(\Phi^u_x(s)) = \tplt^{'}_x(\Phi^u_x(s)) + \varphi'_{x,\tplt}(s) \tplt^{u, \perp}_x(\Phi^u_x(s)),
	\quad \forall s \in (-\rhone, \rhone).
	\earyst
	Assume that Lemma \ref{lem choose straight sections} holds for $(\tplt'_x, \varphi'_{x,\xi})$ instead of $(\tplt^{u, 0}_x, \varphi^u_{x,\xi})$. Then
	for any $x \in M$  we have
	\aryst
	\varphi^u_{x, \tplt'_x} \in Poly^{K}.
	\earyst
	Consequently for any $x \in M$ and any $t > 0$,\footnote{For each $t > 0$, we can let $\xi'_x = L^t_{g^{-t}(x)}\tplt^{u, 0}_{g^{-t}(x)}$ for every $x \in M$.} we can deduce that
	\aryst 
\varphi^u_{x, L^t_{g^{-t}(x)}\tplt^{u, 0}_{g^{-t}(x)}} \in Poly^{K}.
\earyst	
\end{cor}

By enlarging $K$ if necessary, we may suppose that $K$ satisfies the condition of Lemma \ref{lem choose straight sections}. %largeness of K--2%
\begin{defi}
	For any $x \in M$, we denote
	\aryst
	\varphi^{u,s}_{x} = \varphi^u_{x, \tplt^{u, s}_x}.
	\earyst
	The set of {\em stable templates} at $x$ is defined by
	\aryst
	\cGs_x = \{ c  \varphi^{u,s}_{x}+  P \mid  c \in \R, P \in Poly^K,  P(0) = 0  \}.
	\earyst
	By Corollary \ref{cor uniqueness of varphi}, we see that for every $x \in M$,  $\cGs_x$ is independent of the choice of $\{ \tplt^{u, 0}_w \}_{w \in M}$ satisfying Lemma \ref{lem choose straight sections}.
	In a similar way,
	we define $\{ \varphi^{s,u}_x \}_{x \in M}$ and the space of {\em unstable templates} $\cGu_x$ for every $x \in M$ by reversing the flow  and switching the roles of stable and unstable subspaces.  Clearly we also have $\varphi^{u,s}_x(0) = \varphi^{s,u}_x(0) = 0$ for all $x \in M$.
\end{defi}

Stable and unstable templates are related to the approximation of temporal functions, as the following lemmas suggest.

\begin{lemma} \label{lem vcgx}
	For any $x \in M$, there exist $\varkappa \in \{\pm 1\}$ and $P \in Poly^K$ with $P(0) = 0$ such that 
	\aryst
	v^c_{g,x} = \varkappa \varphi^{u, s}_{x} + P
	\earyst
	where $v^c_{g,x}$ is given by \eqref{eq defvcgvug}.
\end{lemma}
\begin{proof}
	Define function $\widehat \tplt^s_x: (-\varrho_1, \varrho_1) \to T^*\R^3$ by  
	\aryst
	\widehat \tplt^s_x(s) = (D\iota_x)^*_{(s,0,0)} \tplt^{u,s}_x(\Phi^u_x(s)) \in T^*\R^3.
	\earyst
	By definition, the linear map $\widehat \tplt^s_x(s)$ vanishes on $\R \times \{(0,0)\}$ for any $s \in (-\varrho_1, \varrho_1)$.  
	
	By Definition \ref{def normalcharts} and \eqref{eq defvcgvug}, we have for any $s \in (-\varrho_1, \varrho_1)$ that
	\aryst
	E^s(\Phi^u_x(s)) &=&
	\R (D\iota_x)_{(s,0,0)}(v^u_{g,x}(s), 1, v^c_{g,x}(s)), \\
	V_{g}(\Phi^u_x(s)) &=& \R (D\iota_x)_{(s,0,0)}(0,0,1).
	\earyst
	Thus we have for any $s \in (-\varrho_1, \varrho_1)$ that
	\aryst
	&&\widehat \tplt^s_x(s)(v^u_{g,x}(s), 1, v^c_{g,x}(s)) = \tplt^s_x(\Phi^u_x(s))((D\iota_x)_{(s,0,0)}(v^u_{g,x}(s), 1, v^c_{g,x}(s))) = 0, \\
	&&\widehat \tplt^s_x(s)(0, 0, 1) = \tplt^s_x(\Phi^u_x(s))((D\iota_x)_{(s,0,0)}(0,0,1)) = 1.
	\earyst
	By linearity, we have
	\ary \label{eq calculatevcgx 1}
	\widehat \tplt^s_x(s) = dt - v^c_{g,x}(s) dy, \quad \forall s \in (-\rhone, \rhone).
	\eary
	
    Let $\widetilde \tplt^\perp_x, \widetilde \tplt^0_x$
    be given by \eqref{eq tildexiperpxtilde0x} and \eqref{eq tildexiperpxtilde0x2} respectively. Recall that $\{ \widetilde \xi^{\perp}_x \}_{x \in M}$ satisfies Lemma \ref{lem tpltperp} in place of $\{  \xi^{u, \perp}_x \}_{x \in M}$; and $\{ \widetilde \tplt^0_x \}_{x \in M}$ satisfies Lemma \ref{lem choose straight sections} in place of $\{  \tplt^{u,0}_x \}_{x \in M}$.
    Then by Lemma \ref{lem tpltperp} there is $\varkappa \in \{\pm 1\}$ such that
	\ary \label{eq tpltperpequalstpltuperp}
	\widetilde \tplt^\perp_x = - \varkappa  \tplt^{u,\perp}_x,
	\eary
and by Corollary \ref{cor uniqueness of varphi} we obtain
	\aryst
	P =
	\varphi^u_{x, \widetilde \tplt^0_x} \in Poly^K.
	\earyst
	Combining the above with \eqref{eq tildexiperpxtilde0x} and \eqref{eq tpltperpequalstpltuperp} we obtain
	\ary 
	\widehat \tplt^s_x(s)   \nonumber
	&=& (D\iota_x)^*_{(s,0,0)}(\tplt^{u,0}_x(\Phi^u_x(s)) + \varphi^{u,s}_x(s) \tplt^{u,\perp}_x(\Phi^u_x(s))) \nonumber \\
	&=& (D\iota_x)^*_{(s,0,0)}(\widetilde \tplt^0_x(\Phi^u_x(s)) + (\varphi^{u,s}_x(s) - P(s)) \tplt^{u,\perp}_x(\Phi^u_x(s))) \nonumber \\
	&=& dt - \varkappa (\varphi^{u,s}_x(s) - P(s)) dy. \label{eq calculatevcgx 2}
	\eary
	The lemma then follows from combining \eqref{eq calculatevcgx 1} and \eqref{eq calculatevcgx 2}.
\end{proof}

\begin{lemma}\label{lem partialzphix}
	For any $x \in M$, there exist $\varkappa \in \{ \pm 1\}$ and $P \in Poly^K$ with $P(0) = 0$ such that
	\aryst
	\partial_z \phi_x(0, \cdot)  = \varkappa \varphi^{s,u}_x + P
	\earyst
	where $\phi_x$ is given by Lemma \ref{lem def of gammax}. 
\end{lemma}
\begin{proof}  
	Notice for any $x \in M$ there exists a function $w^s_{g,x}: (-\rhone, \rhone) \to \R$ such that 
	\aryst
	(D\iota_x)_{(0, s, 0)}(1, w^s_{g,x}(s), \partial_z \phi_x(0, s)) \in E^u(\Phi^s_x(s)), \quad \forall s \in (-\rhone, \rhone)
	\earyst
	where $\phi_x$ is given by Lemma \ref{lem def of gammax}.
	The proof then follows from a similar argument as in Lemma \ref{lem vcgx}.
\end{proof}

We define function $\cD^s: M \to \R_{\geq 0}$ as follows. For each $x \in M$, we set
\aryst
\cD^s(x) = \inf_{\substack{P \in Poly^K, \\ J = [0,1) \mbox{\small{ or }} (-1,0]}} \| \varphi^{u,s}_x - P \|_{J}.
\earyst
We define $\cD^u$ analogously, using the unstable templates in place of stable templates.

%\clr
%We now show that 
% $\cT^s_x$ depends continuously on $x$. For this purpose, for any integer $d \geq 1$, 
%we denote by $\mathbb{V}_d$ the set of $d$-dimensional linear subspaces of $C^0(-\varrho_1, \varrho_1)$. For any ${\cal V} \in \mathbb{V}_d$ with a basis $v_1, \cdots, v_d \in {\cal V}$, for any sufficiently small $\eta > 0$, 
%we denote by $U_{v_1, \cdots, v_d, \eta}$ the set of $d$-dimensional linear subspaces generated by functions $u_1, \cdots, u_d \in C^0(-\varrho_1, \varrho_1)$ with $\max_{1 \leq i \leq d}\norm{v_i - u_i}_{(-\varrho_1, \varrho_1)} < \eta$ (notice that $u_1, \cdots, u_d$ are linearly independent if $\eta$ is sufficiently small). We equipe $\mathbb V_d$ with the topology generated by all such $U_{v_1, \cdots, v_d, \eta}$. 

\begin{lemma}\label{rem templatec0onx}
	The functions $\cD^s$ and $\cD^u$ are continuous.
\end{lemma} 
\begin{proof}
	By Lemma \ref{lem vcgx}, for any given normal coordinate system, 
	 we have
	\aryst
\cD^s(x) = \inf_{\substack{P \in Poly^K, \\ J = [0,1) \mbox{\small{ or }} (-1,0]}} \| v^{c}_{g,x} - P \|_{J}, \quad \forall x \in M.
	\earyst
	We claim that for any $x_0 \in M$, we can arrange it so that $v^c_{g,x}$ depends continuously on $x$ near $x_0$.
	Indeed, we can first assume that the  non-stationary normal coordinates along the stable manifolds, resp. unstable manifolds, depend continuously on $x$ near $x_0$; and then assume that the coordinate system $\{ \bar\iota_x: (-\rOut,\rOut)^3 \to M \}_{x \in M}$ below Definition \ref{def normalcharts} satisfies that $\bar\iota_x$ depends continuously on $x$ near $x_0$. By the proof of Lemma \ref{lem normalcoordinatesystem}, it is clear that $\check h_x$ depends continuously on $\bar\iota_x$. Thus $\check h_x$ depends continuously on $x$ near $x_0$. Consequently, there is a normal coordinate system $\iota_x$ depending continuously on $x$ near $x_0$. Then $v^c_{g,x}$ depends continuously on $x$ near $x_0$ as well. This concludes the proof.
\end{proof}
\clb

\begin{defi} \label{def txn}
	Given an integer $n > 0$ and $x \in M$, we denote
	\aryst
	\cT_{x,n} = \{ h_1 y \varphi^{u,s}_{g^n(x)}(z) + h_2 z \varphi^{s,u}_x(y) + P(z,y) \mid h_1, h_2 \in \R, P \in Poly^{K,K},  P(0,\cdot) = P(\cdot, 0) \equiv 0 \}.
	\earyst
\end{defi}
We notice that for any $x \in M$, any integer $n > 0$ and any $R \in \cT_{x,n}$, we have $R(\cdot, y ) \in \cT^s_{g^n(x)}$ for every $y \in (-\varrho_1, \varrho_1)$.

The following is a corollary of Proposition \ref{prop template approximation 0}.

\begin{cor}\label{prop template approximation}
For all sufficiently large $K > 1$,
there exist $\delta_0 \in (0,1/2)$ and a sequence $\{ D_n > 0 \}_{n \geq 1}$ with $\lim_{n \to \infty} D_n = 0$  such that for all sufficiently small $\epsilon > 0$, for any $x \in M$, for any integer $n \geq 1$ satisfying $\norm{Dg^n|_{E^s(x)}}, \norm{Dg^n|_{E^u(x)}}^{-1} < \epsilon$, there exist $R \in \cT_{x,n}$, $\varkappa \in \{\pm 1\}$, and functions $a_2, \cdots, a_{K}: (-\rhone, \rhone) \to \R$ satisfying \eqref{eq aiupperbound} such that for any $y \in (-\varrho_1, \varrho_1)$
\aryst
|  \Psi_{x}(\varkappa \Lambda_n(x)^{-1} z,y) - R(z, y) - \sum_{i=2}^{K} a_i(y) z^i | &<& C'((\epsilon|y|)^{1+\delta_0} +  \epsilon^2),  \quad \forall z \in (-\varrho_1, \varrho_1), \\
\| R(\cdot, y) \|_{(-\varrho_1, \varrho_1)} &<& D_n  |y|^{\delta_0}.
\earyst
\end{cor}

\begin{comment}
	
By combining Definition \ref{def uniplus}, Lemma \ref{lem vcgx} and Lemma \ref{lem partialzphix}, we obtain the following. 
\begin{lemma}
	Let $r, K \geq 2$ be integers, let $\kappa > 0$,
	and let $J = [0,1)$ or $(-1, 0]$.
The following are true.
\enmt
\item for any $g \in \cU^r_{+, K, \kappa}$, for any $x \in M$, for any  $P \in Poly^K$, we have $\norm{v^c_{g,x} + P}_J > \kappa$;
\item for any $g \in \cU^r_{-, K, \kappa}$, for any $x \in M$, for any  $P \in Poly^K$, we have $\norm{	\partial_z \phi_x(0, \cdot) + P}_J > \kappa$.
\eenmt	
\end{lemma}

\end{comment}

\section{From non-smooth distribution to non-integrability} \label{sec proof of thm genericity of uni}

\subsection{A regularity result}
\begin{defi} \label{def uniplus}
	Given $\eta > 0$ and integers $r \geq 2$, $K \geq 1$, we denote by  $\cU^r_{+, K, \eta}$ the subset of flows in $\cU^r$ verifying that for any $x \in M$, for any polynomial $P \in Poly^K$, we have $\norm{\varphi^{u,s}_{x} + P}_J > \eta$
	for $J = [0,1)$ and $J=(-1, 0]$.
	
	We set 
	\aryst
	\cU^r_{+, K} =  \bigcup_{\eta > 0} \cU^r_{+, K, \eta}.
	\earyst
	
	In a similar way, we define $\cU^r_{-, K, \eta}$ and $\cU^r_{-,K}$ using unstable templates instead of stable templates.
\end{defi}

We begin with an important property of the stable/unstable templates. 
\begin{lemma} \label{lem alltrueallfalse}
	Let $g \in \cU^r \setminus \cU^r_{+, K}$ for some integers $ K \geq 1$ and $r \geq 2$. Then for every $x \in M$, we have $\varphi^{u,s}_{x} \in Poly^K$. 	In another words, $\cU^r_{+, K}$ is the subset of flows in $\cU^r$ verifying that $0 \notin \cT^s_x$ for any $x \in M$.
	We have a similar statement for $\cU^r_{-, K}$.
\end{lemma}
\begin{proof}  
	Let $c \in (0,1]$, we denote by $X_c$ the set of $x \in M$ such that
	\aryst
	\varphi^{u,s}_x|_{(-c,0]} \ \mbox{ or } \ \varphi^{u,s}_x|_{[0,c)} \in Poly^K.
	\earyst
	Since the uniform norm of $\varphi^{u,s}_x$ on $[-1,1]$ is uniformly bounded for all $x \in M$,  we deduce that for any $c \geq 1/2$
	\aryst
		\varphi^{u,s}_x|_{(-c,0]} \ \mbox{ or } \ \varphi^{u,s}_x|_{[0,c)}  \in Poly^K_{C'}.
	\earyst
	By Lemma \ref{rem templatec0onx}, $X_1$ is a closed set.
	By Corollary \ref{cor uniqueness of varphi}, there exists $t_0  > 0$ such that for any $t  > t_0$ we have
	\ary \label{eq xonehalfandxone}
	g^t(X_{1/2}) \subset X_1.
	\eary
	By hypothesis, $X_1 \neq \emptyset$. Let $x \in X_1$. Then by Lemma \ref{lem choose straight sections}, we notice that $X_{1/2}$ contains an open set $W$ of $W^u_{g,loc}(x)$ close to $x$. By \eqref{eq xonehalfandxone}, $X_1$ contains the union of $g^t(W)$ for all $t > t_0$.
	 But since $g$ is transitive, this union is dense, and consequently we have $X_1 = M$. 
\end{proof}

We have  the following observation.
\begin{lemma}\label{lem smooth bundle}
	For any integers $K \geq 1, r \geq 3$ and $g \in \cU^r \setminus \cU^r_{+, K}$, the subspace $E^{s}$ is $C^{1+\delta}$. Similarly, for any integers $L \geq 1, r \geq 3$ and   $g \in \cU^r \setminus \cU^r_{-, L}$, the subspace $E^{u}$ is $C^{1+\delta}$.
\end{lemma}
\begin{proof}
	Let $g \in \cU^r \setminus \cU^r_{+, K}$ for some integers $K \geq 1$ and $r \geq 3$.
	By Lemma \ref{lem alltrueallfalse}, we see that for every $x \in M$, $E^s \oplus E^u$ is uniformly $C^r$ restricted to $W^u_{g}(x, \varrho_1)$. On the other hand, we know that $E^{cs}_g$ is $C^{1+\delta}$ everywhere. Consequently, $E^s$, as the intersection between $E^u \oplus E^s$ and $E^{cs}_g$, is uniformly $C^{1+\delta}$ restricted to $W^u_{g}(x, \varrho_1)$ for every $x \in M$. It is clear that $E^s$ is uniformly $C^2$ restricted to $W^{cs}_g$. Moreover, $W^{cs}_g$ and $W^u_g$ are continuous foliations with uniformly $C^r$ leaves. Then by Journ\'e's lemma in \cite{Jou}, we conclude that $E^s$ is $C^{1+\delta}$ everywhere on $M$. 
	
	The case where $g \in \cU^r \setminus \cU^r_{-, K}$ follows from a similar argument.
\end{proof}

By Lemma \ref{lem alltrueallfalse} and Lemma \ref{lem smooth bundle}, we obtain the following.
\begin{cor}\label{cor case i to iii}
	If $g \in \cU^\infty$ belongs to Class I, then for any integer $K \geq 1$, $g \in \cU^\infty_{+ , K} \cap \cU^\infty_{-, K}$.
\end{cor}
As an immediate consequence, we have the following.
\begin{cor}
	Any $g \in \cU^\infty$ in Class {\it I} is topologically mixing.
\end{cor}	
\begin{proof}
	Let $K$ be the integer given by Corollary \ref{prop template approximation}. 
	By Corollary \ref{cor case i to iii}, we have $g \in \cU^\infty_{+, K}$.
	Then by Corollary \ref{prop template approximation} for any $x \in M$ and any $y \in W^s_g(x,1)$, the function $\Psi_{x,y}$ is not identically zero in any neighborhood of $x$ in $W^u_g$. Consequently, $E^s$ and $E^u$ are not jointly integrable. By Theorem \ref{thm mixing equals nonintegrable}, $g$ is topologically mixing.
\end{proof}

\subsection{Non-integrability condition revisited}

In this subsection, we define the non-integrability condition for Anosov flows in Class {\it II}.

For any $x \in M$, any $z \in W^u_{g, loc}(x)$, we denote by $H^u_{x, z} : W^{cs}_{g,loc}(x) \to W^{cs}_g(z)$ the {unstable holonomy} map.
Since $E^u$ is $C^{1+\delta}$, we deduce that $H^u_{x, z}$  has uniformly bounded $C^{1+\delta}$ norm for all $x \in M$ and $z \in W^u_{g, loc}(x)$. 

We let $\{ \xi^{u, s}_x \}_{x \in M}$ be given by \eqref{eq def xis} and let $\{ \xi^{u, \perp}_x \}_{x \in M}$ be given by Lemma \ref{lem tpltperp}. 
%Define $\{ \xi^{ s, u}_x \}_{x \in M}$ and $\{ \xi^{s, \perp}_x \}_{x \in M}$ in a similar way.
For each $x \in M$ we let $\xi^{\cal U}_x$ denote the unique element in $\Tplt^{u}_1(x)$ satisfying
\ary \label{eq def of xiU}
DH^u_{x, \bar x}(E^s(x)) \subset 
{\rm Ker}\xi^{\cal U}_x(\bar x), \quad \forall \bar x \in W^u_g(x, \varrho_1).
\eary

\begin{defi} \label{def staletempletincaseii}
	For any $x \in M$, we define $\phi^{u,s}_x : (-1,1) \to \R$ by equation
	\ary \label{eq xiusxiUxiuperp}
	 \xi^{u,s}_x(\Phi^u_x(s)) = \xi^{\cal U}_x(\Phi^u_x(s)) + \phi^{u,s}_x(s) \xi^{u,\perp}_x(\Phi^u_x(s)), \quad \forall s \in (-1,1).
	\eary
\end{defi}

\begin{rema} \label{rem templatec0onx2}
	Similar to Lemma \ref{rem templatec0onx}, we can show that the function $\cD: M \to \R_{\geq 0}$ defined by 
	\ary \label{eq definecD}
	\cD(x) = \min(Osc_{(-1,0]} \phi^{u,s}_x,   Osc_{[0,1)}\phi^{u,s}_x  ), \quad \forall x \in M
	\eary
	is continuous.
\end{rema}

We define $V^u_{g,x}, V^c_{g,x}: (-1,1) \to \R$ for all $x \in M$ by
\ary \label{eq dhuescoordinates}
DH^u_{x, \Phi^u_x(z)}(E^s(x))
 = \R (D\iota_x)_{(z,0,0)}(V^u_{g,x}(z), 1, V^c_{g,x}(z)), \forall z \in (-1,1).
\eary
It is known that both $V^u_{g,x}, V^c_{g,x}$ are H\"older continuous functions, and the H\"older exponent and the H\"older norms are bounded uniformly depending only on $g$.

\begin{lemma} \label{lem approximation of temporal case ii}
	There is $\theta \in (0,1)$ depending only on $g$ such that for any $x \in M$, there is $\varkappa \in \{ \pm 1 \}$ such that  for any $y, z \in (-1,1)$ we have
	\aryst
	\Psi_x(z, y) &=& \varkappa y \phi^{u,s}_x(z)  + O(|y|^{1+\theta}) \\
\mbox{and} \quad 	\norm{\phi^{u,s}_x}_\theta &<& C'.
	\earyst 
\end{lemma}
\begin{proof}  
	We fix $z,y \in (-1,1)$.
	As in the proof of Lemma \ref{lem linear approximation by stable templates},
	we denote by $\iota_x(z', 0, 0)$ the unique intersection between  $W^{cs}_{g,loc}(\iota_x(\gamma_x(z,y), 0))$ and $W^u_{g,loc}(x)$. 
	By \eqref{eq straighten 4} in Lemma \ref{lem def of gammax}, we deduce that 
	\ary \label{eq zminusz'}
	|z-z'| = O( |y|^\delta).
	\eary
	
	We denote the $\iota_x$-coordinate of the intersection between $W^s_g(\iota_x(z',0,0))$ and $W^{cu}_g(\iota_x(0,y,0))$  by
	\aryst
	(\cX_x(z,y), \cY_x(z,y), t_1);
	\earyst 
	and denote the $\iota_x$-coordinate of the intersection between $W^{cs}_g(\iota_x(z',0,0))$ and $W^{u}_g(\iota_x(0,y,0))$  by
	\aryst
	(\cX_x(z,y), \cY_x(z,y), t_2).
	\earyst  
	By definition, we have
	\ary \label{eq psit1t2}
	\Psi_x(z,y) =  t_1 - t_2.
	\eary
	
	Since $W^s_g(\iota_x(z',0,0))$ is a uniformly $C^r$ curve and $v^c_{g,x}$ is uniformly H\"older, we have, by \eqref{eq straighten 3} and \eqref{eq zminusz'}, that
	\ary \label{def t1}
	t_1 &=& \cY_x(z,y)v^c_{g,x}(z') + O( |\cY_x(z,y)|^{2}) \nonumber \\
	&=& y v^c_{g,x}(z) + O( |y|^{1+\delta}).
	\eary
	Since $H^{u}_{x, \Phi^u_x(z')}$ has  uniformly bounded  $C^{1+\delta}$ norm, and $V^c_{g,x}$ is uniformly H\"older, again by \eqref{eq straighten 3} and \eqref{eq zminusz'}, we deduce that
	\ary \label{def t2}
	t_2 &=&\cY_x(z,y) V^c_{g,x}(z') + O( |\cY_x(z,y)|^{1+\delta}) \nonumber \\
	&=&  y V^c_{g,x}(z) + O( |y|^{1+\delta}).
	\eary
	Recall that
	\aryst
	E^s(\Phi^u_x(z)) =
	\R (D\iota_x)_{(z, 0, 0)}(v^u_{g,x}(z), 1, v^c_{g,x}(z)).
	\earyst
	Then by \eqref{eq tildexiperpxtilde0x} and \eqref{eq tpltperpequalstpltuperp} there is $\varkappa \in \{\pm 1\}$ such that
		\ary  \label{eq tpltperpxdiotax=pm1}
		\tplt^{u,\perp}_x((D\iota_x)_{(z,0,0)}(v^u_{g,x}(z), 1, v^c_{g,x}(z))) = \varkappa.
		\eary 
	By evaluating \eqref{eq xiusxiUxiuperp} at $(D\iota_x)_{(z, 0, 0)}(v^u_{g,x}(z), 1, v^c_{g,x}(z))$ and by \eqref{eq tpltperpxdiotax=pm1}, we obtain
	 \aryst
	 0 = \xi^{\cal U}_{x}(\Phi^u_x(z))[(D\iota_x)_{(z,0,0)}(v^u_{g,x}(z), 1, v^c_{g,x}(z))] + \varkappa  \phi^{u,s}_x(z).
	 \earyst
	 On the other hand, by \eqref{eq def of xiU}, \eqref{eq dhuescoordinates} and $\xi^{\cal U}_x \in \Tplt^u_1(x)$ we obtain that
	 \aryst
	 \xi^{\cal U}_{x}(\Phi^u_x(z))[(D\iota_x)_{(z,0,0)}(1,0,0)] &=& 0, \\
	  \xi^{\cal U}_{x}(\Phi^u_x(z))[(D\iota_x)_{(z,0,0)}(V^u_{g,x}(z), 1, V^c_{g,x}(z))] &=& 0, \\
	  \xi^{\cal U}_{x}(\Phi^u_x(z))[(D\iota_x)_{(z,0,0)}(0, 0, 1)] &=& 1.
	 \earyst
	 Thus we have
	 \ary \label{eq expressionofphius}
	 \phi^{u,s}_x(z) = \varkappa (
	 V^c_{g,x}(z) - v^c_{g,x}(z)).
	 \eary
	 We conclude the proof by \eqref{eq psit1t2}, \eqref{eq expressionofphius}, \eqref{def t1} and \eqref{def t2}.
\end{proof}

\begin{lemma}\label{lem 01lawforphius}
	If there exists $x \in M$ such that $\phi^{u,s}_x|_{J} = 0$ for interval $J = [0,1)$ or $(-1,0]$, then for every $y \in M$ we have $\phi^{u,s}_y = 0$.  In this case, $E^{u}$ and $E^s$ are jointly integrable.
\end{lemma}
\begin{proof}
	The proof is similar to that of Lemma \ref{lem alltrueallfalse}. 
	Denote by $X$ the set of $x \in M$ such that $\phi^{u,s}_x = 0$. Then by definition it is clear that $X$ is a forward invariant closed set. Moreover, if $x \in X$, then there exist $\varrho > 0$ and some $z \in W^u_{g, loc}(x)$ such that 
	$W^u_{g}(g^1(z), \varrho) \subset X$. But since $g$ is transitive, the union of the forward iterates of any local unstable leave is dense in $M$. By Remark \ref{rem templatec0onx2} we have $X = M$. This gives the first statement.

To show the joint integrability, we notice that for any $x \in M$ and any $z \in W^u_{g}(x, 1)$, we have  $DH^u_{x, z}(E^s(x)) = E^s(z)$.	
	It is known that for any $x \in M$, $E^s$ is a $C^{r-1}$ line field on $W^{cs}_g(x)$ which uniquely integrates to the stable foliation $W^s_g$ on $W^{cs}_g(x)$. Then for any $x \in M$, for any $z \in W^u_{g}(x, 1)$, $H^u_{x,z}$ maps $W^s_g$ on $W^{cs}_g(x)$ to $W^s_g$ on $W^{cs}_g(z)$. This implies that $E^s$ and $E^u$ are jointly integrable.
\end{proof}

\begin{cor}\label{lem jointus}
   If $E^u$ and $E^s$ are not jointly integrable, then there is $\kappa > 0$ such that for any $x \in M$, for $J = [0,1)$ or $(-1,0]$, we have $Osc_{J} \phi^{u,s}_x > \kappa$.
\end{cor}
\begin{proof} 
	By Remark \ref{rem templatec0onx2}, the function $\cD$ defined by \eqref{eq definecD} is continuous. Then by Lemma \ref{lem 01lawforphius}, we see that $Osc_{J} \phi^{u,s}_x > 0$ for every $x \in M$ and for $J = [0,1)$ or $(-1,0]$ whenever $E^u$ and $E^s$ are not jointly integrable. The corollary then follows from the compactness of $M$.
\end{proof}

%\subsection{Proof of Proposition \ref{prop decay near axis}}

\section*{Standing assumptions}
From Section \ref{sec construct partition} to Section \ref{sec Uniform non-integrability on uniform set}, we fix an Anosov flow $g$ and a H\"older potential $F$.
The following notations will be used throughout these sections. We let
\enmt
\item[$\bullet$] $K$ be an integer depending only on $g$ so that Corollary \ref{prop template approximation} holds for $K$;
\item[$\bullet$] $\theta \in  (0,1)$ be a small constant depending only on $g$ and $F$. In the course of the proof we will reduce $\theta$ if necessary, but will only do this finitely many times.  
In particular, all the smallness requirements about $\theta$ in Section \ref{sec markov} are satisfied, e.g., $F \in C^\theta(M, \R)$; and Lemma \ref{lem approximation of temporal case ii} holds for $\theta$.
\eenmt

\section{Choosing the correct scales} \label{sec construct partition}

\subsection{Slow growing sequence and stable sequence}

In this section, we will introduce the scales at which the cancellations for $\cL_{a,b}$ takes place.  

\begin{defi}\label{def slowseq}
We say that a sequence of $L^\infty$ functions $\{ h^\epsilon: M \to \Z_+ \}_{\epsilon > 0}$ is {\em slowly growing} for an Anosov flow $g$ if there is $\kappa > 0$ such that for all sufficiently small $\epsilon$
\aryst %\label{eq slowgrowingaxiom0}
\inf_{x \in M} h^\epsilon(x) > \kappa|\log \epsilon|,
\earyst
and for any integer $p > 0$, for any sufficiently large integer $k > 0$, for any sufficiently small $\epsilon$
	\ary \label{eq slowgrowingproperty}
	h^\epsilon(g^{-k}(x)) \leq h^\epsilon(x) +  k - p, \quad \forall x \in M.
	\eary 
\end{defi}

We have the following.
\begin{lemma} \label{eq lyapunovfunctionstable}
     For any slowly growing sequence $\{ h^\epsilon: M \to \Z_+ \}_{\epsilon > 0}$,
     the sequence of functions
	$\{ \Lambda^\epsilon: U \to \R_+ \}_{\epsilon > 0}$ defined by
	\ary \label{eq lambdah}
	\Lambda^\epsilon(x) = \sup_{y \in W^s_g(x,1)} \Lambda_{h^\epsilon(y)}(x)
	\eary
	is stable (given in Definition \ref{def stableandtame}).
\end{lemma}
\begin{proof} 
	By Definition \ref{def slowseq}, there is $\kappa > 0$ such that for all sufficiently small $\epsilon$, there exists $C_\epsilon > 0$ such that 
	\aryst
	C_\epsilon > 
	\Lambda^\epsilon(x) > \epsilon^{-\kappa}, \quad \forall x \in U.
	\earyst 
	Indeed, the left inequality above follows from that $h^\epsilon \in L^\infty(M)$, and the right inequality above follows from the first inequality in Definition \ref{def slowseq}. 
	
	Let $p > 0$ be a large integer to be determined. Let $n > 1$ be a large integer, and denote $t = \tau_n(x)$. Clearly, $t \geq n \tau_0$. 
	By \eqref{eq wsgx1 is long} and the Markov property,  we have $g^t(W^s_{g}(x, 1)) \subset W^s_{g}(\sigma^n(x), 1)$ by letting $n$ be sufficiently large.
	Take an arbitrary $y \in W^s_{g}(x, 1)$. By letting $n$ be sufficiently large, and letting $\epsilon$ be sufficiently small,  we have $h^\epsilon(y) \leq h^\epsilon(g^t(y)) + t - p$ for all sufficiently small $\epsilon$. Consequently, by distortion estimate, we have 
	\aryst
	\Lambda^\epsilon(x) e^{p\chi_0} &\leq& C' \Lambda^\epsilon(\sigma^n(x))\sup_{y \in W^s_{g}(x, 1)}\norm{Dg^t|_{E^u(y)}}, \\
	\sup_{y \in W^s_{g}(x, 1)}\norm{Dg^t|_{E^u(y)}} &\leq& C'\inf_{y \in W^s_{g}(x, 1)}\norm{Dg^t|_{E^u(y)}}.
	\earyst
	By letting $p$ be sufficiently large, we obtain for all sufficiently small $\epsilon$ that
	\aryst
	\Lambda^\epsilon(x) \leq \frac{1}{2}\Lambda^\epsilon(\sigma^n(x)) \inf_{y 
		\in W^s_{g}(x, 1)}\norm{Dg^{\tau_n(x)}|_{E^u(y)}}.
	\earyst
	By \eqref{eq wsgx1 is long}, we can iterate the above inequality to show that $\{ \Lambda^\epsilon \}_{\epsilon > 0}$ is stable.
\end{proof}

\subsection{Choosing the correct scales - Class $I$}
In this subsection, we assume that $g \in \cU^\infty$ is in Class $I$. 
By Corollary \ref{cor case i to iii}, we see that $g \in \cU^\infty_{+,K} \cap \cU^\infty_{-, K}$.  By Lemma \ref{lem alltrueallfalse}, we can choose $\kappa_1 > 0$ such that $g \in \cU^\infty_{+, K, \kappa_1} \cap \cU^\infty_{-, K, \kappa_1}$. 
 
	Let $C_3 > 1$ and $\delta_5 > 0$ be determined in due course.   
\begin{defi} \label{def hepsilon}
	Given a sufficiently small $\epsilon > 0$, for any $x \in M$, we let
	\enmt
	\item $k^\epsilon(x)$ be the smallest integer $n \geq 1$ such that
	 $ \norm{Dg^n|_{E^s(x)}}, \norm{Dg^n|_{E^u(x)}}^{-1} < \epsilon$;
	\item  
	the {\em matching time of order $\epsilon$} at $x$, denoted by $\varsigma^\epsilon(x)$, be the smallest integer $n \geq k^\epsilon(x)$ satisfying
	that there is $\varkappa \in \{\pm 1\}$ such that for every $y \in (-1,1)$, there exists $\varphi \in \cGs_{g^n(x)}$ such that
	\ary
	\norm{\Psi_x(\varkappa \Lambda_n(x)^{-1} \cdot, y) -  \varphi}_{(-2,2)} &\leq& C_3 ((\epsilon|y|)^{1+\delta_5} + \epsilon^{2}), \label{eq def matchingtime 1} \\
	\norm{\varphi}_{(-2,2)} &\leq& \max(\epsilon |y|^{\delta_5/2}, C_3\epsilon|y|).\label{eq def matchingtime 2}
	\eary
	\eenmt
	For every $x \in U$,
	the {\em matching scale of order $\epsilon$} at $x$ is defined by
	\ary \label{def lambda epsilon}
	\Lambda^\epsilon(x) = \sup_{y \in W^s_g(x,1)}  \Lambda_{\varsigma^\epsilon(y)}(x).
	\eary
\end{defi}

The following lemma shows that for suitable choices of $C_3$ and $\delta_5$, the sequence of matching time of order $\epsilon$ for $\epsilon > 0$ is slowly growing. We will show at the end of this section, the matching scale of order $\epsilon$ for $\epsilon > 0$ is a stable sequence.

\begin{lemma} \label{lem varsigmaslowgrowing}
	Let $\delta_5 > 0$ be sufficiently small depending only on the exponents in \eqref{eq rate stable unstable},
	let $K > 0$ be a sufficiently large integer depending only on $g$ and $\delta_5$, and let $C_3 > 1$ be sufficiently large depending on the implicit constants in Lemma \ref{lem psixnorm1} and Corollary \ref{prop template approximation}.  
	Then the sequence of functions $\{ \varsigma^\epsilon \}_{\epsilon > 0}$ in Definition \ref{def hepsilon} is slowly growing.
\end{lemma}
\begin{rema}
	In the rest of the paper, we will always choose $C_3$ and $\delta_5$ relative to an Anosov flow $g \in \cU^\infty$ in Class $I$ such that Lemma \ref{lem varsigmaslowgrowing} holds. We will enlarge $C_3$ and reduce $\delta_5$ if necessary, but we will only do this finitely many times, and both $C_3$ and $ \delta_5$ depend only on $g$.
\end{rema}

\begin{proof}[Proof of Lemma \ref{lem varsigmaslowgrowing}]
	We first notice the following lemma.
	\begin{lemma} \label{lem lower bound by uni for templates}
		There exist an increasing function $c_0: (0,1) \to (0,1]$ and an increasing function $c_1: [0,1] \to [0,1]$ with $\lim_{s \to 0} c_1(s) = 0$ such that
		for any $\varrho > 0$, any $x \in M$, any $\varphi \in \cGs_x$ or $\cGu_x$ and any $s_0 \in (-\rhone + \varrho, \rhone - \varrho)$, we have
		\aryst
		\|\varphi\|_{B(s_0, \varrho)} &\geq& 
c_0(\varrho) \| \varphi \|_{(-\rhone, \rhone)}, \\
        \|\varphi\|_{B(0, \varrho)} &\leq& c_1(\varrho) \| \varphi \|_{(-\rhone, \rhone)}.
		\earyst
	\end{lemma} 
\begin{proof}

	Fix $x \in M$.
	Without loss of generality, let us assume that $\varphi \in \cGs_x \setminus Poly^K$ and  
	\aryst
	\varphi = 
	\varphi^{u,s}_x + P
	\earyst
	where $P \in Poly^K$ and $P(0) = 0$. The argument for $\varphi \in Poly^K$ is straightforward; and the argument for $\varphi \in \cGu_x$ is similar to the one below.
		For simplicity, let us also assume that $s_0 = 0$. The argument for general $s_0$ is similar to the one below.

		We collect some basic properties of polynomials: for any $\varrho \in (0,\varrho_1), D > 0$ and any integer $L > 1$,
	\enmt
	\item
	there is a constant $C_{\varrho, D, L} > 0$  such that every $P \in Poly^L \setminus Poly^L_{C_{\varrho, D, L}}$ satisfies
	\aryst \label{eq poly prop 1}
	\norm{P}_{(-\varrho, \varrho)} > D.
	\earyst
	\item 
	there is a constant $B_{\varrho, L} > 0$ such that every $P \in Poly^L$ satisfies
	\aryst \label{eq poly prop 2}
	\norm{P}_{(-\varrho, \varrho)} \geq B_{\varrho, L} \norm{P}_{(-\varrho_1, \varrho_1)}.
	\earyst
	\item there is a constant $E_{L} > 0$ such that  every $P \in Poly^L$ satisfies
	\aryst \label{eq poly prop 3}
	\norm{P}_{(-\varrho, \varrho)} \leq E_{L} |\varrho| \norm{P}_{(-\varrho_1, \varrho_1)}.
	\earyst
	\eenmt
	
	We set
	\aryst
	D = 1 + \sup_{x \in M}
	\norm{\varphi^{u,s}_x}_{C^\theta(-\rhone, \rhone)}.
	\earyst
	
	 Let us first detail the proof of the first inequality. 
	We first assume that $P \in Poly^K \setminus Poly^K_{C_{\varrho, 10D, K}}$.
	Then we have
	\aryst
	\norm{P}_{(-\varrho, \varrho)} > 10D, \quad
	\norm{P}_{(-\varrho, \varrho)} \geq B_{\varrho, K} \norm{P}_{(-\rhone, \rhone)}.
	\earyst
	Then
	\aryst
	\norm{\varphi}_{(-\rhone,\rhone)} \leq 2\norm{P}_{(-\rhone,\rhone)}, \quad  \norm{P}_{(-\varrho,\varrho)}\leq 2\norm{\varphi}_{(-\varrho, \varrho)}.
	\earyst
	Thus we have
	\aryst
	\norm{\varphi}_{(-\varrho, \varrho)}  
	> \frac{1}{2}\norm{P}_{(-\varrho, \varrho)} 
	> \frac{1}{2}B_{\varrho, K} \norm{P}_{(-\rhone, \rhone)} 
	> \frac{1}{4}B_{\varrho, K} \norm{\varphi}_{(-\rhone, \rhone)}. 
	\earyst
	
	Now we consider the case where $P \in Poly^K_{C_{\varrho, 10D, K}}$.
	Clearly there is $M_{\varrho, D, K} > 0$ such that
	we have
	\aryst
	\norm{\varphi}_{(-\rhone,\rhone)} \leq \norm{P}_{(-\rhone,\rhone)} + D < M_{\varrho, D, K}.
	\earyst
	It remains to show that there exists $\kappa_\varrho > 0$ such that
	for any $P \in Poly^K$, we have
	\aryst
	\norm{\varphi}_{(-\varrho, \varrho)} > \kappa_\varrho.
	\earyst
	Define $\lambda_m$ and $\mu_m$ as in \eqref{eq xmfmmumlambam}.
	We let $q \geq 1$ be an integer such that
	\aryst
	|\lambda_0 \cdots \lambda_{q-1}| < \varrho \leq |\lambda_0 \cdots \lambda_{q-2}|.
	\earyst
	It is clear that we have
	\aryst
	|\mu_0 \cdots \mu_{q-1}|^{-1} \geq \varrho^{C'}/C'.
	\earyst
	By Corollary \ref{lem stable templates transformation},  there exists $Q \in Poly^K$ such that for every $z \in (-\varrho_1, \varrho_1)$
	\aryst
	\varphi(z) &=& \varphi^{u,s}_{x}(z) + P(z) \\
	&=& \mu_0^{-1} \cdots \mu_{q-1}^{-1} (\varphi^{u,s}_{g^q(x)}(\lambda_0^{-1} \cdots \lambda_{q-1}^{-1} z) +  Q(z) ).
	\earyst
	Thus by $g \in \cU^r_{+, 	K, \kappa_1}$ and Definition \ref{def uniplus} we have
	\aryst
	\norm{\varphi}_{(-\varrho, \varrho)} &\geq& |\mu_0 \cdots \mu_{q-1}|^{-1} \norm{\varphi^{u,s}_{g^q(x)} + Q(\lambda_0 \cdots \lambda_{q-1} \cdot)}_{(-1, 1)} \\
	&\geq&  \varrho^{C'} \kappa_1/C'.
	\earyst
    This concludes the proof of the first inequality.

	We now sketch the proof of the second inequality.
	Since $\varphi^{u,s}_x(0) = 0$ and $\norm{\varphi^{u,s}_x}_\theta < D$, we have
	\ary \label{eq vaprhiupperboundonrho}
	\norm{\varphi}_{(-\varrho, \varrho)} &\leq& \norm{\varphi^{u,s}_x}_{(-\varrho, \varrho)}  + \norm{P}_{(-\varrho, \varrho)} \leq D|\varrho|^\theta + E_L |\varrho| \norm{P}_{(-\varrho_1, \varrho_1)}.
	\eary
	If $\norm{P}_{(-\varrho_1, \varrho_1)} > 10D$, then $\norm{\varphi}_{(-\varrho_1, \varrho_1)} \geq \frac{1}{2}\norm{P}_{(-\varrho_1, \varrho_1)}$, and consequently 
	\aryst
	\norm{\varphi}_{(-\varrho, \varrho)} &\leq& (|\varrho|^\theta + E_L |\varrho|)  \norm{P}_{(-\varrho_1, \varrho_1)} \leq 2(|\varrho|^\theta + E_L |\varrho|)  \norm{\varphi}_{(-\varrho_1, \varrho_1)}.
	\earyst
	If $\norm{P}_{(-\varrho_1, \varrho_1)} \leq 10D$, then by \eqref{eq vaprhiupperboundonrho} and $g \in  \cU^\infty_{+,K, \kappa_1}$, we have
	\aryst
	\norm{\varphi}_{(-\varrho, \varrho)} \leq 10(|\varrho|^\theta + E_L |\varrho|)D \leq 10(|\varrho|^\theta + E_L |\varrho|)D\kappa_1^{-1} \norm{\varphi}_{(-\varrho_1, \varrho_1)}.
	\earyst
	Combine the above estimates, we conclude the proof of the second inequality.
	\clb	
\end{proof}

We fix $x \in M$ and $x_m$, $\lambda_m$ and $\mu_m$ as in \eqref{eq xmfmmumlambam}.
Let $\epsilon > 0$ be a small constant. We set 
\aryst
n = \varsigma^\epsilon(x).
\earyst 
By Definition \ref{def hepsilon}, we have
\aryst
|\lambda_0 \cdots \lambda_{n-1}|, |\mu_0 \cdots \mu_{n-1}|^{-1} < \epsilon.
\earyst
Then clearly we have
\aryst
n >  (C')^{-1}|\log \epsilon|.
\earyst
By Corollary \ref{prop template approximation}, 
 for any $n' > 0$  sufficiently large depending on $g$ and $\epsilon$, we have
\aryst
\norm{\Psi_x(\Lambda_{n'}(x)^{-1}\cdot, y)}_{(-2,2)} \leq  C'(n'e^{-n'\chi_0} + D_{n'}) <  \epsilon^{2}, \quad \forall y \in (-1,1).
\earyst  
This implies that $\varsigma^\epsilon \in L^\infty(M)$.

We now verify \eqref{eq slowgrowingproperty} for $\varsigma^\epsilon$ in place of $h^\epsilon$.
Fix an arbitrary $y \in (-1,1)$, an integer $p > 0$, and let $k > 0$ be a large integer to be determined depending on $g$ and $p$. Denote 
\aryst
m = n - p, \quad 
\bar y = (\mu_{-k} \cdots \mu_{-1})^{-1}y, \quad
\varrho = |\lambda_m \cdots \lambda_{n-1}|.
\earyst
 Note that $m > 0$ since we assume that $\epsilon$ is small.

By \eqref{lem finvariance}, there exist $\varkappa_1, \varkappa_2 \in \{\pm 1\}$ such that for any $z \in (-\varrho_1, \varrho_1)$
\ary 
\label{eq psixthreeequalities}  
\ \ \   \Psi_{x_{-k}}(\Lambda_{k+m}(x_{-k})^{-1} z, y) 
&=& \Psi_{x}( \varkappa_1 \Lambda_{m}(x)^{-1}  z, \bar y)  \nonumber \\
&=& \Psi_{x}( \varkappa_2 \Lambda_n(x)^{-1} \varrho^{-1} z, \bar y).
\eary
By letting $k$ be sufficiently large depending only on $g$ and $p$, we have
\aryst
|\lambda_{-k} \cdots \lambda_{m-1}|, |\mu_{-k} \cdots \mu_{m-1}|^{-1} < \epsilon.
\earyst
Then by Definition \ref{def hepsilon}, it remains to show that \eqref{eq def matchingtime 1} and \eqref{eq def matchingtime 2} hold for $(x_{-k}, k+m)$ in place of $(x,n)$.

If  $|y|^{\delta_5} < \epsilon$, then by Lemma \ref{lem psixnorm1} and by letting $\delta_5 < \delta_2/3$ ($\delta_2$ is given by Lemma \ref{lem psixnorm1}), we obtain for all sufficiently small $\epsilon$ that
\aryst
\norm{\Psi_{x_{-k}}(\varkappa_3 \Lambda_{k+m}(x_{-k})^{-1} \cdot, y)}_{(-2,2)} < C' |y|^{\delta_2} \leq C' \epsilon^{\delta_2/\delta_5} < \epsilon^2.  
\earyst
In this case we may take $\varphi = 0$ in \eqref{eq def matchingtime 1} and \eqref{eq def matchingtime 2}.

Now assume that $|y|^{\delta_5} \geq \epsilon$ and  $\delta_5 < \delta_0$ where $\delta_0$ is given by Corollary \ref{prop template approximation}.
By Lemma \ref{lem psixnorm1}, Corollary \ref{prop template approximation},
% (more precisely, by $\lim_{k \to \infty} D_{k+m} = 0$), 
and  by letting $k$ be sufficiently large depending only on $g, p, C_3$, there exist $\varkappa_3 \in \{\pm 1\}$, $\varphi_1 \in \cGs_{x_m}$ such that 
\ary \label{eq compare 1}
\norm{\Psi_{x_{-k}}(\varkappa_3 \Lambda_{k+m}(x_{-k})^{-1} \cdot, y) - \varphi_1}_{(-2,2)} &\leq& C_3((\epsilon|y|)^{1+\delta_0} + \epsilon^2). 
\eary

By Definition \ref{def hepsilon}, there exist $\varkappa_4 \in \{\pm 1\}$ and $\varphi_2 \in \cGs_{x_n}$ such that 
\ary \label{eq compare 2}
\norm{\Psi_x(\varkappa_4 \Lambda_n(x)^{-1}\cdot, \bar y)  -  \varphi_2}_{(-2,2)} &<& C_3 ( (\epsilon|\bar y|)^{1+\delta_5} + \epsilon^{2}), \\
\norm{\varphi_2}_{(-1,1)} <  \max(\epsilon |\bar y|^{\delta_5 /2}, C_3\epsilon|\bar y|) &<& e^{-k\chi_0\delta_5/2} \max(\epsilon |y|^{\delta_5/2}, C_3\epsilon|y|). \label{eq compare 3}
\eary
Then by \eqref{eq psixthreeequalities}, \eqref{eq compare 1} and \eqref{eq compare 2} we have
\ary \label{eq varphivarrhovarphi}
\norm{ \varphi_1(\varkappa_3 \varrho \cdot ) - \varphi_2(\varkappa_2\varkappa_4 \cdot) }_{(-2,2)} < 2  C_3( (\epsilon|y|)^{1+\delta_5} + \epsilon^{2} ).
\eary
Since  $|y|^{\delta_5} \geq \epsilon$, we have for all sufficiently small $\epsilon$ that
\ary \label{eq minmax2}
(\epsilon|y|)^{1+\delta_5} + \epsilon^{2} <\epsilon^{\delta_5/2} (\epsilon|y|^{\delta_5/2} +  C_3 \epsilon |y|).
\eary

It is clear that $\varrho > e^{-p\chi_*}$.
By
Lemma \ref{lem lower bound by uni for templates}, \eqref{eq compare 3}, \eqref{eq varphivarrhovarphi} and \eqref{eq minmax2} we obtain
\aryst
\norm{\varphi_1}_{(-2,2)} &\leq& c_0(\varrho)^{-1} \norm{\varphi_1}_{(-\varrho, \varrho)} \leq c_0(\varrho)^{-1} ( \norm{\varphi_1(\varkappa_3 \varrho  \cdot) - \varphi_2(\varkappa_2\varkappa_4 \cdot)}_{(-2,2)} + \norm{\varphi_2}_{(-2, 2)} ) \\
&\leq& c_0(\varrho)^{-1}( 2  C_3( (\epsilon|y|)^{1+\delta_5} + \epsilon^{2} ) + e^{-k\chi_0\delta_5/2} \max(\epsilon |y|^{\delta_5/2}, C_3\epsilon|y|)) \\
&\leq& c_0(e^{-p\chi_*})^{-1}( 4  C_3 \epsilon^{\delta_5/2} + e^{-k\chi_0\delta_5/2} )\max(\epsilon |y|^{\delta_5/2}, C_3\epsilon|y|).
\earyst
By letting $k$ be sufficiently large depending on $g, \delta_5$ and $p$, and by letting $\epsilon$ be small depending on $g, C_3$ and $p$, we have
\aryst
\norm{\varphi_1}_{(-2,2)} < \max(\epsilon |y|^{\delta_5/2}, C_3\epsilon|y|).
\earyst
In this case, we may let $\varphi = \varphi_1$ in \eqref{eq def matchingtime 1} and \eqref{eq def matchingtime 2}.

Since $y$ is taken arbitrarily in $(-1,1)$, we have
\aryst
\varsigma^\epsilon(x_{-k}) \leq k + m = \varsigma^\epsilon(x) + k - p.
\earyst
Thus $\{ \varsigma^\epsilon \}_{\epsilon > 0}$ is slowly growing.
\end{proof}

 We obtain the following corollary. 
\begin{cor} \label{cor lambdaepsilonstable}
	For any $g \in \cU^\infty$ in Class $I$,
	the sequence of functions
	$\{ \Lambda^\epsilon \}_{\epsilon > 0}$ is stable and tame.
\end{cor}
\begin{proof} 
	By Lemma \ref{eq lyapunovfunctionstable} and Lemma \ref{lem varsigmaslowgrowing}, we see that the sequence of functions 
	$\{ \Lambda^\epsilon \}_{\epsilon > 0}$ is stable.
	Then by Definition \ref{def hepsilon}, there exist an integer $n \geq 1$ and $\varkappa \in \{  \pm 1 \}$ such that for any $y \in (-1,1)$, there exists $\varphi \in \cT^s_{g^n(x)}$ with
\ary \label{tame 1}
\norm{\Psi_x(\varkappa \Lambda^\epsilon(x)^{-1} \cdot, y) -  \varphi}_{(-1,1)} &<& 2C_3\epsilon^{1+\delta_5},  \\
\norm{\varphi}_{(-1,1)} &<& C_3\epsilon|y|^{\delta_5/2}. \label{tame 12}
\eary
Write $\varphi = h\varphi^{u,s}_{g^n(x)} + P$ for some $h \in \R$ and $P \in Poly^K$.
By $g\in \cU^\infty_{+, K, \kappa_1}$  we have
\aryst
\norm{\varphi}_{(-1,1)} \geq h\kappa_1
\earyst
and consequently, by \eqref{tame 12}, we have $h < C'C_3\kappa_1^{-1}\epsilon|y|^{\delta_5/2}$.
Again by \eqref{tame 12}, we see that
\aryst
\|P\|_{(-1,1)} \leq \| h \varphi^{u,s}_{g^n(x)} \|_{(-1,1)} + \| \varphi \|_{(-1,1)} \leq C' C_3 \epsilon|y|^{\delta_5/2}.
\earyst
Thus
\ary \label{tame 2}
\norm{\varphi}_\theta &\leq&  \kappa_1^{-1}C'C_3 \epsilon|y|^{\delta_5/2}. 
\eary

By the fact that $W^{cs}_g$ is a $C^1$ foliation, we can find $R \in C^\theta(J^{\Lambda^\epsilon}_x)$ such that \eqref{tame 1} and \eqref{tame 2} hold for $\Psi_{x, \Phi^s_x(y)}(\Phi^{\Lambda^\epsilon}_x(\cdot))$ and $\epsilon R$ in place of $\Psi_x(\Lambda^\epsilon(x)^{-1}\cdot , y)$ and $\varphi$ respectively.
Consequently $\{ \Lambda^\epsilon \}_{\epsilon > 0}$ is tame.
\end{proof}

\subsection{Choosing the correct scales - Class {\it II}} \label{subsec scaleclassii}
In this subsection, we assume that $g \in \cU^\infty$ is in Class {\it II}. In another words, $E^u$ is $C^{1+\delta}$ for some $\delta > 0$.

\begin{defi}\label{def hepsilon 1}
	For any $\epsilon > 0$ we define $\vartheta^\epsilon: M \to \Z_+$ by
	\ary
	\label{eq def of vartheta}
	\vartheta^\epsilon(x) = \inf \{ k \geq 1 \mid  \norm{Dg^k|_{E^s(x)}} < \epsilon \}.
	\eary
	We define
\aryst
\Lambda^\epsilon(x) = \Lambda_{\vartheta^\epsilon(x)}(x).
\earyst
\end{defi}

\begin{lemma}\label{lem varsigmaslowgrowing case ii}
	By letting $\theta$ be sufficiently small depending only on $g$, the following is true.
	For any $x \in M$, there is a function $R \in C^\theta(-2,2)$ 
	such that
	\aryst
	\norm{R}_\theta &<& C', \\
	\epsilon^{-1}\Psi_{x}( \Lambda^\epsilon(x)^{-1}z, y) &=& y R(z) + O( \epsilon^{\delta}), \quad \forall y \in (-1,1), z \in (-2, 2).
	\earyst
	Moreover, if $E^s$ and $E^u$ are not jointly integrable, then there is  $\kappa_3 > 0$ such that for all sufficiently small $\epsilon > 0$, for $J_0 = [0,1)$ or $(-1,0]$, we have
	\aryst
	Osc_{J_0}R > \kappa_3.
	\earyst
\end{lemma}
\begin{proof}
	Denote $n = \vartheta^\epsilon(x)$ and $x_n  =g^n(x)$. 
	As before we denote $\lambda_i = \lambda_{g^i(x)}$ and $\mu_i = \mu_{g^i(x)}$ for all integer $i$.
	By \eqref{lem finvariance}, there is $\varkappa \in \{  \pm 1 \}$ such that for all $z \in (-2,2)$, $y \in (-1,1)$ we have
	\aryst
	\Psi_x(\Lambda^\epsilon(x)^{-1}z, y) = 
	\Psi_{x}(\varkappa \lambda_0 \cdots \lambda_{n-1}z, y) = 
	\Psi_{x_n}(\varkappa z, (\mu_0 \cdots \mu_{n-1})^{-1}y).
	\earyst
	By Lemma \ref{lem approximation of temporal case ii}, we have
	\aryst
	\epsilon^{-1}\Psi_{x_n}(z, (\mu_0 \cdots \mu_{n-1})^{-1}y) =  cy \phi^{u,s}_{x_n}(z) + O( \epsilon^{\delta}).
	\earyst
	where we have
	\aryst
	|c| = |\epsilon^{-1}(\mu_0 \cdots \mu_{n-1})^{-1}| \sim 1.
	\earyst
	By Lemma \ref{lem approximation of temporal case ii} and by letting $\theta$ be sufficiently small depending only on $g$, we have $\norm{\phi^{u,s}_{x_n}}_\theta <  C'$.
	Then  the first statement of the lemma is satisfied for
	\aryst
	R(z) = c \phi^{u,s}_{x_n}(\varkappa z), \quad \forall z \in (-2,2).
	\earyst
	The last statement then follows from Corollary \ref{lem jointus}.
\end{proof}

Similar to Corollary \ref{cor lambdaepsilonstable}, we have the following.
\begin{cor} \label{lemma lambdaepsilonstable}
	For any $g \in \cU^\infty$ in Class {\it II}, the sequence of functions 
	$\{ \Lambda^\epsilon \}_{\epsilon > 0}$ is stable and tame.
\end{cor}

\section{Uniform non-integrability on uniform set}\label{sec Uniform non-integrability on uniform set}

\subsection{ Uniform non-integrability in Class $I_F$} \label{subsec nonexpanding I}
 We first assume that $g$ is in Class $I$.
 By Corollary \ref{cor case i to iii}, there is $\kappa_1 > 0$ so that
$g \in \cU^\infty_{+, K, \kappa_1} \cap \cU^\infty_{-, K, \kappa_1}$.

We first state a property of the functions defined in Definition \ref{def txn}.
\begin{lemma} \label{lem oscillation exists}
	For any $\varrho \in (0, \rhone)$, there exists $\kappa > 0$ depending only on $g, K, \varrho$ such that
	for any $(z_0, y_0) \in (-\varrho_1 + \varrho, \varrho_1 - \varrho)^2$,
	 any integer $n > 0$, any $x \in M$,
	 and any $R \in \cT_{x,n}$, we have
	\aryst
	Osc_{B(z_0, \varrho) \times B(y_0, \varrho)} R  \geq \kappa \norm{R}_{(-\varrho_1, \varrho_1)^2}.	
	\earyst
\end{lemma}

\begin{proof}
	Denote $x_n  = g^n(x)$.
	Take an arbitrary 
	\aryst
	R(z,y) = h_1 y \varphi^{u,s}_{x_n}(z) + h_2 z \varphi^{s,u}_x(y) + P(z,y) \in \cT_{x,n}.
	\earyst
    We may assume without loss of generality that $\norm{R}_{(-\rhone, \rhone)^2} = 1$.
	
	We let $C > 1$ be a large constant to be determined depending only on $g, \varrho_1$ and $K$.  
	Assume that
	\ary
	\norm{\varphi^{u,s}_{w}}_{(-\rhone, \rhone)},
	 \norm{\varphi^{s,u}_w}_{(-\rhone, \rhone)}  < C^{1/10}, \quad \forall w \in M. \label{eq upperboundc7110}
	\eary
	We divide the proof into two cases.

	\noindent{$(1)$} We first assume that $|h_1| > C^{-1}$ or $|h_2| > C^{-1}$.
	Without of loss of generality, let us assume that $|h_1| > C^{-1}$ (the other case can be handled in a similar way). 
	By $g \in \cU^\infty_{+, K, \kappa_1}$ and by Lemma \ref{lem lower bound by uni for templates}, we have for any $y$ with $\varrho/2 < |y- y_0| < \varrho$ that 
	\aryst
	Osc_{B(z_0, \varrho)} ( R(\cdot,y) - R(\cdot, y_0) ) \geq \frac{c_0(\varrho) \varrho \kappa_1}{2C}.
	\earyst	

	\noindent{$(2)$} Now assume that  $|h_1|, |h_2| \leq C^{-1}$. Then  by \eqref{eq upperboundc7110} we have
	\ary \label{eq differencerp}
	\norm{R - P}_{(-\varrho_1, \varrho_1)^2} < 2C^{-1/10}\varrho_1.
	\eary
	By letting $C$ be sufficiently large depending on $\varrho_1$, we have
	\aryst
	\norm{P}_{(-\varrho_1, \varrho_1)^2} \geq 	\norm{R}_{(-\varrho_1, \varrho_1)^2} -	\norm{R-P}_{(-\varrho_1, \varrho_1)^2}  \geq \frac{1}{2}.
	\earyst
	Then there is a constant $c_{K, \varrho} > 0$ depending only on $K,\varrho$ such that
	\aryst
	Osc_{B(z_0, \varrho) \times B(y_0, \varrho)} P  > c_{K,\varrho}.	
	\earyst
	Lemma \ref{lem oscillation exists} follows from \eqref{eq differencerp} and by letting $C$ be sufficiently large.
\end{proof}

We consider the following subsets of $U$.

\begin{defi}\label{def uniform set}

	Given $\epsilon, \kappa > 0$ and an integer $n > 0$, we define the $(n, \kappa)$, resp. $(\epsilon, n, \kappa )$-{\em  non-expanding uniform set} by
	\aryst 
	\Omega( n, \kappa) &=& \{ x  \mid   \det Dg^{i}(x) < e^{ i \kappa}  \mbox{ for any integer } i > n \}, \\
		\mbox{ resp. } \ \ 	\Omega(\epsilon, n, \kappa) &=& \{ x  \mid  \det Dg^{i}(x) < e^{ i \kappa}  \mbox{ for any integer } i \in (n, \varsigma^\epsilon(x)]  \}.
	\earyst
\end{defi}
We will show in Section \ref{sec recurrence} that $\Omega(n, \kappa)$ is $n_1$-recurrent for some constant $n_1 > 0$ whenever $g$ is in Class $I_F$ and $n$ is large. 

%\begin{defi} \label{def of Omega}
%For any integer $n > 1$, any $C_1, \kappa  > 0$, and for all sufficiently small $\epsilon > 0$, we set $\Lambda^\epsilon$ as in Definition \ref{def hepsilon} for some sufficiently large $C_3$. We denote by $\Omega(C_1, \epsilon, n, \kappa)$
%	the set of $x \in U$ such that $W^u_g(x, C_1 \Lambda^\epsilon(x)^{-1}) \cap \Omega(n, \kappa) \neq \emptyset$.
%\end{defi}

We have the following lemma. Recall that we let $\delta$ denote a positive constant depending only on $g$.
\begin{lemma} \label{lem psi to r}
	For any sufficiently small $\kappa > 0$ and any integer $n_0 > 0$ the following is true.
	 For any $C_6 > 0$, for all sufficiently small $\epsilon > 0$,
	for any $x \in U$ and  integer $n > 0$ satisfying either of the conditions
	\enmt
	\item  $x \in \Omega(\epsilon, n_0, \kappa)$ and $k^\epsilon(x) - C_6 \leq n \leq \varsigma^\epsilon(x) + C_6$,
\mbox{ or }	\item $x \in \Omega(n_0, \kappa)$ and $n \geq k^\epsilon(x) - C_6$,
	\eenmt  
	there exist $R \in \cT_{x,n}$ and $\varkappa \in \{\pm 1\}$ such that
	\aryst
	\norm{\Psi_{x}(\varkappa \Lambda_n(x)^{-1} \cdot, \cdot) - R}_{(-\rhone, \rhone)^2} = O(\epsilon^{1+\delta}).
	\earyst
	Moreover, there is $C'_3 > 0$ depending only on $g$ such that if in addition to the above we also have  
	\ary \label{eq psi small in uniform}
%	n \geq k^\epsilon(x) \quad  \mbox{and} \quad 
	\norm{\Psi_{x}(\Lambda_n(x)^{-1}\cdot, \cdot)}_{(-2,2)^2} < \epsilon/C'_3,
	\eary
	then  $n \geq \varsigma^\epsilon(x)$.
\end{lemma}
\begin{proof}
	Given $x \in U$, let us denote $\lambda_i, \mu_i$ as  in \eqref{eq xmfmmumlambam}.
%	Let 
%	\aryst
%	h_1 = (\mu_0 \cdots \mu_{n-1})^{-1} \quad  \mbox{and} \quad  h_2 = \lambda_0 \cdots \lambda_{n-1}.
%	\earyst
	By  Corollary \ref{prop template approximation}, there exist
    functions $a_2, \cdots, a_{K} : (-\varrho_1, \varrho_1) \to \R$  satisfying \eqref{eq aiupperbound}, $\varkappa \in \{ \pm 1\}$ and $R \in \cT_{x,n}$ such that for any $z,y \in (-\varrho_1, \varrho_1)$
	\ary \label{eq psixn app}
	\Psi_{x}(\varkappa \Lambda_n(x)^{-1} z, y) &=& R(z,y)  + \sum_{i=2}^{K} a_i(y) z^i + O((\epsilon|y|)^{1+\delta} + \epsilon^{2}). 
	\eary
	
	We have the following.
	\begin{claim} \label{claim for ai}
	For $2 \leq i \leq K$, for $y \in (-\rhone, \rhone)$, we have
	\ary \label{eq upperboundforaiigeq2}
	|a_i(y)| = O(\epsilon^{1+ \delta}|y|).
	\eary
	
\end{claim}
\begin{proof}[Proof of the claim.]
	 Under either (1) or (2) of the lemma, we have $|\lambda_0 \cdots \lambda_{n-1}| \leq  C'' \clb \epsilon$ and
  	\aryst
	|\mu_0 \cdots \mu_{m-1}\lambda_0 \cdots \lambda_{m-1}|^{-1} \sim  \det Dg^m( x) \leq C''e^{m\kappa}, \quad \forall 1 \leq m \leq n
	\earyst
	for some $C''$ depending on $g, C_6$ and $n_0$.
	 Moreover, we may assume without loss of generality that $n < C' |\log \epsilon|$, for otherwise the claim is trivial by \eqref{eq aiupperbound}.

Let $\eta_0$ be given by Proposition \ref{prop template approximation 0}.	
		By \eqref{eq aiupperbound},  for any $2 \leq i \leq K$, for any $y \in (-\rhone, \rhone)$, we have
	\ary
	&&	|a_i(y)| 
	\leq C'|y|\sum_{m=0}^{\lfloor (1-\eta_0)n \rfloor-1} |\mu_0 \cdots \mu_{m-1}|^{-1}|\lambda_m \cdots \lambda_{n-1}|^i  \nonumber \\
	&\leq& C'|y|\sum_{m=0}^{\lfloor (1-\eta_0)n \rfloor-1} |\lambda_0 \cdots \lambda_{n-1}||\mu_0 \cdots \mu_{m-1}\lambda_0 \cdots \lambda_{m-1}|^{-1}|\lambda_m \cdots \lambda_{n-1}|^{i -1}  \nonumber \\
	&\leq& C' C'' \epsilon^{1+\delta} e^{n\kappa}|y|. \nonumber \label{eq boundforai}
	\eary
	Here in the last inequality we have used that $i \geq 2$ and $m \leq (1-\eta_0)n$. 
	By letting $\kappa$ be sufficiently small depending only on $g$ and $\eta_0$  we obtain the claim.
\end{proof}
By \eqref{eq psixn app} and Claim \ref{claim for ai}, we obtain the former part of Lemma \ref{lem psi to r}.
	
	To see the last statement, we assume that \eqref{eq psi small in uniform} holds, 
	and notice that by  \eqref{eq psixn app} and Claim \ref{claim for ai}
	\ary \label{eq rissmall}
	\norm{R}_{(-2,2)^2} \leq \norm{\Psi_{x}(\Lambda_n(x)^{-1}\cdot, \cdot)}_{(-2,2)^2} + C' \epsilon^{1+\delta} < 2 \epsilon/C'_3
	\eary
	for all sufficiently small $\epsilon > 0$.
	By $g \in \cU^\infty_{+, K, \kappa_1} \cap \cU^\infty_{-, K, \kappa_1}$, $R \in \cT_{x,n}$ and by letting $C'_3$ be sufficiently large, we have
	\aryst
	n \geq k^\epsilon(x).
	\earyst
    Moreover, by letting $C_3$ in Definition \ref{def hepsilon} be sufficiently large and by $R(\cdot, 0) \equiv 0$, we see that for any $y \in (-1,1)$ 
	\ary \label{eq rissmall 1}
    \norm{R(\cdot, y)}_{(-2,2)} \leq \frac{1}{2} \max(\epsilon|y|^{\delta_5/2}, C_3\epsilon|y|).
    \eary	
	For any $y \in (-1,1)$, by taking 
	\aryst
	\varphi(z) = R(z,y) + \sum_{i=2}^{K} a_i(y)z^i,
	\earyst
	 and taking $C_3$  and $C'_3$ large,
	we see that \eqref{eq def matchingtime 1} and \eqref{eq def matchingtime 2} follow from \eqref{eq psixn app}, \eqref{eq rissmall 1} and Claim \ref{claim for ai}. 
\end{proof}

	The following is the main lemma of this section.  Recall that the notion $C$-{\rm UNI} is introduced in Definition \ref{def tameanduni}.

\begin{lemma} \label{lem oscillation lower bound and smoothness upper bound} 
	For any sufficiently small $\kappa > 0$,  for any $C_1, n_0 > 0$, 
	$C_1$-{\rm UNI} holds on $\Omega(n_0,\kappa)$ at scales $\{ \Lambda^\epsilon \}_{\epsilon > 0}$.
\end{lemma}
\begin{proof}
	To verify $C_1$-{\rm UNI}, we take an arbitrary $x \in U$ with 
	\ary \label{eq takexnearomega}
	W^u_g(x, C_1\Lambda^\epsilon(x)^{-1}) \cap \Omega(n_0, \kappa) \neq \emptyset.
	\eary   
	By Definition \ref{def hepsilon}, there exists $\hat x \in W^s_g( x, 1)$ so that for 
	$n =  \varsigma^\epsilon(\hat x)$ we have 
	\aryst
	\Lambda^\epsilon( x) = \Lambda_n( x ) \sim \Lambda_n(\hat x).
	\earyst 
	By \eqref{eq local metric comparable},
	there exists a unique $\hat y \in (-2, 2)$ so that
	\ary \label{eq defofhaty}
	x = \Phi^s_{\hat x}(\hat  y).
	\eary
	Moreover, by  \eqref{eq takexnearomega} and by distortion estimate, we have that 
\aryst
\hat x \in \Omega(\epsilon,n_0 + C'_2, 2\kappa)
\earyst
 for some $C'_2 > 0$ depending on $g, \kappa$ and $C_1$. 

We denote $\hat x_n = g^n(\hat x)$, and for each integer $m$ we abbreviate $\lambda_{g^m(\hat x)}$ and $\mu_{g^m(\hat x)}$ as $\hat \lambda_m$ and $\hat \mu_m$ respectively.  Without loss of generality, let us assume that $J_0 = [0,1)$.

	Recall that by Definition \ref{def hepsilon}, for all sufficiently small $\epsilon$, we have
	\ary
	|\hat\lambda_0 \cdots \hat\lambda_{n-1}|, |\hat\mu_0 \cdots \hat\mu_{n-1}|^{-1} &<& \epsilon, \\
	\label{eq psixnupperbound}
	\norm{\Psi_{\hat x}(\Lambda_n(\hat x)^{-1}\cdot, y)}_{(-2,2)}  &<& 2C_3\epsilon, \quad \forall y \in (-1,1).
	\eary
	By  Lemma \ref{lem psi to r}, there exist $\varkappa \in \{\pm 1\}$ and
	\ary
	&& R(z,y) = h_1 y \varphi^{u,s}_{\hat x_n}(z) + h_2 z \varphi^{s,u}_{\hat x}(y) + Q(z,y) \in \cT_{\hat x,n} 
	\eary
	where 
	\ary
	   && Q \in Poly^{K,K} \quad \mbox{and} \quad	h_1, h_2 = O(\epsilon), \label{eq h1h2aresmall}
	\eary
	such that
	\ary \label{eq psixnrty}
\Psi_{\hat x}( \varkappa \Lambda_n(\hat x)^{-1} z, y) &=& R(z,y)  + O(\epsilon^{1+\delta}), \quad \forall z \in (-\rhone, \rhone).
	\eary
	
		Let $C'_3$ be given by Lemma \ref{lem psi to r}. 
	We have  that
	\aryst 
	\norm{\Psi_{\hat x}(\Lambda_{n-1}(\hat x)^{-1}\cdot, \cdot)}_{(-2,2)^2} > \epsilon/C'_3,
	\earyst
	for otherwise we would have $\varsigma^\epsilon(\hat x) \leq n-1$ by  Lemma \ref{lem psi to r}.
	Then by Lemma \ref{lem oscillation exists} and Lemma \ref{lem psi to r}, we see that
	\aryst 
	\norm{\Psi_{\hat x}(\Lambda_{n}(\hat x)^{-1}\cdot, \cdot)}_{(-2,2)^2} > \epsilon/C'_3C',
	\earyst
	Along with \eqref{eq psixnupperbound} and \eqref{eq psixnrty}, we see that 
	for all sufficiently small $\epsilon$
	\ary \label{eq rangeofr}
	3C_3\epsilon > \norm{R}_{(-2,2)^2} > \epsilon/C'_3C'.
	\eary
	Then there is $C_4 > 0$ depending only on $g$ and $C_3$ such that 
	\ary \label{eq tame in y 1}
	Q \in Poly^{K,K}_{C_4 \epsilon}.
	\eary

By Lemma \ref{lem oscillation exists} and \eqref{eq rangeofr}, there exists  $\bar y \in B(\hat y, \frac{\rmp}{4})$ such that
\aryst
\norm{R(\cdot, \bar y) - R(\cdot, \hat y)}_{(-2,2)}  > (C_5C'_3)^{-1}\epsilon
\earyst
for some $C_5 > 0$ depending only on $g, \varrho_2$.
By \eqref{eq h1h2aresmall}, \eqref{eq rangeofr} and \eqref{eq tame in y 1},
for any $y \in (-\varrho_1, \varrho_1)$ we have
\aryst
	\norm{R(\cdot, y) - R(\cdot, \bar y)}_{(-2,2)}  <  C' \epsilon |y - \bar y|^\theta + C' C_4 \epsilon |y - \bar y|.
\earyst
Thus for some $\varrho_3 \in (0, \varrho_2)$ depending only on $g, \varrho_2, C_4, C'_3$, for any $y \in B(\bar y, \varrho_3)$ we have
	\ary \label{eq differenceofr}
	\norm{R(\cdot, y) - R(\cdot, \hat y)}_{(-2,2)}  > (2C_5C'_3)^{-1}\epsilon.
	\eary
	Notice that by \eqref{eq local metric comparable} and \eqref{eq defofhaty}, we have $\Phi^s_{\hat x}(\bar y) \in W^s_g(x, \varrho_2)$.

		By \eqref{eq psixnrty}, \eqref{eq differenceofr} and Lemma \ref{lem oscillation exists} and by distortion estimate, for all sufficiently small $\epsilon$ we have for any $y \in B(\bar y,  \rhothree)$ that 
\aryst
&& \norm{\Psi_{\hat x}(  \Lambda_n(x)^{-1}\cdot, y)  - \Psi_{\hat x}(  \Lambda_n(x)^{-1}\cdot, \hat y)}_{(-2,2)} \\
&\geq& \norm{\Psi_{\hat x}(  \Lambda_n(\hat x)^{-1}\cdot, y)  - \Psi_{\hat x}(  \Lambda_n(\hat x)^{-1}\cdot, \hat y)}_{(-(C')^{-1}, (C')^{-1})} \\
&\geq& \norm{R(\cdot, y) - R(\cdot, \hat y)}_{(-(C')^{-1}, (C')^{-1})} - O(\epsilon^{1+\delta}) \\
&\geq& (C'C_5C'_3)^{-1}\epsilon.
\earyst
	Our claim then follows from Lemma \ref{lem shift invariance of temporal function}.
\end{proof}

\begin{comment}
The following corollary of Lemma \ref{lem psi to r} will also be used at the end of this section. 
	
	\begin{cor} \label{cor psi to r}
		For $\kappa_1, n_0$ in Lemma \ref{lem psi to r}, 
		for any $x \in \Omega_-(n_0, \kappa)$, for all sufficiently small $\epsilon > 0$, 
		for any integer $n \geq 1$ such that $|\lambda_0 \cdots, \lambda_{n-1}|$, $|\mu_0^{-1} \cdots \mu_{n-1}^{-1}| < \epsilon$, for any $\varrho_4 > 0$, there exists $\epsilon_6 > 0$ depending only on $g, K, L, \varrho_4$  such that
		\aryst
		\norm{\Psi_{x}(\Lambda_{n}(x)^{-1}\cdot, \cdot)}_{(-\rhone, \rhone)^2} < 
		C(\varrho_4)^{-1}\norm{\Psi_{x}(\Lambda_{n}(x)^{-1}\cdot, \cdot)}_{(-\varrho_4, \varrho_4)^2} + O(\epsilon^{1+\delta}).  
		\earyst
\end{cor}

\end{comment}

\subsection{ Uniform non-integrability in Class {\it II}}

For $g$ in Class {\it II}, we have the following analogous statement.
\begin{lemma}\label{lem oscillation lower bound and smoothness upper bound 3}
	 Define $\{ \Lambda^\epsilon \}_{\epsilon > 0}$ by Definition \ref{def hepsilon 1}. Then for any $C_1 > 0$,
	$C_1$-{\rm UNI} holds on $U$ at scales $\{ \Lambda^\epsilon \}_{\epsilon > 0}$.
\end{lemma}
\begin{proof}
	This is an immediate consequence of Lemma \ref{lem varsigmaslowgrowing case ii}.
\end{proof}

\subsection{$\{\Lambda^\epsilon\}_{\epsilon > 0}$ is adapted to $\Omega$}

We will only detail the case where $g$ is in Class $I_F$. The other case follows easily from distortion estimate.

\begin{lemma}\label{lem verifyadapted}
	 For any $\kappa > 0$, for any sufficiently large integer $n_1 > 0$, for any integer $n > 0$,
	the sequence $\{ \Lambda^\epsilon \}_{\epsilon > 0}$, defined in Definition \ref{def hepsilon}, is $n_1$-adapted to 	$\Omega( n, \kappa)$.
\end{lemma}
\begin{proof}
	Take  $x \in \Omega(n,\kappa)$ and an inverse branch $v \in \sigma^{-n_1}_{x}$.
	We denote $x_* = v(x)$.
	We take an arbitrary $\hat x \in U$ such that 
	\aryst
	y \in W^u_g(x_*, 4\Lambda^\epsilon(y)^{-1}). 
	\earyst
	In the following, we will show that $\Lambda^\epsilon(x) < C \Lambda^\epsilon(y)$ for some $C > 0$ independent of $x,y$ and $\epsilon$.
	
	By construction, we have 
	\aryst
	x' = g^t(x_*) \in W^s_{g}(x, 1) \quad \mbox{where} \quad t = \tau_{n_1}(x_*).
	\earyst
	Denote
	\aryst
	m = \lfloor t \rfloor \quad \mbox{and} \quad
	 \tilde x = g^m(x_*).
	\earyst
	By distortion estimate, there is an integer $n_0 > 0$ depending only on $g, n, \kappa$ and $n_1$ such that $x_*, \tilde x, x' \in \Omega( n_0, 2\kappa)$.

	By definition, there exists $\bar x \in W^s_g(y, 1)$ such that 
	\ary \label{def of n0andbarx}
	\varsigma^\epsilon(\bar x) = N_0 : =  \sup_{x'' \in W^s_g(y, 1)} \varsigma^\epsilon(x'').
	\eary
	By  distortion estimate, we see that 
	\ary
	N_0 \geq k^\epsilon(\bar x) &\geq& \max(k^\epsilon(x_*), k^\epsilon(\tilde x), k^\epsilon(x')) - t -  C'.  \nonumber
	\eary
Then by Lemma \ref{lem psi to r}, there exist $R \in \cT_{x_*,N_0}$, $R' \in \cT_{x',N_0}$ and $\tilde R \in \cT_{\tilde x,N_0}$ such that
	\ary
	\norm{\Psi_{x_*}(\Lambda_{N_0}(x_*)^{-1}\cdot, \cdot) - R}_{(-\rhone, \rhone)^2} &=& O(\epsilon^{1+\delta}), \label{eq r} \\
\norm{\Psi_{x'}(\Lambda_{N_0}(x')^{-1}\cdot, \cdot) - R'}_{(-\rhone, \rhone)^2} &=& O(\epsilon^{1+\delta}), \label{eq r'} \\
    \norm{\Psi_{\tilde x}(\Lambda_{N_0}(\tilde x)^{-1}\cdot, \cdot) - \tilde R}_{(-\rhone, \rhone)^2} &=& O(\epsilon^{1+\delta}). \label{eq tilde r}
	\eary

	By \eqref{eq def matchingtime 1} and \eqref{eq def matchingtime 2} in Definition \ref{def hepsilon} and by \eqref{def of n0andbarx}, we have
	\ary \label{lem psixupperbound1}
	\norm{\Psi_{\bar x}(\Lambda_{N_0}(\bar x)^{-1}\cdot, \cdot)}_{(-1,1)^2} < 2C_3\epsilon.
	\eary
	By $\bar x \in W^s_g(y, 1)$ and by distortion estimate, we have
	\ary \label{eq comparetwolambda1}
	\Lambda^\epsilon(y) =
 \Lambda_{N_0}(y) \sim \Lambda_{N_0}(\bar x),
	\eary
	and there is $\varrho > 0$ depending only on $g$ such that, by letting $\varrho_1$ be large depending only on $g$, we have
	\ary \label{eq comparetwopsi1}
	\norm{\Psi_{y}(\Lambda_{N_0}(y)^{-1}\cdot, \cdot)}_{(-\varrho, \varrho)^2} \leq \norm{\Psi_{\bar x}(\Lambda_{N_0}(\bar x)^{-1}\cdot, \cdot)}_{(-\varrho_1, \varrho_1)^2}.
	\eary
	Similarly, by $y \in W^u_g(x, 4\Lambda^\epsilon(y)^{-1})$ and by distortion estimate, we have
	\ary \label{eq comparetwolambda2}
	\Lambda_{N_0}(x_*) \sim \Lambda_{N_0}(y).
	\eary
	Then by making $\varrho$ smaller and making $\varrho_1$ larger if necessary, both depending only on $g$, we have
    \ary \label{eq comparetwopsi2}
    \norm{\Psi_{x_*}(\Lambda_{N_0}(x_*)^{-1}\cdot, \cdot)}_{(-\varrho, \varrho)^2} \leq \norm{\Psi_{y}(\Lambda_{N_0}(y)^{-1}\cdot, \cdot)}_{(-\varrho_1, \varrho_1)^2}
    \eary

	By \eqref{eq r}-\eqref{eq tilde r}, \eqref{eq comparetwolambda1} and Lemma \ref{lem oscillation exists}, 
	we can compare, up to error $O(\epsilon^{1+\delta})$, the following: (1) LHS of \eqref{eq comparetwopsi1} and RHS of \eqref{eq comparetwopsi2}; (2) LHS of \eqref{lem psixupperbound1} and RHS of \eqref{eq comparetwopsi1}; (3) LHS of \eqref{eq comparetwopsi2} and $\norm{\Psi_{x_*}(\Lambda_{N_0}(x_*)^{-1}\cdot, \cdot)}_{(-\rhone, \rhone)^2}$. Then by \eqref{lem psixupperbound1}, \eqref{eq comparetwopsi2} and \eqref{eq comparetwopsi1}, we obtain
	\ary \label{eq boundingpsixlambdan0}
	\norm{\Psi_{x_*}(\Lambda_{N_0}(x_*)^{-1}\cdot, \cdot)}_{(-\rhone, \rhone)^2}  \leq C'C_3 \epsilon.
	\eary

	By Lemma \ref{lem shift invariance of temporal function}(1) and by letting $n_1$ be sufficiently depending only on $g$, there exist $C > 1$ such that 
	\aryst
	 \norm{\Psi_{x'}(\Lambda_{N_0}(x')^{-1}\cdot, C^{-1}\cdot)}_{(-\rhone, \rhone)^2} 
	< \norm{\Psi_{\tilde x}(\Lambda_{N_0 - m}(\tilde x)^{-1}\cdot, \cdot)}_{(-\rhone, \rhone)^2}.
	\earyst
	Then by \eqref{eq r'} and Lemma \ref{lem oscillation exists}, we have
	\ary \label{eq boundingpsixi'} \\
	\nonumber
	\norm{\Psi_{x'}(\Lambda_{N_0}(x')^{-1}\cdot, \cdot)}_{(-\rhone, \rhone)^2} \leq C' \norm{\Psi_{\tilde x}(\Lambda_{N_0 - m}(\tilde x)^{-1}\cdot, \cdot)}_{(-\rhone, \rhone)^2} + O(\epsilon^{1+\delta}).
	\eary
	By \eqref{lem finvariance}, \eqref{eq tilde r}, \eqref{eq boundingpsixlambdan0} and Lemma \ref{lem oscillation exists},  there exists $C'_{n_1} > 0$ depending only on $g$ and $n_1$ such that
	\ary
	&& \norm{\Psi_{\tilde x}(\Lambda_{N_0 - m}(\tilde x)^{-1}\cdot, \cdot)}_{(-\rhone, \rhone)^2}  \label{eq psin2minusn} \\
	&<& C'C'_{n_1} \norm{\Psi_{\tilde x}(\Lambda_{N_0 - m}(\tilde x)^{-1}\cdot, \mu_{0}^{-1} \cdots \mu_{m -1}^{-1}\cdot)}_{(-\rhone, \rhone)^2} + O(\epsilon^{1+\delta})  \nonumber  \\
	&<& 
	C'C'_{n_1} \norm{\Psi_{x_*}(\Lambda_{N_0 }(x_*)^{-1}\cdot, \cdot)}_{(-\rhone, \rhone)^2}  + O(\epsilon^{1+\delta})  \nonumber \\
	&<& C' C'_{n_1} C_3 \epsilon.  \nonumber
	\eary
	
	By letting $p$ be a large integer depending only on $g$ and $n_1$, we have
	\ary \label{eq pislarge1}
	N_1 = N_0 + p > \sup_{x'' \in  W^s_g(x, 1)} k^\epsilon(x'').
	\eary
	By letting $\rhone$ be large depending only on $g$, we have
	\ary \label{eq upperboundpsix''}
    \sup_{x'' \in W^s_g(x, 1)} \norm{\Psi_{x''}(\Lambda_{N_1}(x'')^{-1}\cdot, \cdot )}_{(-1,1)^2} 
	\leq	\norm{\Psi_{x'}(\Lambda_{N_1}(x')^{-1}\cdot, \cdot )}_{(-\rhone, \rhone)^2}.
	\eary
	
	By \eqref{eq pislarge1}, we have $\Lambda_{N_1}(x')^{-1} < e^{-\chi_0p} \Lambda_{N_0}(x')^{-1}$. Then by \eqref{eq r'} and Lemma \ref{lem lower bound by uni for templates}, we have
	\aryst
	\norm{\Psi_{x'}(\Lambda_{N_1}(x')^{-1}\cdot, \cdot )}_{(-\rhone, \rhone)^2}
	\leq c_1(e^{-\chi_0 p}) \norm{\Psi_{x'}(\Lambda_{N_0 }(x')^{-1}\cdot, \cdot )}_{(-\rhone, \rhone)^2} + O(\epsilon^{1+\delta}).
	\earyst
	By \eqref{eq boundingpsixi'} and \eqref{eq psin2minusn}, we have
	\aryst
	\norm{\Psi_{x'}(\Lambda_{N_1}(x')^{-1}\cdot, \cdot )}_{(-\rhone, \rhone)^2}
	\leq c_1(e^{-\chi_0 p }) C'_{n_1} C' C_3\epsilon + O(\epsilon^{1+\delta}).
	\earyst
	We let $C'_3$ be as in  Lemma \ref{lem psi to r}. 
	By letting $p$ be sufficiently large depending only on $g, n_1, C'_3$ and $C_3$, and by letting $\epsilon$ be sufficiently small, we have
	\aryst
	\norm{\Psi_{x'}(\Lambda_{N_1}(x')^{-1}\cdot, \cdot )}_{(-\rhone, \rhone)^2} &<& 
	\epsilon/C'_3.
	\earyst
	By \eqref{eq pislarge1}, \eqref{eq upperboundpsix''} and Lemma \ref{lem psi to r}, 
	this implies that 
	\aryst
	N_1 \geq \sup_{x'' \in  W^s_g(x, 1)} \varsigma^\epsilon(x'').
	\earyst
	Hence $\Lambda_{N_1}(x) \geq \Lambda^\epsilon(x)/C'$. 
	On the other hand, it is clear that there exists $C_{n_1} > 0$ depending only on $g$ and $n_1$ such that
	\aryst
	\Lambda_{N_0}(x_*) \geq C_{n_1}^{-1} \Lambda_{N_1}(x).
	\earyst
    The proof of Lemma \ref{lem verifyadapted} then follows from \eqref{eq comparetwolambda1} and \eqref{eq comparetwolambda2}.
\end{proof}

\section{The recurrence property} \label{sec recurrence}

In this section, we fix a H\"older potential $F$ and  assume that $g$ is in Class $I_F$.
 
\begin{lemma} \label{lem fast recurrence for expanding set}
	For any  $\kappa > 0$, for any integer $n_1 > 0$, for all sufficiently large integer $n > 0$, the set $ \Omega(n, \kappa)$ (defined in Definition \ref{def uniform set}) is $n_1$-recurrent with respect to $\nu_U$.
\end{lemma}
\begin{proof}
	We will abbreviate $\Omega(n, \kappa)$ as $\Omega$.
	Define for any integer $L > 0$ and $\eta \in (0,1)$
	\aryst
		B_{L, \eta} &=& \{  x \in U \mid  |\{ 1 \leq j \leq L \mid \sigma^{jn_1}(x) \in \Omega \}| < \eta L  \}.
	\earyst		
	We have
	\aryst
	\inf_{y \in W^s_{g, loc}(x)} \det Dg^{\tau(x)}(y) > e^{-C'}.
	\earyst
	By \eqref{eq wsgx1 is long}, and
	by letting $n$ be large depending on $g, \kappa$ and $n_1$,  we can see that for each $x \notin \Omega$ 
	there exists $m \geq 1$ such that
	\aryst
	\inf_{y \in W^s_{g, loc}(x)}
	\det Dg^{\tau_{mn_1}(y)}(y) > e^{mn_1 \kappa/2}.
	\earyst 
	Then it is straightforward to see that for any sufficiently small $\eta > 0$, for any sufficiently large integer $L > 0$, for any $x \in B_{L, \eta}$, there exists an integer  $m > L$ such that 
	\aryst
	\det Dg^{\tau_{mn_1}(x)}(x) > e^{(1-\eta)m n_1 \kappa/2 - \eta Ln_1 C'} > e^{m n_1 \kappa/4}.
	\earyst
	Then we conclude the proof by Lemma \ref{lem subexponential growth on average}.
\end{proof}

\section{Proof of Proposition \ref{prop main 1}} \label{sec checking the conditions}

In this section, we use the results we obtained from Section \ref{sec construct partition} to Section \ref{sec recurrence} to deduce Proposition \ref{prop main 1}.

\noindent{$\bullet$} Assume that $g$ is in Class $I_F$.
We choose $K$, $C_3, \delta_5$ so that Lemma \ref{lem varsigmaslowgrowing} is applicable.
We define $\{ \Lambda^\epsilon \}_{\epsilon > 0}$  by Definition \ref{def lambda epsilon} for $C_3,\delta_5$.
By Corollary \ref{cor lambdaepsilonstable}, the sequence of functions  $\{ \Lambda^\epsilon \}_{\epsilon > 0}$ is stable and tame.

We let $\kappa$ satisfy the condition of Lemma \ref{lem oscillation lower bound and smoothness upper bound}, and let $n_1 > 0$ be a sufficiently large integer so that Lemma \ref{lem verifyadapted} is applicable. 
 By Lemma \ref{lem fast recurrence for expanding set}, the set $\Omega = \Omega( n, \kappa)$ is $n_1$-recurrent with respect to $\nu_U$ for some integer $n$ sufficiently large depending only on $g,\kappa$ and $n_1$.
 By Lemma \ref{lem verifyadapted}, we can see that the sequence of functions $\{ \Lambda^\epsilon \}_{\epsilon > 0}$ is $n_1$-adapted to $\Omega$.
This verifies Proposition \ref{prop main 1}(1),(2).

The conclusion of Proposition \ref{prop main 1} then follows from
 Lemma \ref{lem oscillation lower bound and smoothness upper bound}.% and Lemma \ref{lem omegaisuniform}.

\smallskip

\noindent{$\bullet$} Assume that $g$ is in Class {\it II}. 
For all sufficiently small $\epsilon > 0$, we define $\vartheta^\epsilon$ and $\{ \Lambda^\epsilon \}_{\epsilon > 0}$  by Definition \ref{def hepsilon 1}. 
 By Corollary \ref{lemma lambdaepsilonstable}, the sequence of functions  $\{ \Lambda^\epsilon \}_{\epsilon > 0}$ is stable and tame.
   We set $\Omega = U$. The conclusion of Proposition \ref{prop main 1} then follows by combining Corollary \ref{lemma lambdaepsilonstable} and Lemma \ref{lem oscillation lower bound and smoothness upper bound 3}.

\section{Proof of Proposition \ref{prop main 2}} \label{sec proofofpropmaincriterion}

In this section, we will give the proof of Proposition \ref{prop main 2}  based on several lemmas whose proofs compass Section \ref{sec partition} to \ref{sec proof of lemma prop exp rec in vol exp}.

\subsection{An inductive scheme for the decay} \label{sec: An inductive scheme for the decay}

We give a variant of the inductive scheme of Dolgopyat (see \cite[Lemma 10'']{Dol}) for controlling the iterates of $\widetilde\cL_{a,b}$ using {\it majorant sequence} and certain recurrence estimate.

Let $\theta \in (0,1)$ be as in Proposition \ref{prop main 2}.
In the course of the proof, we may reduce $\theta$ if necessary, but will only do this finitely many times. This will not affect the generality of the result.

\begin{prop} \label{prop L2 iteration}
	Under the hypothesis of Proposition \ref{prop main 1} for some  sufficiently large $C_1 > 1$,
	there exist $C_8 > 1$, $\kappa_3,  \kappa_4, \delta_1  > 0$ and an integer $n_1 > 0$ such that the following is true.
	For  any $a \in \R$ with $|a|$ sufficiently small, for any $b \in \R$ with $|b|$ sufficiently large,  define $f^{(a,b)}$ as in Section \ref{sec smoothing} for $\delta_1$, and denote
	\ary \label{def clwclcm}
	\epsilon = |b|^{-1}, \quad \widetilde\cL = \widetilde\cL_{a,b},\quad \cM = \cM_{a,b} \quad \mbox{(see \eqref{eqmabdef} and \eqref{eq deflab}).}
	\eary
	Then
	for any $u \in C^\theta(U)$,  there is a sequence of functions $\{ H_n \}_{0 \leq n \leq  \lfloor  \ln |b|  \rfloor}$ in $C^0(U, \R_+)$ such that $H_0 \leq  \norm{u}_{\theta,b}$, and
	\enmt
	\item for any $0 \leq n \leq \ln |b|$  we have
	\aryst
	|\widetilde\cL^{nn_1 + \lfloor C_8 \ln |b| \rfloor}u(x)| \leq H_n(x), \quad \forall x \in U;
	\earyst
	\item for any $1 \leq n \leq \ln |b|$ there is a 
	%$\cP^+_0$-measurable, $\cP^+_1$-dense 
	subset $\Omega_n \subset U$	
	such that
	\aryst
	H^2_{n}(x) \leq \begin{cases}
		(1 - \kappa_4)\cM^{n_1} H^2_{n-1}(x), & \mbox{if } x \in \Omega_n,	\\
		 \clb \cM^{n_1} H^2_{n-1}(x), & \mbox{otherwise};
	\end{cases}
	\earyst
	\item for any $\frac{1}{2}\ln|b| \leq  n 	\leq  \ln|b|$, we have
	\aryst
	\nu_U(\{ x \in U  \mid  |\{1 \leq j \leq n \mid \sigma^{jn_1 + \lfloor C_8 \ln |b| \rfloor}(x) \in \Omega_j \}| < \kappa_3 n  \}) < e^{- n \kappa_3}.
	\earyst 
	\eenmt
\end{prop}

We are ready to deduce Proposition \ref{prop main 2} from Proposition \ref{prop L2 iteration}.
\begin{proof}[Proof of Proposition \ref{prop main 2}]
	Let $C_1 > 1$ be sufficiently large so that Proposition \ref{prop L2 iteration} is applicable.
	We let 
	$n_1, \kappa_3, \kappa_4, \delta_1 > 0$ be given by Proposition \ref{prop L2 iteration}. We define $f^{(a,b)}$, $\widetilde\cL$, $\cM$ as in Proposition \ref{prop L2 iteration}, and abbreviate $\cL = \cL_{a,b}$.
	
	Take an arbitrary integer $L \in ( \frac{1}{2}\ln|b|,  \ln |b|)$. We define a $U$-valued random process $X$ on the space $(U, \nu_U)$ 
	by  $\{ X_{m}(x) = \sigma^{m n_1 + \lfloor C_8 \ln |b| \rfloor}(x) \}_{0 \leq m \leq L}$  where $x$ has distribution $\nu_U$. 
	We see that for any $1 \leq m \leq L$, $X_m = \sigma^{n_1}(X_{m-1})$; and the marginal distribution of $X_m$ equals $\nu_U$ since $\nu_U$ is $\sigma$-invariant.
	Then by the Gibbs property, for any $0 \leq m \leq L-1$, we have
	\ary \label{eq xmtoxm+1}
\PP(X_m = y \mid X_{m+1} = z) = \begin{cases}
	0, & z \neq \sigma^{n_1}(y), \\
	e^{\hat f_{n_1}(y)}, & z = \sigma^{n_1}(y).
\end{cases}
	\eary
	Let $\kappa_* > 0$ be a small constant so that
	\ary \label{eq kappastarsmall}
	1 + \kappa_* < (\frac{1-\kappa_4/4}{1-\kappa_4/2})^{\kappa_3}.
	\eary
	We define a $\R$-valued random process $G = \{G_m\}_{0 \leq m \leq L}$ on the same probability space of $X$ 
	by setting
	\aryst
	G_0(x) = H^2_0(x), \quad 
	G_{m+1}(x) = \begin{cases}
		(1-\kappa_4/2)G_m(x),   & \mbox{if } X_{m+1} \in \Omega_{m+1},  \\
		(1+ \kappa_*)G_m(x),   & \mbox{otherwise, }
	\end{cases}  \forall 0 \leq m \leq L-1.
	\earyst

	\begin{claim}
		For any $0 \leq m \leq L$ we have
		\ary \label{eq e f m x m h2m x m}
		\mathbb E( G_m \mid X_{m}) \geq H^2_m(X_m).
		\eary
	\end{claim}
\begin{proof}
	It is direct to see that the claim holds for $m=0$.	Now assume that we have shown the claim for some $0 \leq m < L$.
	
	For any $z \in U$, we have
	\aryst
	\mathbb E(G_{m+1} \mid X_{m+1} = z) &=& \sum_{y \in \sigma^{-n_1}(z)} \mathbb E(G_{m+1} \mid X_{m+1} = z, X_{m} = y) \PP(X_m = y \mid X_{m+1} = z).
	\earyst
	We first assume that $z \in \Omega_{m+1}$.  Notice that for all $a,b \in \R$ with $|a|$ sufficiently small and $|b|$ sufficiently large  we have that
	\aryst
	\mathbb E(G_{m+1} \mid X_{m+1} = z) &=& (1-\kappa_4/2)\sum_{y \in \sigma^{-n_1}(z)} \mathbb E(G_m \mid X_{m+1} = z, X_{m} = y) \PP(X_m = y \mid X_{m+1} = z) \\
\mbox{(by \eqref{eq xmtoxm+1})}	&=&(1-\kappa_4/2)\sum_{y \in \sigma^{-n_1}(z)} \mathbb E(G_m \mid X_{m} = y) e^{\hat f_{n_1}(y)} \\
\mbox{(by induction)}	&\geq&(1-\kappa_4/2) \sum_{y \in \sigma^{-n_1}(z)} H^2_m(y) e^{\hat f_{n_1}(y)} \\
\mbox{(by \eqref{eq hatftofab})} 	&\geq&(1-\kappa_4) \sum_{y \in \sigma^{-n_1}(z)} H^2_m(y) e^{f^{(a,b)}_{n_1}(y)} \\
	&=& (1-\kappa_4) (\cM^{n_1} H^2_m)(z)  = H_{m+1}^2(z).
	\earyst
	The last equality follows from Proposition \ref{prop L2 iteration}(2).
	
	Now assume that $x \notin \Omega_{n+1}$. By a similar argument as above, we deduce that
	\aryst
\mathbb E(G_{m+1} \mid X_{m+1} = z)	&\geq& (1 + \kappa_*) \sum_{y \in \sigma^{-n_1}(z)} H^2_m(y) e^{\hat f_{n_1}(y)} \\ 	&\geq& 
(\cM^{n_1} H^2_m)(z)   = H_{m+1}^2(z).
\earyst	
The last equality follows from  Proposition \ref{prop L2 iteration}(2), and the last inequality follows from \eqref{eq hatftofab}, since we have
\aryst
(1+\kappa_*) e^{\hat f_{n_1}(y)} \geq e^{f^{(a,b)}_{n_1}(y)}, \quad \forall y \in U
\earyst
whenever $|a|$ is sufficiently small and $|b|$ is sufficiently large.
\end{proof}
	By definition, for any $x \in U$ we have
	\ary \label{eq h n less h 0}
	G_{L}(x) \leq (1+\kappa_*)^L G_0(x).
	\eary
	By Proposition \ref{prop L2 iteration}(3), and by letting $\epsilon$ be sufficiently small, we see that
	\aryst
	\nu_U(\{ x \in U  \mid  |\{1 \leq j \leq L \mid X_j(x) \in \Omega_j \}| < \kappa_3 L  \}) < e^{- L \kappa_3}.
	\earyst
	On the other hand, by \eqref{eq kappastarsmall}, for any $x \in U$ satisfying
	\aryst
	|\{1 \leq j \leq L \mid X_j(x) \in \Omega_j \}| \geq \kappa_3 L, 
	\earyst
	we have 
		\ary \label{eq h n strictly less h 0}
	G_L(x) \leq (1-\kappa_4/2)^{\kappa_3 L}(1+\kappa_*)^L G_0(x) \leq (1-\kappa_4/4)^{\kappa_3 L} G_0(x).
	\eary
	Let us denote $m = \lfloor C_8 \ln |b| \rfloor$. Then
	 by \eqref{eq c0differenceLandtildeL}, Proposition \ref{prop L2 iteration}(1),(3),  \eqref{eq e f m x m h2m x m}, \eqref{eq h n less h 0} and \eqref{eq h n strictly less h 0}, for any $u \in C^\theta(U)$ we have
	\aryst
	\norm{\cL^{L n_1 + m} u}^2_{L^2(U, d\nu_U)} &<& 2\norm{\widetilde\cL^{L n_1 + m}u -  \cL^{Ln_1 + m}u}^2_{C^0} + 2\norm{\widetilde\cL^{L n_1 + m}u}^2_{L^2(U, d\nu_U)}  \\
	&<&2C'Ln_1 |b|^{-\delta_1\theta/4}\norm{u}^2_{C^0} + 2\norm{H_L}^2_{L^2(U, d\nu_U)} \\
	&<&2C'Ln_1 |b|^{-\delta_1\theta/4}\norm{u}^2_{C^0} + 2\norm{G_L}_{L^1(U, d\nu_U)} \\
	&\leq& 2C'Ln_1 |b|^{-\delta_1\theta/4}\norm{u}^2_{C^0} +  2( (1- \kappa_4/4)^{L\kappa_3} + (1+\kappa_*)^Le^{-L\kappa_3}) \norm{G_0}_{C^0} \\
	&\leq& 2C'Ln_1 |b|^{-\delta_1\theta/4}\norm{u}^2_{C^0} +  2( (1- \kappa_4/4)^{L\kappa_3} + (1+\kappa_*)^Le^{-L\kappa_3}) \norm{H_0}^2_{C^0} \\
	&\leq&[2C'Ln_1 |b|^{-\delta_1\theta/4} +  2( (1- \kappa_4/4)^{L\kappa_3} + (1+\kappa_*)^Le^{-L\kappa_3})  ] \norm{u}^2_{\theta, b}.  
	\earyst
	Then it is clear that we have \eqref{eq decay tol2fromtheta} for $n =  C\ln |b|$ where $C$ is sufficiently large depending on $\kappa_3$ and $\kappa_4$. It is  an exercise to conclude the proof for all sufficiently large $n$.
\end{proof}

\subsection{Cancellation on uniform set}
 Recall that for any function $\Lambda: U \to \R_+$ and any $x \in U$, $J^\Lambda_x, \Phi^\Lambda_x$ are defined by \eqref{def jx} and \eqref{def wx} respectively. 

\begin{defi}\label{def coneK}
	We denote by ${\cal K}_\Lambda$
	the set of $C^1$ functions $h: U \to \R_+$ such that for any $x \in U$, any $z \in J^\Lambda_x$, we have
	\ary  \label{def property of the cone} 
	| \partial_z  [\log h \circ \Phi^\Lambda_x](z)| \leq 1.
	\eary
\end{defi}
 Consequently, for every $h \in \cal K_\Lambda$ we have
\ary
\label{eq ratio h}
(C')^{-1} < \frac{h(\Phi^\Lambda_x(z))}{h(x)}  < C', \quad \forall z \in J^\Lambda_x.
\eary

\begin{defi} \label{def hlambdaehc}
	Let $\Lambda, H, G: U \to \R_+$ be three arbitrary functions.
	We denote by  $\cH^\Lambda_{H, G}$ the set of functions $u \in C^\theta(U)$ such that 
	for any $x \in U$, there exists $u_x \in C^\theta(J^\Lambda_x)$ satisfying
	\aryst
	|u_x|_\theta &\leq& G(x), \\
	|u (\Phi^\Lambda_x(s)) - u_x(s)| &\leq& H( \Phi^\Lambda_x(s)), \quad \forall s \in J^\Lambda_x.
	\earyst
\end{defi}

We now state the inductive step for verifying Proposition \ref{prop L2 iteration}(1),(2). 
\begin{lemma} \label{lem dim cancellation}
	Under the hypothesis of Proposition \ref{prop main 1} for some sufficiently large $C_1 > 1$, 
	there exist $C_8 > 1$, $\delta_1, \eta_1 > 0$ such that for all sufficiently large integer $n_1 > 0$, for all sufficiently small $\kappa_5 > 0$, for any $a \in \R$ with $|a|$ sufficiently small, and any $b \in \R$ with $|b|$ sufficiently large, the following is true. Define $f^{(a,b)}$ as in Section \ref{sec smoothing} for $\delta_1$, define $\epsilon = |b|^{-1}$, $\widetilde \cL$, $\cM$ as in \eqref{def clwclcm}, and define $\Omega$, $\{\Lambda^\epsilon\}_{\epsilon > 0}$ by Proposition \ref{prop main 1} for $C_1, n_1$.
	Then for any $u \in C^\theta(U)$
	\enmt
	\item 
	we have
	\aryst
	u_0 := \widetilde \cL^{\lfloor C_8 \ln |b| \rfloor} u \in \cH^{\Lambda^\epsilon}_{\epsilon^{\eta_1} \norm{u}_{C^0},  C_8\norm{u}_{\theta, b}};
	\earyst
	\item
	there is a $\nu_U$-measurable finite partition $\cP_0$ of $U$ with 
	$\cP_1 = \sigma^{-n_1}(\cP_0)  \succeq  \cP_0$  such that
	  for any integer $n \geq 0$, for any $H_n \in \cK_{ \Lambda^\epsilon}$ satisfying
	\ary \label{eq nn3uniform bound}
	|\widetilde \cL^{nn_1}u_0(x)| &\leq& H_n(x) \leq \norm{u}_{\theta,b}, \quad \forall x \in U, \\
	\label{eq nn3holder bound}
	\widetilde \cL^{nn_1}u_0 &\in& \cH^{\Lambda^\epsilon}_{(n+1)\epsilon^{\eta_1}\norm{u}_{\theta, b},  C_8H_n}, \\
	\label{eq lowerboundofhn}
	H_n &\geq& (n+1) \epsilon^{\eta_1/2} \norm{u}_{\theta,b},
	\eary
	there exists a function $P_n \in C^1(U, [1 - \kappa_5,1]) \cap \cK_{ \Lambda^\epsilon}$ 
	such that 
	\enmt
	\item for every $x \in \Omega$, there exist $y \in \cP_0(x)$ and an inverse branch $v \in \sigma^{-n_1}_y$ such that we have
	\ary \label{eq valueofpnonsomesmallinterval}
	P_{n}(z) = 1 - \kappa_5, \quad \forall z \in v(\cP_1(y));
	\eary
	\item the function $H_{n+1} = \cM^{n_1}(P_n H_n)$ satisfies
	\ary
	\label{eq n+1n3incone}
	H_{n+1} &\in& \cK_{\Lambda^\epsilon}, \\
	\label{eq n+1n3uniform bound}
	|\widetilde \cL^{(n+1)n_1}u_0(x)| &\leq& H_{n+1}(x) \leq \norm{u}_{\theta,b}, \quad \forall x \in U, \\
	\label{eq n+1n3holder bound}
	\widetilde \cL^{(n+1)n_1}u_0 &\in& \cH^{\Lambda^\epsilon}_{(n+2)\epsilon^{\eta_1}\norm{u}_{\theta, b}, C_8 H_{n+1}}.
	\eary
	\eenmt
	\eenmt
\end{lemma}
The proof of Lemma \ref{lem dim cancellation} is deferred to Section \ref{sec proof lem dim cancellation} after some preparations in Section \ref{sec partition} and \ref{sec: bounds for smoothness}.

\begin{defi}
	Given a subset $A \subset U$ and a $\nu_U$-measurable finite  partition $\cP$ of $U$. We say that $A$ is $\cP$-{\em measurable} if $A$ is measurable modulo $\nu_U$-null subsets with respect to the $\sigma$-algebra generated by $\cP$.
	Given a $\cP$-measurable set $B \subset U$,
	we say that a subset $A \subset B$ is $\cP$-{\em dense} if for $\nu_U$-almost every $x \in B$, we have $\nu_U(\cP(x) \cap A) > 0$. Here $\cP(x)$ denotes the atom of $\cP$ containing $x$. 
\end{defi}

As an immediate corollary of Lemma \ref{lem dim cancellation}, we have the following.

\begin{cor}\label{cor construct omega n+1}
	For any $\kappa_5 > 0$ and any integer $n_1 > 0$, there is $\kappa_4 > 0$ depending only on $g$, $F$, $\kappa_5$ and $n_1$ such that the following is true.
	For $H_n$, $H_{n+1}$, $\cP_0$, $\cP_1$, $\Omega$ in Lemma \ref{lem dim cancellation} (for $\kappa_5$ and $n_1$), 
	denote by $\widetilde \Omega$ the minimal (with respect to inclusion) $\cP_0$-measurable subset of $U$ containing $\Omega$. Then
	there is a $\cP_1$-measurable and $\cP_0$-dense subset $\Omega_{n+1}$ of $\widetilde \Omega$ such that 
	\aryst
	H^2_{n+1}(x) \leq \begin{cases}
		(1 - \kappa_4)\cM^{n_1} H^2_{n}(x), & \mbox{if } x \in \Omega_{n+1},	\\
		\cM^{n_1} H^2_{n}(x), & \mbox{otherwise}.
	\end{cases} 
	\earyst	
\end{cor}
\begin{proof}
	Let $P_n$ be given by Lemma \ref{lem dim cancellation}.
	By Cauchy's inequality, we have
	\aryst
	(\cM^{n_1}(P_{n} H_{n}))^2 \leq \cM^{n_1} P_{n}^2 \cM^{n_1} H_{n}^2.
	\earyst
	It is clear that 
	\aryst
	\cM^{n_1} P^2_n \leq 1.
	\earyst
	By Lemma \ref{lem dim cancellation}, there is  $\kappa_4 > 0$ depending on $g$, $\kappa_5$ and $n_1$ such that 
	\aryst
	(\cM^{n_1} P^2_n)(z) < 1-\kappa_4, \quad \forall z \in \cP_1(y)
	\earyst
	where $y$ is any point given by Lemma \ref{lem dim cancellation}(1).
	
	It suffices to take $\Omega_{n+1}$ to be the union of all $\cP_1(y)$ given by Lemma \ref{lem dim cancellation}(1). By definition, $\Omega_{n+1}$ is $\cP_1$-measurable and $\cP_0$-dense in $\widetilde \Omega$.
\end{proof}

\subsection{Fast recurrence to uniform set}

We can verify Proposition \ref{prop L2 iteration}(3) using the following lemma.

\begin{lemma} \label{prop exp rec in vol exp}
	For any integer $n_1 > 0$ and a $n_1$-recurrent subset $\Omega \subset U$,
	there exist $\kappa_3, C_4 > 0$ depending on $g, F, n_1$ and $\Omega$
	such that the following is true.
    Let $\cP_0$ be a $\nu_U$-measurable finite partition such that $\cP_1 = \sigma^{-n_1}(\cP_0)  \succeq  \cP_0$, 
	let  $\widetilde \Omega$ be a $\cP_0$-measurable subset containing $\Omega$,  let 
	$\{ \Omega_n \}_{n \geq 1}$ be a sequence of $\cP_1$-measurable and $\cP_0$-dense subsets of  $\widetilde\Omega$.  Then for any $L \geq C_4$, we have
	\aryst
	\nu_U( \{ x \in U \mid |\{ 1 \leq j \leq L \mid \sigma^{jn_1}(x) \in \Omega_j \}| < \kappa_3 L \} )< e^{-L\kappa_3}.
	\earyst
\end{lemma}

We defer the proof of Lemma \ref{prop exp rec in vol exp} to Section \ref{sec proof of lemma prop exp rec in vol exp}.

\

\subsection{Concluding the proof of Proposition \ref{prop L2 iteration}}

Now we can provide the proof of Proposition \ref{prop L2 iteration}.
\begin{proof}[Proof of Proposition \ref{prop L2 iteration}]
	Under the hypothesis of  Proposition \ref{prop main 1} for a sufficiently large $C_1 > 0$, we let $C_8 > 1$, $\delta_1, \eta_1 > 0$ be given by Lemma \ref{lem dim cancellation}; let integer $n_1 > 0$ be sufficiently large so that Proposition \ref{prop main 1} are Lemma \ref{lem dim cancellation} are applicable; let $\{\Lambda^\epsilon\}_{\epsilon > 0}$ and the $n_1$-recurrent subset $\Omega$ with respect to $\nu_U$ be given by Proposition \ref{prop main 1} for $C_1, n_1$; and let $\kappa_5 > 0$ be sufficiently small.
	Given $a$ with $|a|$ sufficiently small and $b$ with $|b|$ sufficiently large, we define function $H_0 \equiv \norm{u}_{\theta, b}$. Then Proposition \ref{prop L2 iteration}(1) holds for $n=0$. 
	By Lemma \ref{lem dim cancellation}(1), we have \eqref{eq nn3uniform bound}-\eqref{eq lowerboundofhn} for $n=0$.
	
	Define $\cP_0, \cP_1$ by Lemma  \ref{lem dim cancellation}.  
	We can apply Lemma \ref{lem dim cancellation} to $u$ and $H_0$ to obtain $P_0$ and $H_1$ satisfying \eqref{eq n+1n3incone} to \eqref{eq n+1n3holder bound} for $n=0$.
	Let $\kappa_4$ be given by Corollary \ref{cor construct omega n+1} depending only on $g, F, \kappa_5$ and $n_1$. Let  $\widetilde \Omega$ be the minimal  $\cP_0$-measurable subset containing $\Omega$.  Then Corollary \ref{cor construct omega n+1} gives a $\cP_1$-measurable $\cP_0$-dense subset $\Omega_1$ of $\widetilde \Omega$ such that Proposition \ref{prop L2 iteration}(2) holds for $n=1$.
	
	Iterate the above argument,
	we can apply Lemma \ref{lem dim cancellation} and Corollary \ref{cor construct omega n+1} successively to define $P_n, H_n, \Omega_n$ for all $n \geq 1$ as long as we have \eqref{eq lowerboundofhn}.
	In particular, by letting $\kappa_5 < \eta_1/C'$, we have \eqref{eq lowerboundofhn} whenever $n \leq \ln |b|$  by inequalities 
	\aryst
	H_{m+1} = \cM^{n_1}(P_m H_{m}) \geq (1-\kappa_5) \cM^{n_1} H_m, \quad \forall 0 \leq m \leq n-1.
	\earyst
	Proposition \ref{prop L2 iteration}(1), (2) then follow from Lemma \ref{lem dim cancellation} and  Corollary \ref{cor construct omega n+1}.
	
	Let $ \kappa_3$ be given by Lemma \ref{prop exp rec in vol exp}.
	By our construction,  $\Omega_n$ is $\cP_1$-measurable and $\cP_0$-dense in $\widetilde \Omega$ for each $n \geq 1$. 	
		Note that for any $U$-random variable $X$ with distribution $\nu_U$, we know that $\sigma^{\lfloor C_8 \ln |b| \rfloor}X$ has distribution $\nu_U$ as above.
	Then we can apply Lemma \ref{prop exp rec in vol exp} to show that Proposition \ref{prop L2 iteration}(3) is satisfied for $\kappa_3$.  This concludes the proof. 
\end{proof}

\section{Partition} \label{sec partition}

\begin{defi}
	For any integer $k \geq 1$, for any $x \in U$ satisfying  
	\aryst
	\sigma^i(x) \in Int (U_{\a_i}), \quad \forall 0 \leq i \leq k
	\earyst	
	for a sequence $(\a_0 \cdots \a_{k}) \in \cA^{k+1}$,
	there exist a unique inverse branch  
	\aryst
	v_i : U_{\a_{i}} \to U_{\a_{i-1}}
	\earyst
	for each $1 \leq i \leq k$.
	In this case, we define {\em the  $k$-th cylinder at} $x$ by	
	\aryst
	U^x_{k} = v_{1} \cdots v_{k}(U_{\a_k}).
	\earyst
\end{defi}
It is clear that for any integer $k \geq 1$, the $k$-th cylinder is well-defined at a $\nu_U$-full measure subset.
 Clearly, for any integer $k \geq 1$, for any $x \in U$ such that
	the $k$-th cylinder at $x$ is well-defined, we know that the $k$-th cylinder  is also well-defined at any $z \in Int(U^x_k)$, and
	\aryst
	U^x_k = U^z_k.
	\earyst

\begin{lemma} \label{lem regular uniform set}
	There exists  $C_1 > 1$ depending only on $g$ such that the following is true.
	Given a stable sequence of functions $\{ \Lambda^\epsilon: U \to \R_{+} \}_{\epsilon > 0}$ , there exist $\kappa > 0$ and
	an integer $n' > 1$ such that for all sufficiently small $\epsilon > 0$, there is a finite partition $\cP$ of a  $\nu_U$-full measure subset of $U$ into cylinders
	\aryst 
	\cP = \{ U_i \mid 1 \leq i \leq l \},
	\earyst
	such that 
	\enmt
	\item each $U_i \in \cP$ contains some $z \in U_i$ satisfying
	\ary \label{eq tame lyapunov function gives good partition}
	U_i &\subset& W^u_{g}(z, C_1 \Lambda^\epsilon(z)^{-1});
	\eary
	\item for every $x \in U_i$, we have
	\ary
	\mbox{either } \ W^u_{g}(x, [0,\Lambda^\epsilon(x)^{-1})) \ \mbox{ or } \ W^u_{g}(x, (- \Lambda^\epsilon(x)^{-1}, 0]) \subset  U_i; \label{eq partition contain half neighborhood}
	\eary
	\item for any $n'' > n'$, for $\nu_U$-a.e. $x \in U$, every $v \in \sigma^{-n''}_x$ induces a $C^1$-contraction from $\cP(x)$ to $\cP(v(x))$ (in particular, $\sigma^{-n''}(\cP) \succeq \cP$). Moreover  we have
	\ary \label{eq stable lyapunov function gives increasing partition}
\norm{ D v|_{\cP(x)}} diam(\cP(x)) < C' e^{-n'' \kappa}diam(\cP(v(x))) .
	\eary
	\eenmt	
\end{lemma}

\begin{proof} 
	For any integer $k \geq 1$, for any $x \in U$ at which the $k$-th cylinder is well-defined, we define
	\aryst
	D(x, k) = \inf_{z \in U^x_k} \inf \{  \varrho \mid U^x_k \subset W^u_g(z, \varrho)  \}.
	\earyst
	We have seen that $D(x,k)$ depends on $x$ through $U^x_k$.
	By \eqref{contraction} and distortion estimate, we know that there is $C > 1$ depending on $g$ such that for $\nu_U$-a.e. $x \in U$, 
	\aryst
	D(x, k+1) \geq C^{-1} D(x,k), \quad \forall k \in \N. 
	\earyst
	By \eqref{contraction}, $diam(U^x_k)$ goes to $0$ as $k$ goes to infinity. Hence, $D(x,k)$ goes to $0$ as $k$ goes to infinity.
	
	Let $\epsilon$ be a small constant.	For $\nu_U$-a.e. $x \in U$, we set
	\aryst
	k_x = \inf \{ k \geq 1 \mid D(x,k) \leq C \sup_{y \in U^{x}_{ k} } \Lambda^\epsilon(y)^{-1}  \}.
	\earyst
	It is clear that $k_x$ is finite since $\Lambda^\epsilon \in L^\infty(U)$. Hence $diam(U^x_k)$ is uniformly lower bounded.
	Moreover, for $\nu_U$-a.e. $x $,  we have
	\ary \label{eq lower bound d x k}
	C\sup_{y \in U^x_{k_x}}\Lambda^\epsilon(y)^{-1} \geq D(x, k_x) \geq C^{-1} D(x, k_x-1) \geq  \sup_{y \in U^x_{k_x}}\Lambda^\epsilon(y)^{-1}.
	\eary

	For $\nu_U$-a.e. $x$, we have  for any $y \in U^x_{k_x}$,
	\ary \label{eq 1}
	U^x_{k_x} = U^y_{k_y};
	\eary
	and for any $y \notin U^x_{k_x}$ at which $U^y_{k_y}$ is defined, 
	\aryst
	Int (U^x_{k_x}) \cap Int (U^y_{k_y}) = \emptyset.
	\earyst
	Then
	\aryst
	\cP  = \{ U^x_{k_x} \mid x \in U \}
	\earyst
	is a finite partition of a $\nu_U$-full measure subset of $U$ into cylinders.
	By \eqref{eq lower bound d x k}, we obtain \eqref{eq partition contain half neighborhood}.
	
	For a $\nu_U$-typical $x \in U$, we select an arbitrary $z \in U^x_{k_x}$ such that
	\aryst
	\Lambda^\epsilon(z)^{-1} > \frac{1}{2}\sup_{y \in U^x_{k_x}}\Lambda^\epsilon(y)^{-1}. 
	\earyst 
	Then by \eqref{eq lower bound d x k}, we have
	\aryst
	2C\Lambda^\epsilon(z)^{-1} \geq D(x, k_x).
	\earyst
	By the definition of $D(x, k_x)$, there exists $y \in U^x_{k_x}$ such that
	\aryst
	U^x_{k_x} \subset W^u_g(y, 2C\Lambda^\epsilon(z)^{-1}).
	\earyst
	In particular, for all sufficiently small $\epsilon$, we have $z \in W^u_{g}(y, 1)$. Then by \eqref{eq local metric comparable} 
	\aryst
	W^u_g(y, 2C\Lambda^\epsilon(z)^{-1}) \subset W^u_g(z, 4C\Lambda^\epsilon(z)^{-1}).
	\earyst
	As a result, we obtain \eqref{eq tame lyapunov function gives good partition} with $C_1 = 4C$.

	We now show \eqref{eq stable lyapunov function gives increasing partition}. Let $n''$ be a large integer and denote  $\cP'' = \sigma^{-n''}(\cP)$. Let $x$ be a $\nu_U$-typical point in $U$. Then we have $\cP''(x) = v(\cP(\sigma^{n''}(x)))$ where $v \in \sigma^{-n''}$ satisfies that $x = v(\sigma^{n''}(x))$.  By construction, for some $y \in \cP''(x)$, we have $\cP''(x) = v(\cP(\sigma^{n''}(y))) =  g^{-\tau''}(W^u_{g}(g^{\tau''}(y), J))$ where $\tau'' = \tau_{n''}(y) \in (n''\tau_0, n''\tau_*)$; and some  interval $J$ containing $0$ satisfying $|J| \leq C'C_1 \Lambda^\epsilon(\sigma^{n''}(y))^{-1}$ (here we  use the $C^1$ bound of the center-stable holonomy maps between $W^u_{g, loc}(\sigma^{n''}(y))$ and $W^u_{g, loc}(g^{\tau''}(y))$). 
	By the hypothesis that $\Lambda$ is stable and by distortion estimate, there exists $\kappa > 0$ such that for all $n''$ sufficiently large depending on $g$ and $\{\Lambda^\epsilon\}_{\epsilon > 0}$, for all sufficiently small $\epsilon$,
	\aryst
	diam(\cP''(x))
	&\leq&C'\norm{Dv|_{\cP(\sigma^{n''}(x))}}diam(\cP(\sigma^{n''}(x) ))\\
	 &\leq& C'\norm{Dg^{\tau''}|_{E^u(y)}}^{-1}|J| \\
	&\leq& C'C_1\norm{Dg^{\tau''}|_{E^u(y)}}^{-1} \Lambda^\epsilon(\sigma^{n''}(y))^{-1} \\
	&\leq& C'C_1e^{-n'' \kappa} \Lambda^\epsilon(y)^{-1} \\
	&\leq& C'C_1 e^{-n'' \kappa}diam(\cP(y)) < diam(\cP(y)).
	\earyst
	By definition, $\cP(y) = \cP(x)$ and  both $\cP''(x)$ and $\cP(x)$ are cylinders at $x$. Consequently, $\cP''(x) \subset \cP(x)$. 
	
	%	Item \eqref{eq stable lyapunov function gives increasing partition}
	%	 follows from \eqref{eq tame lyapunov function gives good partition}, \eqref{eq stable function} and the uniform expansion of $\sigma$.
\end{proof}

\clb

\section{Bounds for smoothness}\label{sec: bounds for smoothness}

In this section, we will use the following notation. For each integer $n > 0$, for each collection of functions $\{ \varphi_v: Dom(v) \to \C \}_{v \in \sigma^{-n}}$, we denote for each $\alpha \in I$ and $x \in U_\a$ that
\aryst 
\sideset{}{^*}\sum_{v \in \sigma^{-n}} \varphi_v(x) =  \sum_{v \in \sigma^{-n}(\alpha)} \varphi_v(x).
\earyst

\subsection{Mapping between the cones} \label{subsec cone of functions}
For each $C^1$ function $\psi: U \to \R_+$, we define a  linear operator $ \cI_\psi: C^1(U, \R_+) \to C^1(U, \R_+)$ by
\aryst
\cI_\psi h(x) &=& h(x) \psi(x).
\earyst

\begin{lemma} \label{lem cone of functions}
	For any $g \in \cU^2$, any $F \in C^\theta(M)$, and any stable sequence of functions $\{ \Lambda^\epsilon : U \to \R_+ \}_{\epsilon > 0}$, there exist $\delta_1 \in (0,1)$ and an integer $n_4 > 0$ such that the following is true.  For all $a, b \in \R$ with $|a|$ sufficiently small, and  with $|b|$ sufficiently large, we define $f^{(a,b)}$ as in Section \ref{sec smoothing} for $\delta_1$, and
	set $\epsilon, \widetilde \cL$ and $\cM$ as in \eqref{def clwclcm}.
	Then for any sufficiently small $\epsilon$, for any $m > n_4$, any $\psi \in \cK_{\Lambda^\epsilon}$, we have
	\aryst
	\cM_{a,b}^{m} \cI_\psi( \cal K_{ \Lambda^\epsilon})  \subset \cal K_{ \Lambda^\epsilon}.
	\earyst
\end{lemma}
\begin{proof}  
		We fix $n, \kappa > 0$ such that $\{ \Lambda^\epsilon \}_{\epsilon > 0}$ is $( n, \kappa)$-stable.
	We take  
	\ary \label{eq delta3delta4}
	\delta_1 = \kappa/2.
	\eary

	Let $n'$ be given by Lemma \ref{lem regular uniform set}, and let $n_4 > n'$ be a large integer to be determined later. 	In the following, we let $\epsilon$ be a small constant, and abbreviate $\Lambda^\epsilon$, $J^{\Lambda^\epsilon}_x$, $\Phi^{\Lambda^\epsilon}_x$ as $\Lambda$, $J_x$, $\Phi_x$ respectively. 

	Take an arbitrary $h \in \cK_{ \Lambda}$, and let $m > n_4$ be an integer. 
	We set
	\aryst
	H = \cM_{a,b}^{m} \cI_\psi h.
	\earyst
	Let $x \in U$ and let $v \in \sigma^{-m}_x$.
	By Lemma \ref{lem regular uniform set} and by letting $n_4$ be sufficiently large depending only on $g$ and $\{ \Lambda^\epsilon \}_{\epsilon > 0}$, the map $B_v = \Phi_{v(x)}^{-1} v \Phi_x$ satisfies $B_v(J_x) \subset J_{v(x)}$, and we have
	\ary \label{eq contraction of B}
	\|DB_v\|_{C^0} < C' e^{-n_4\kappa}.
	\eary	
	 We have for any $z \in J_x$ that
	\aryst
	H(\Phi_x(z)) &=& \cM_{a,b}^{m} \cI_\psi h(\Phi_x(z)) \\
	&=& \sideset{}{^*}\sum_{v \in \sigma^{-m}} (e^{f^{(a,b)}_{m}} h \psi)(\Phi_{v(x)}(B_v(z))).
	\earyst
	Then for any $z \in J_x$ we have
	\aryst
	&& \partial_z[H(\Phi_x(z))] \\
	&=& \sideset{}{^*}\sum_{v \in \sigma^{-m}} [e^{f^{(a,b)}_{m}} h \psi](\Phi_{x}(z)) DB_v(z) \\
	&& \cdot (D[f^{(a,b)}_{m} \circ \Phi_{v(x)}]  + D[\log h \circ  \Phi_{v(x)}]  + D[\log \psi \circ  \Phi_{v(x)}] )(B_v(z)).
	\earyst
	Since $\{ \Lambda^\epsilon \}_{\epsilon > 0}$ is $(n, \kappa)$-stable, we have $\norm{D\Phi_x} \leq C' \Lambda^\epsilon(x)^{-1} \leq C'|b|^{-\kappa}$.
	Then by Lemma \ref{lem f a b}, \eqref{contraction} and \eqref{eq delta3delta4}, we have
	\aryst
    DB_v(z)D[ f^{(a,b)}_{m} \circ \Phi_{v(x)}](B_v(z)) =  D[ f^{(a,b)}_{m} v \Phi_x ](z)  \leq C'|b|^{\delta_1}\norm{D\Phi_{x}} < C' |b|^{\delta_1 - \kappa} < 1/2.
	\earyst
	Since $h, \psi \in \cK_{\Lambda}$, we have
	\aryst
	 |D[ \log h \circ  \Phi_{v(x)}]| \leq  1, \quad  |D[ \log \psi \circ  \Phi_{v(x)}]| \leq  1.
	\earyst
	Summarizing the above estimates, we obtain
	\aryst
\partial_z[H(\Phi_x(z))]  \leq \frac{1}{2} H(\Phi_x(z)) + \frac{1}{2} H(\Phi_x(z)) =  H(\Phi_x(z)), \quad \forall z \in J_x.
	\earyst
	This concludes the proof.
\clb
\end{proof}

\subsection{Decay of H\"older semi-norms}

%Using the above definition, we may summarize Corollary \ref{lem template approximation of return time function} as 
%\ary  \label{eq epsiloninversetaumv}
%\epsilon^{-1}\tau_m \circ v \in \cH^{\Lambda^\epsilon}_{\epsilon^\delta/C', C'{\bf 1}, 1}
%\eary 
%for some $C' > 1$ depending only on $g$, and for all $m \geq 1$ and all $v \in \sigma^{-m}$.

\begin{lemma} \label{lem decayofholderseminorms}
	For a stable and tame sequence of functions $\{ \Lambda^\epsilon: U \to \R_+\}_{\epsilon > 0}$, there is $\eta_1 > 0$ such that any sufficiently large $C_8 > 1$ and integer $n_1 > 0$ the following is true.
	For all $a \in \R$ with $|a|$ sufficiently small, and all $b \in \R$ with $|b|$ sufficiently large, we define $f^{(a,b)}$ as in Section \ref{sec smoothing} for $\delta_1$ given by Lemma \ref{lem cone of functions}, and set $\epsilon$, $\widetilde \cL$ and $\cM$ as in \eqref{def clwclcm}.
	Then for any sufficiently small $\epsilon > 0$ we have the following:
	\enmt
	\item for any integer $m \geq \lfloor C_8 \ln |b| \rfloor$, for any $u \in C^\theta(U)$ we have
	\aryst
	\widetilde \cL^m u \in \cH^{\Lambda^\epsilon}_{\epsilon^{\eta_1} \norm{u}_{C^0},  C_8\norm{u}_{\theta, b}};
	\earyst
	\item for any functions $W, H \in \cK_{ \Lambda^\epsilon}$ and  $u \in \cH^{\Lambda^\epsilon}_{W,  C_8H}$ satisfying
	\ary \label{eq ulessthanh}
	|u(x)| \leq H(x), \quad \forall x \in U,
	\eary
	we have
	\aryst
	\widetilde\cL^{n_1}u \in \cH^{\Lambda^\epsilon}_{ \cM^{n_1}(W +  \epsilon^{\eta_1} H) ,\frac{C_8}{4}\cM^{n_1} (H+W)}.
	\earyst
	\eenmt 
\end{lemma}
\begin{proof}
	We fix $n, \kappa, C,  \eta_1 > 0$ such that $\{ \Lambda^\epsilon \}_{\epsilon > 0}$ is $( n, \kappa)$-stable and $(C, \eta_1)$-tame.
	
	Let $m \geq 1$ be an aribtrary integer. We set
	\aryst
	\hat u = \widetilde \cL^{m} u, \quad  \rho =e^{f^{(a,b)}_{m}}.
	\earyst
	By definition, we have
	\aryst
	\hat u = \sideset{}{^*}\sum_{v \in \sigma^{-m}}  ( u \rho e^{ib \tau_{m}}) \circ v.
	\earyst

We let $\epsilon > 0$ be a small constant, and, as in Lemma \ref{lem cone of functions}, we abbreviate $\Lambda^\epsilon$, $J^{\Lambda^\epsilon}_x$, $\Phi^{\Lambda^\epsilon}_x$ as $\Lambda$, $J_x$, $\Phi_x$ respectively. 

	Fix $x \in U$. We set
	\aryst
	\bar u = \hat u \circ \Phi_x, \ \ u_v =  u \circ v  \Phi_x, \ \
	\rho_v = \rho \circ v \Phi_x, \ \ \theta_v = e^{ib\tau_{m}}\circ v  \Phi_x.
	\earyst
	Then 
\aryst
\bar u = \sideset{}{^*}\sum_{v \in \sigma^{-m}}  u_v  \rho_v  \theta_v .
\earyst
	By our notations, we have 
	\aryst
	\sideset{}{^*}\sum_{v \in \sigma^{-m}} \rho_v = 1. 
	\earyst

	Since $\{ \Lambda^\epsilon\}_{\epsilon > 0}$ is $(C, \eta_1)$-tame,  for every $x \in U$ and every $v \in \sigma^{-n_1}_x$, 
	there exists $\Theta_v \in C^\theta(J_x)$ such that
	\ary \label{eq thetavThetav}
	\norm{\Theta_v}_{C^0} = 1	, \quad	
	|\Theta_v|_\theta < 2\pi C, \quad \norm{ \theta_v - \Theta_v}_{C^0} < \epsilon^{\eta_1}.
	\eary
	Since $\{\Lambda^\epsilon\}_{\epsilon > 0}$ is $(n, \kappa)$-stable,
	by letting $m$ be sufficiently large, the map 
	\aryst
	B_v = \Phi_{v(x)}^{-1} v \Phi_x
	\earyst satisfies $B_v(J_x) \subset J_{v(x)}$. 
	 	Hence
	\ary \label{eq rhovc0bound}
	\|\rho_v  \|_{J_x} \leq 
	\|\rho\circ \Phi_{v(x)}  \|_{J_{v(x)}}. 
	\eary
		By Lemma \ref{lem f a b}, we have 
	\ary \label{eq holderboundforrho}
	|\rho \circ v|_\theta \leq C' \norm{\rho \circ v}_{C^0}.
	\eary
	Then we have
	\ary \label{eq rhovthetarhovc0}
	| \rho_v	|_{\theta} \leq \norm{D\Phi_x}^\theta |\rho \circ v|_\theta \leq C' \Lambda(x)^{-\theta} \|\rho_v  \|_{C^0}.% \leq  C' \Lambda(x)^{-\theta} \|\rho\circ  \Phi^u_{v(x)}(\Lambda(v(x))^{-1}\cdot)  \|_{C^0}.
	\eary
	
	By Lemma \ref{lem regular uniform set} and by enlarging $m$ if necessary we obtain
	\ary \label{eq contraction of B 2}
	\|DB_v\|_{C^0} < C' e^{-m \kappa}.
	\eary

	By  direct computation, we can choose $C_6 > 1$ depending only on $g$ such that
	\ary \label{eq rhovproportional}
	C^{-1/2}_6 < \frac{\rho_v (s)}{ \rho_v (0)}  < C^{1/2}_6, \quad \forall s	\in J_{v(x)}. 
	\eary
	
	\noindent{Proof of item (1)}:
	Set
\ary 
u' = \sideset{}{^*}\sum_{v \in \sigma^{- m}} u_v \rho_v \Theta_v.
\eary
Then by \eqref{eq thetavThetav} we have
\aryst
|\bar u - u'| < \epsilon^{\eta_1} \norm{u}_{C^0}.
\earyst

By  \eqref{eq thetavThetav}, \eqref{eq holderboundforrho}, \eqref{eq contraction of B 2} and \eqref{eq rhovproportional} 
we have
\aryst
|u'|_\theta &\leq&  \sideset{}{^*}\sum_{v \in \sigma^{-m}}|u_v|_\theta \norm{\rho_v}_{C^0} + \norm{u_v}_{C^0} |\rho_v|_\theta + \norm{u_v}_{C^0}\norm{\rho_v}_{C^0}|\Theta_v|_\theta \\
&\leq& C'e^{-m\kappa \theta} |u|_\theta + C'C\norm{u}_{C^0}.
\earyst
By letting $m > C_8 \ln |b|$ and letting $C_8$ be sufficiently large depending only on $g$, $\kappa$ and $C$, we have
\aryst
|u'|_\theta < C_8 \norm{u}_{\theta, b}.
\earyst	
This concludes the proof of item (1).
	
	\
	
	\noindent{Proof of item (2)}:
		Fix an arbitrary $v \in \sigma^{- m}$.
	By hypothesis, we have that
	\ary \label{eq uvhcircv} 
	|u_v(s)| \leq H( v(\Phi_x(s))), \quad \forall s \in J_x.
	\eary
	By \eqref{eq ratio h} and by enlarging $C_6$ if necessary we have
	\ary \label{eq ratio h v 1 v 2}
	C^{-1/2}_6 < \frac{H( \Phi_{v(x)}(s))}{ H(v(x))}  < C^{1/2}_6, \quad \forall s	\in J_{v(x)}. 
	\eary
	
	Since $u \in \cH^{\Lambda^\epsilon}_{W, C_8 H}$,
	there is $\tilde u_v \in C^\theta(J_{v(x)})$ such that
	\ary \label{eq def cH epsilon H1}
	|\tilde u_v|_\theta &\leq& C_8 H(v(x	)), \\
	| u (\Phi_{v(x)}(s))  - \tilde u_v(s) | &\leq& W(\Phi_{v(x)}(s)),  \quad \forall s \in J_{v(x)}.
	\eary
	Thus for every $s \in J_x$, we have
	\ary
	|u_v(s) - \tilde u_v \circ B_v(s)| 
	&=& |(u \circ \Phi_{v(x)}  - \tilde u_v) \circ B_v(s)| \nonumber \\
	&\leq&  W (\Phi_{v(x)}(B_v(s))) = W ( v (\Phi_x(s))). \label{eq uvtildeuvcircb}
	\eary
	By \eqref{eq contraction of B 2} and \eqref{eq def cH epsilon H1}, we have
	\ary \label{eq tilde u v circ B theta}
	|\tilde u_v \circ B_v|_{\theta} \leq  C' e^{-m\kappa\theta} |\tilde u_v|_{\theta}  
	\leq  C_8C' e^{-m\kappa\theta} H(v(x)).
	\eary
	Set
	\ary \label{def defofu'}
	u' = \sideset{}{^*}\sum_{v \in \sigma^{-m}} \tilde u_v \circ B \rho_v \Theta_v.
	\eary
	Then by  \eqref{eq uvtildeuvcircb}, \eqref{eq uvhcircv}, \eqref{eq ulessthanh} and \eqref{eq thetavThetav}, for every $s \in J_x$ we have
	\ary
	\label{eq conditionh1}
	&& |\bar u(s)  - u'(s)| \leq \sideset{}{^*}\sum_{v \in \sigma^{-m}} | ( u_v  \rho_v   \theta_v)(s) - (\tilde u_v \circ B \rho_v \Theta_v)(s) |  \\
	&\leq& \sideset{}{^*}\sum_{v \in \sigma^{-m}}(|u_v(s) - \tilde u_v \circ B(s)| |\Theta_v(s)| \rho_v(s) +  |\theta_v(s) - \Theta_v(s)| |u_v(s)| \rho_v(s) )  \nonumber \\
	&\leq& \sideset{}{^*}\sum_{v \in \sigma^{-m}} \rho_v(s)(W( v( \Phi_x(s))) + \epsilon^{\eta_1} |u_v(s)|)   \nonumber  \\
	&\leq& [\cM^{m}(W +  \epsilon^{\eta_1} H)]( \Phi_x(s) ).  \nonumber
	\eary
	
	By \eqref{eq thetavThetav}, \eqref{eq tilde u v circ B theta}, \eqref{eq rhovthetarhovc0},
	we  have
	\aryst
	|\tilde u_v \circ B \rho_v \Theta_v|_\theta 
	&\leq&  |\tilde u_v \circ B|_{\theta}\norm{\rho_v \Theta_v}_{C^0} +   |\rho_v|_\theta \norm{\tilde u_v \circ B \Theta_v}_{C^0} + |\Theta_v|_\theta \norm{\tilde u_v \circ B \rho_v}_{C^0}  \\
	&\leq& C'C_8C_6e^{- m\kappa\theta} H(v(x))\rho_v(0) + 
	C'C_6 \Lambda(x)^{-\theta}  (H + W)(v(x)) \rho_v(0) \\
	&+& C C_6(H + W)(v(x))\rho_v(0) \\
	&\leq&C_6(C'C_8e^{- m\kappa\theta} + C' \Lambda(x)^{-\theta}  + C ) (H+W)(v(x))\rho_v(0)	.
	\earyst
	By taking $|b|$ sufficiently large depending on $g$, we can make $\Lambda(x)^{-\theta}$ arbitrarily small. Then by making $|b|$ sufficiently large depending on $g, C$, and by making   $m,  C_8$ sufficiently large depending on $g, C, C_6$, we have
	\ary \label{eq conditionh2}
	|\tilde u_v\circ B \rho_v  \Theta_v |_{\theta}
	&\leq& \frac{C_8}{4}((H+W) \rho) (v(x)).
	\eary
	By \eqref{def defofu'} and \eqref{eq conditionh2}, we obtain
	\aryst
	|u'|_\theta < \frac{C_8}{4} \cM^{m} (H+W)(x).
	\earyst
	By \eqref{eq conditionh1} and \eqref{eq conditionh2}, we obtain $\hat u \in \cH^\Lambda_{ \cM^{m}(W +  \epsilon^{\eta_1} H) ,\frac{C_8}{4}\cM^{m} (H+W)}$. This concludes the proof of item (2).
\end{proof}

\section{Proof of Lemma \ref{lem dim cancellation}}\label{sec proof lem dim cancellation}

Notice that Lemma \ref{lem dim cancellation}(1) follows from Lemma \ref{lem decayofholderseminorms}(1). It remains to show Lemma \ref{lem dim cancellation}(2).
We will verify Lemma \ref{lem dim cancellation}(2) by an argument similar to the one in \cite[Lemma 10'']{Dol} as follows. 

Let $n_1$ be a large integer to be determined in due course.
By hypothesis, there exist constants $C_5 > 0$, $\kappa_6 \in (0, \frac{1}{10})$
such that
\enmt
\item 
$\{ \Lambda^\epsilon \}_{\epsilon > 0}$ is $(n_1, C_5)$-adapted to $\Omega$;
\item 
 $(C_1, \kappa_6)$-{\rm UNI} holds on $\Omega$ at scales $\{ \Lambda^\epsilon \}_{\epsilon > 0}$.
\eenmt
We fix a constant $\delta_1 > 0$ given by Lemma \ref{lem cone of functions}, and fix a constant $\eta_1 > 0$ given by Lemma \ref{lem decayofholderseminorms}.

Let $a,b, u, n, H_n$ be given by Lemma \ref{lem dim cancellation}.
Recall that
$f^{(a,b)}$ is defined for $\delta_1$;
and $\epsilon$, $\widetilde \cL$ and $\cM$ are given by \eqref{def clwclcm}. 
We let $C_1$ be sufficiently large so that Lemma \ref{lem regular uniform set} is applicable, and define $\cP_0$ by Lemma \ref{lem regular uniform set} applied to $\Lambda^\epsilon$.
We will abbreviate $\Phi^{\Lambda^\epsilon}_x, J^{\Lambda^\epsilon}_x$ as $\Phi^\epsilon_x, J^\epsilon_x$ respectively.

Let $C_8 > 1$ be a large integer to be determined depending on $g$ and $\{ \Lambda^\epsilon \}_{\epsilon > 0}$. Recall that $u_0 =  \cL^{\lfloor C_8 \ln |b| \rfloor}u$.
Denote 
\aryst
\tilde u = \widetilde \cL^{nn_1}u_0.
\earyst
  By  the assumptions \eqref{eq nn3uniform bound} and \eqref{eq lowerboundofhn}, we have for all $x \in U$ that
\ary \label{eq tildeulessthanH}
|\tilde u(x)| \leq H_n(x) \leq \norm{u}_{\theta, b}, \quad H_n \geq (n+1)\epsilon^{\eta_1/2}\norm{u}_{\theta, b}.
\eary
Assume that for some $\kappa_5 \in (0,1/2)$ we have found $P_n \in C^1(U, [1-\kappa_5,1]) \cap \cK_{\Lambda^\epsilon}$ satisfying \eqref{eq valueofpnonsomesmallinterval} and \eqref{eq n+1n3uniform bound}. 
By letting  $n_1$ be sufficiently large,  we may apply Lemma \ref{lem cone of functions} to see that for all sufficiently large $n_1$, for all sufficiently small $\epsilon$,  
\aryst
H_{n+1} = \cM^{n_1}(P_n H_n) \in \cK_{\Lambda^\epsilon}.
\earyst
Moreover, by $\kappa_5 < 1/2$ we have   
\ary \label{eq hn+1geq12cmhn}
H_{n+1}  \geq \frac{1}{2}\cM^{n_1} H_n.
\eary 
Let $C_8$ and $n_1$ be sufficiently large depending on $g$ and $\{ \Lambda^\epsilon \}_{\epsilon > 0}$ so that Lemma \ref{lem decayofholderseminorms} is applicable.
Then by applying Lemma \ref{lem decayofholderseminorms}
to $W = (n+1) \epsilon^{\eta_1} \norm{u}_\theta$ and $H = H_n$, we 
obtain that
\aryst
\widetilde\cL^{(n+1)n_1}u_0 &\in& \cH^{\Lambda^\epsilon}_{\cM^{n_1}((n+1) \epsilon^{\eta_1} \norm{u}_\theta +  \epsilon^{\eta_1} H_n ),  \frac{C_8}{4}\cM^{n_1} (H_{n} + (n+1) \epsilon^{\eta_1} \norm{u}_{\theta, b})  } \\
&\subset&
\cH^{\Lambda^\epsilon}_{\cM^{n_1}((n+2) \epsilon^{\eta_1} \norm{u}_{\theta, b}  ),  \frac{C_8}{4}\cM^{n_1} (H_{n} + (n+1) \epsilon^{\eta_1} \norm{u}_{\theta, b})  } \\
& \subset& \cH^{\Lambda^\epsilon}_{(n+2) \epsilon^{\eta_1} \norm{u}_{\theta, b}, C_8H_{n+1}}.
\earyst
The last two inclusions above follow from   \eqref{eq tildeulessthanH} and \eqref{eq hn+1geq12cmhn}. This gives \eqref{eq n+1n3holder bound}.

We now define $P_n$, by specifying the values of $P_n$ on each atom of $\sigma^{-n_1}\cP_0$. 

Let $\cP_0(x_0)$ be an atom of $\cP_0$ where $x_0$  is in the $\nu_U$-full measure set in Lemma \ref{lem regular uniform set}. If $\cP_0(x_0)$ is disjoint from $\Omega$, then we set
\aryst
P_n|_{\sigma^{-n_1}(\cP_0(x_0))} \equiv 1.
\earyst

Now we assume in the rest of the proof that 
$\cP_0(x_0)$ meets $\Omega$. In this case, we can assume without of generality that
$x_0 \in \Omega$. By Lemma \ref{lem regular uniform set}, there exist $x \in \cP_0(x_0)$ and  an interval $J \subset (-C_1,C_1)$ containing either $[0,1)$ or $(-1,0]$, such that  \clb
\ary \label{eq choose x}
\cP_0(x_0) = \cP_0(x) = W^u_{g}(x, \Lambda^\epsilon(x)^{-1} J).
\eary

Without loss of generality, let us assume that $[0,1) \subset J$, and consequently, $[0,1) \subset J^\epsilon_x$. 
Then $\Phi^\epsilon_x$ induces a diffeomorphism from $[0,1)$ to a subset of $\cP_0(x)$.
By Lemma \ref{lem regular uniform set}(2), \eqref{eq local metric comparable} and \eqref{eq choose x}, we have
\ary \label{eq lambdaxlambdax0}
\Lambda^\epsilon(x) < 2C_1\Lambda^\epsilon(x_0).
\eary

In the following, we denote by $Arg(w)$ the argument of of a complex number $w \in \C \setminus \{0\}$,
i.e., the unique number in $\R / 2\pi \Z$ with $w = |w|e^{i Arg(w)}$. 

\begin{lemma} \label{lem small u v big u v}
	There exists  $C_9 > 1$ depending only on $g$ such that for all $n_1$ sufficiently large depending only on $g$, $\kappa_6$ and $\{ \Lambda^\epsilon \}_{\epsilon > 0}$,  for any $v \in \sigma^{-n_1}_x$, we have 
	\enmt
    \item[ either   $\it{(1)}$]  for any $y \in \Phi^\epsilon_x([0,1))$
	\aryst
	|(e^{f^{(a,b)}_{n_1}} \tilde u) ( v(y) )| \leq \frac{3}{4} (e^{f^{(a,b)}_{n_1}} H_n)(v(y)),	
	\earyst
	\item[ or $\it{(2)}$]  there exists $\omega \in \R/2\pi \Z$ such that for any $y \in \Phi^\epsilon_x([0,1))$
	\aryst
	&&|(e^{f^{(a,b)}_{n_1}} \tilde u)( v(y) )| > C_9^{-1} (e^{f^{(a,b)}_{n_1}} H_n)(v(y))  \\
	&\mbox{and  } & \| Arg((e^{f^{(a,b)}_{n_1}} \tilde u)( v(y)) ) -  \omega \|_{\R / 2\pi \Z} < \frac{1}{100}\kappa_6.
	\earyst
	\eenmt
\end{lemma}
\begin{proof}
	Already in the proof of Lemma \ref{lem decayofholderseminorms}, we have seen that there is $C_6 > 1$ depending only on $g$ such that for any $x \in U$ and any $s \in J^\epsilon_{x}$
	\ary \label{eq ratio h v 1 v 2X}
	C^{-1/2}_6 < \frac{ H_n( \Phi^\epsilon_{x}(s))}{H_n(x)},  \frac{e^{f^{(a,b)}_{n_1} \circ v(\Phi^\epsilon_x(s))}}{e^{f^{(a,b)}_{n_1}(v(x))}} < C^{1/2}_6.
	\eary

	Assume that item (1) does not hold. Then there exist  $v \in \sigma^{-n_1}_x$ and  $s \in [0,1)$ such that $y= \Phi^\epsilon_x(s)$ satisfies that
	\aryst
	|(e^{f^{(a,b)}_{n_1}} \tilde u)(v(y))| \geq \frac{3}{4}(e^{f^{(a,b)}_{n_1}} H_n)(v(y)).
	\earyst
  By letting $n_1$ be sufficiently large, we have 
	\aryst
	v(\cP_0(x)) \subset v(W^u_g(x, C_1\Lambda^\epsilon(x)^{-1})) \cap U \subset W^u_{g}(v(x), \Lambda^\epsilon(x)^{-1}) \cap U.
	\earyst
	Then by \eqref{eq ratio h v 1 v 2X} we have
	\ary \label{eq efabn3tildeulowerbound}
	|(e^{f^{(a,b)}_{n_1}} \tilde u)(v(y))| \geq \frac{3}{4C_6}\norm{( e^{f^{(a,b)}_{n_1}} H_n ) \circ v}_{\Phi^\epsilon_x([0,1))}.
	\eary
	By Definition \ref{def hlambdaehc}, \eqref{eq nn3holder bound}, \eqref{eq tildeulessthanH} and Lemma \ref{lem regular uniform set}, 
	there exists $\kappa > 0$ depending only on $g$ and $\{ \Lambda^\epsilon \}_{\epsilon > 0}$, and there is a function  $\hat u \in C^\theta(J)$ such that
	\ary \label{tildeucircvapp1}
	\norm{\tilde u \circ v(\Phi^\epsilon_x(\cdot))) - \hat u}_{J}  &<& (n+1)\epsilon^{\eta_1} \norm{u}_{\theta, b}  < \epsilon^{\eta_1/2} H_n(v(x)),\\
	\label{tildeucircvapp2}
	|\hat u|_\theta &<&   C_8C' e^{-n_1\kappa\theta} H_n(v(x)). 
	\eary
	Here we deduce \eqref{tildeucircvapp2} by 
	\aryst
	\norm{D[(\Phi^\epsilon_{v(x)})^{-1} v \Phi^\epsilon_x]} \leq C'e^{-n_1\kappa},
	\earyst
	which follows from Lemma \ref{lem regular uniform set}.
	Then by \eqref{eq nn3uniform bound}, \eqref{eq ratio h v 1 v 2X} and \eqref{tildeucircvapp1}, we have
	\ary \label{eqnormujupperbound}
	\norm{\hat u}_{J} \leq 2C_6 H_n(v(x)).
	\eary
	By \eqref{tildeucircvapp1} and \eqref{eq ratio h v 1 v 2X} we have 
	\ary
	\label{eq uftildeucompare} 
	&&\norm{ (\tilde u e^{f^{(a,b)}_{n_1}})\circ v \circ \Phi^\epsilon_x  - \hat u  e^{f^{(a,b)}_{n_1}}\circ v \circ \Phi^\epsilon_x}_{[0,1)}  \\
	&<& C_6\epsilon^{\eta_1/2} (H_ne^{f^{(a,b)}_{n_1}})(v(x)). \nonumber 
	\eary
	By \eqref{eqnormujupperbound}, \eqref{tildeucircvapp2}, \eqref{eq rhovthetarhovc0} and the hypothesis that $\{\Lambda^\epsilon\}_{\epsilon > 0}$ is stable, there is a constant $\eta > 0$ depending on $\{\Lambda^\epsilon\}_{\epsilon > 0}$ such that for all sufficiently small $\epsilon$ we have
	\aryst  
	 && |\hat u e^{f^{(a,b)}_{n_1}}(v (\Phi^\epsilon_x(\cdot)))|_{\theta, [0,1)}  \\
	&<&
	\norm{\hat u}_{[0,1)} |e^{f^{(a,b)}_{n_1}}(v (\Phi^\epsilon_x(\cdot)))|_{\theta, [0,1)}
	+
	|\hat u|_\theta \norm{e^{f^{(a,b)}_{n_1}}(v (\Phi^\epsilon_x(\cdot)))}_{[0,1)} \nonumber \\
	&<& C'(  \Lambda^\epsilon(x)^{-1} )^\theta \norm{e^{f^{(a,b)}_{n_1}} \circ v}_{\Phi^\epsilon_x([0,1))} \norm{\hat u}_{J}  + 
	C_8 C' e^{-n_1\kappa\theta} H_n(v(x)) \norm{e^{f^{(a,b)}_{n_1}} \circ v}_{\Phi^\epsilon_x([0,1))}  \nonumber  \\
	&\leq& 10C^2_6C'( \epsilon^{\theta\eta}  +C_8 e^{-n_1\kappa\theta}) (H_ne^{f^{(a,b)}_{n_1}} )( v(x)). \nonumber
	\earyst
	By letting $n_1$ be sufficiently large depending on $g,  \kappa_6, C_6, C_8, \kappa, \theta$, and by letting $\epsilon$ be sufficiently small (or equivalently, letting $|b|$ be sufficiently large) depending on $g,  \kappa_6, C_6, \eta, \theta$, we have
	\ary \label{eq lower bound for u v phi u x ---}
	\ 	\ \  \  \  \  \   \  \  \  \   \  \  \  	|\hat u e^{f^{(a,b)}_{n_1}} (v(\Phi^\epsilon_x(\cdot))) |_{\theta, [0,1)} < \frac{\kappa_6}{20^4C_6^{2} }  \norm{( e^{f^{(a,b)}_{n_1}} H_n ) \circ v}_{\Phi^\epsilon_x([0,1))}.
	\eary
	Then by \eqref{eq uftildeucompare} and by letting $\epsilon$ be sufficiently small, we obtain
	\ary \label{eq oscoftildeuef}
	Osc_{[0,1)} (\tilde u e^{f^{(a,b)}_{n_1}})( v(\Phi^\epsilon_x(\cdot))) < \frac{\kappa_6}{20^4C_6^{2} }  \norm{( e^{f^{(a,b)}_{n_1}} H_n ) \circ v}_{\Phi^\epsilon_x([0,1))}. 
	\eary
	
	Given an arbitrary $z \in \Phi^\epsilon_x([0,1))$.  We have 
	\ary
	| (\tilde u e^{f^{(a,b)}_{n_1}})(v(z))| 
	&\geq& 	| (\tilde u e^{f^{(a,b)}_{n_1}})(v(y))| - 2\norm{ (\tilde u e^{f^{(a,b)}_{n_1}})\circ v \circ \Phi^\epsilon_x  - \hat u  e^{f^{(a,b)}_{n_1}}\circ v \circ \Phi^\epsilon_x}_{[0,1)} \nonumber \\
	&& - |\hat u  e^{f^{(a,b)}_{n_1}}\circ v \circ \Phi^\epsilon_x|_{\theta, [0,1)}    \nonumber  \\  
(\mbox{by }\eqref{eq efabn3tildeulowerbound}, \eqref{eq uftildeucompare}, \eqref{eq lower bound for u v phi u x ---})	& \geq& \frac{3}{4C_6}(e^{f^{(a,b)}_{n_1}} H_n)(v(z)) - 2C_6\epsilon^{\eta_1/2} (H_ne^{f^{(a,b)}_{n_1}})(v(x))  \nonumber \\
	&-& \frac{\kappa_6}{20^4C_6^{2} }  \norm{( e^{f^{(a,b)}_{n_1}} H_n ) \circ v}_{\Phi^\epsilon_x([0,1))} \nonumber \\
(\mbox{by } \eqref{eq ratio h v 1 v 2X})	& \geq& \frac{3}{4C_6}(e^{f^{(a,b)}_{n_1}} H_n)(v(z)) - 2C_6^2\epsilon^{\eta_1/2} (H_ne^{f^{(a,b)}_{n_1}})(v(z))  \nonumber \\
	&-& \frac{1}{20}C_6^{-1} (e^{f^{(a,b)}_{n_1}} H_n)(v(z)) \nonumber \\
	&\geq& (100C_6)^{-1} (H_ne^{f^{(a,b)}_{n_1}})(v(z)). \label{eq tildeulowerbound}
	\eary
	Assume to the contrary that item (2) fails. Then there would exist $z_1, z_2 \in \Phi^\epsilon_x([0,1))$ such that
	\ary \label{eq differenceofargisnotsmall}
\quad \|Arg((e^{f^{(a,b)}_{n_1}} \tilde u)( v(z_1))) - Arg( (e^{f^{(a,b)}_{n_1}} \tilde u)(v(z_2))) \|_{\R / 2\pi \Z} > \frac{1}{100}\kappa_6.
	\eary
	Then by \eqref{eq uftildeucompare}, \eqref{eq tildeulowerbound} and \eqref{eq differenceofargisnotsmall}  for all sufficiently small $\epsilon > 0$ we have
	\aryst
	|(e^{f^{(a,b)}_{n_1}} \tilde u)(v(z_1)) - (e^{f^{(a,b)}_{n_1}} \tilde u)( v(z_2))| >\frac{\kappa_6}{20^4C_6^2}\norm{(e^{f^{(a,b)}_{n_1}} H_n)\circ v}_{\Phi^\epsilon_x([0,1))}.
	\earyst
	This contradicts \eqref{eq oscoftildeuef}.
\end{proof}

We now define $P_n$ on $\sigma^{-n_1}(\cP_0(x_0))$  in the case where $\cP_0(x_0)$ meets $\Omega$. We have the following two cases to consider.

\smallskip

\noindent{(1)}
If there exists an inverse branch $v \in \sigma^{-n_1}_x$ such that
Lemma \ref{lem small u v big u v}(1) holds, then we may construct $P_n$ as follows. 
Let $\zeta \in C^\infty([0,1), [1/2,1])$ such that $\zeta(s) = 1-\kappa_5$ for all $s \in (1/4,3/4)$, and $\zeta$ equals $1$ on $[0,1/8] \cup [7/8, 1)$. Moreover, we may also require that $\norm{\zeta}_{C^1} \leq 10\kappa_5$. Then for any $z \in \sigma^{-n_1}(\cP_0(x_0))$ we set
\ary \label{eq defnofpn1}
P_n(z) = 
\begin{cases}
	1, & z \notin v(\Phi^\epsilon_x([0,1))), \\
	\zeta((\Phi^\epsilon_x)^{-1}(\sigma^{n_1}(z))), & z \in  v(\Phi^\epsilon_x([0,1))).
\end{cases}
\eary
By Lemma \ref{lem small u v big u v}(1), \eqref{eq tildeulessthanH} and by requiring $\kappa_5 < \frac{1}{4}$, we have
\aryst
|(e^{f^{(a,b)}_{n_1}} \tilde u)( \tilde v(z))| \leq (e^{f_{n_1}^{(a,b)}}P_{n} H_n )(\tilde v(z)), \quad \forall z \in U,  \tilde v \in \sigma^{-n_1}_x.
\earyst
By requiring $n_1$  be sufficiently large depending only on $g$, there exists $y \in \cP_0(x)$ such that $P_n|_{v(\cP_1(y))} = 1-\kappa_5$. In this case, we clearly have \eqref{eq n+1n3uniform bound} for all $x \in \cP_0(x_0)$.

\smallskip

\noindent{(2)}
We now assume that Lemma \ref{lem small u v big u v}(1) fails for every $v \in \sigma^{-n_1}_x$. Then for every $v \in \sigma^{-n_1}_x$, there exists $\omega_v \in \R / 2\pi \Z$ such that Lemma \ref{lem small u v big u v}(2) holds with $\omega_v$ in place of $\omega$.

By Corollary \ref{cor eq b tau k}, for any inverse branch $v \in \sigma^{-n_1}_x$, 
we have
\ary \label{eq taunvtaunv}
\tau_{n_1}(v(z)) - \tau_{n_1}(v(x)) =  \Psi_{x, x^{v}}(z).
\eary

By the hypothesis that $\{\Lambda^\epsilon \}_{\epsilon > 0}$ is $(C, \eta)$-tame for some $C, \eta > 0$, we have 
\aryst 
\|\Psi_{x, \Phi^s_x(y)}(\Phi^\epsilon_x(\cdot))\|_{[0,1)} < C \epsilon |y|^\eta + \epsilon^{1+\eta}.
\earyst 
By \eqref{eq taunvtaunv}, we obtain for every $v_2 \in \sigma^{-n_1}_x$ that
\aryst
Osc_J [b \tau_{n_1}(v_2(\Phi^\epsilon_x(\cdot)))] = O(d(x, x^{v_2}) + \epsilon^\eta).
\earyst 
In particular, by letting $n_1$ be sufficiently large, and letting $\epsilon$ be sufficiently small, there exists some $v_2 \in \sigma^{-n_1}_x$ so that
\aryst
Osc_J [b \tau_{n_1}(v_2(\Phi^\epsilon_x(\cdot))) ] < \frac{1}{100} \kappa_6.
\earyst

By the hypothesis that $(C_1, \kappa_6)$-{\rm UNI} holds on $\Omega$ at scales $\{ \Lambda^\epsilon \}_{\epsilon > 0}$,
there is $\bar y \in (-\varrho_2, \varrho_2)$ such that for any $y \in (\bar y - \kappa_6, \bar y + \kappa_6)$, for any $\omega \in \R / 2\pi \Z$,
there is a sub-interval $J_1 \subset [0,1)$ with $|J_1| >\kappa_6$ such that
\aryst
\inf_{s \in J_1}\|b\Psi_{x, \Phi^s_x(y)}(\Phi^\epsilon_x(s)) - \omega\|_{\R / 2\pi\Z} > \kappa_6.
\earyst
In particular, by \eqref{eq wsgx1 is long} and by letting $n_1$ be sufficiently large depending only on $g$ and $\kappa_6$, we can choose $v_1$ with $x^{v_1} \in W^s_{g}(x, (\bar y-\kappa_6, \bar y+\kappa_6))$, and a sub-interval $J_1 \subset [0,1)$ with $|J_1| >\kappa_6$ such that
\aryst
\inf_{s \in J_1}\|b\Psi_{x, x^{v_1}}(\Phi^\epsilon_x(s)) - \omega\|_{\R / 2\pi\Z} > \kappa_6
\earyst
where 
\aryst
\omega = \omega_{v_2}  - \omega_{v_1} + b\tau_{n_1} \circ v_2(x) - b\tau_{n_1} \circ v_1(x) \mod 2\pi \Z.
\earyst
Then for our choice of $J_1, v_1,v_2, \kappa_6$ we have for every $z \in \Phi^\epsilon_x(J_1) \subset \cP_0(x)$
\aryst
\| Arg((\tilde ue^{f^{(a,b)}_{n_1}} e^{ib\tau_{n_1}})(v_1(z))) - Arg((\tilde ue^{f^{(a,b)}_{n_1}} e^{ib\tau_{n_1}})(v_2(z))) \|_{\R / 2\pi\Z} > \frac{1}{2}\kappa_6.
\earyst
Let us assume that  
\aryst
\norm{(e^{f^{(a,b)}_n}H_n) \circ {v_1}}_{C^0(\cP_0(x))} \leq    \norm{(e^{f^{(a,b)}_n}H_n) \circ {v_2}}_{C^0(\cP_0(x))}.
\earyst
The other case is handled similarly.
 We may define $P_n$ on $\sigma^{-n_1}(\cP_0(x_0))$ by
\ary \label{eq defnofpn2}
P_n(z) = 
\begin{cases}
	1, & z \notin v_1(\Phi^\epsilon_x(J_1)), \\
	\zeta_{J_1}	((\Phi^\epsilon_x)^{-1}(\sigma^{n_1}(z))), & z \in v_1(\Phi^\epsilon_x(J_1)).
\end{cases}
\eary
where $\zeta_{J_1}$ is defined in analogy of $\zeta$: $\zeta_{J_1}$ equals $1-\kappa_5$ in the center half-interval of $J_1$; equals $1$ near the boundary; and has $C^1$ norm bounded by $10\kappa_5/\kappa_6$ (notice that $|J_1| \geq \kappa_6$).
As before, by requiring $n_1$ be sufficiently large depending on $g$ and $\kappa_6$,
there exists $y \in \cP_0(x)$ such that $(\Phi^\epsilon_x)^{-1}(\cP_1(y)) \subset J_1$ and  $d( \Phi^\epsilon_x)^{-1}(\cP_1(y)), \partial J_1) > \frac{1}{4}|J_1|$. Hence $P_n|_{v_1(\cP_1(y))} = 1-\kappa_5$.

By \eqref{eq ratio h v 1 v 2X} and by an elementary observation (see for instance \cite[Proposition 8]{Dol}), we have
for every $z \in \Phi^\epsilon_x(J_1)$ that
\aryst
&& |(\tilde ue^{f^{(a,b)}_{n_1}} e^{ib\tau_{n_1}})(v_1(z)) + (\tilde ue^{f^{(a,b)}_{n_1}} e^{ib\tau_{n_1}})(v_2(z))| \\
&\leq&  (1-\epsilon_5)(e^{f^{(a,b)}_n}H_n)(v_1(z)) + (e^{f^{(a,b)}_n}H_n)(v_2(z)).
\earyst
Consequently, we see that \eqref{eq n+1n3uniform bound} holds for all $x \in \cP_0(x_0)$.

We have defined $P_n$ on the interior of every atom of $\sigma^{-n_1}(\cP_0)$. Moreover, $P_n$ equals to $1$ near the boundary of every atom of $\sigma^{-n_1}(\cP_0)$. It suffices to let $P_n$ be $1$ on the completement of the interior of the union of all atoms of $\sigma^{-n_1}(\cP_0)$. It is clear that $P_n$ is a $C^1$ function.
We can conclude the proof by the following.
\begin{lemma} \label{claim pn in cone}
	By letting $n_1$ be sufficiently large depending on $g$ and $C_1$, 
	by letting $\kappa_5$  be sufficiently small depending on $n_1, C_1$ and $C_5$, 
	we have $P_n \in \cK_{\Lambda^\epsilon}$.
\end{lemma}
\begin{proof}
	We fix an arbitrary $\hat x \in U$.
	Let $x_0 \in U$ and $v \in \sigma^{-n_1}_{x_0}$
	satisfy that $v(\cP_0(x_0))$ intersects $W^u_{g}(\hat x, \Lambda^\epsilon(\hat x)^{-1})$. We set $$U' = v(\cP_0(x_0))  \cap W^u_{g}(\hat x, \Lambda^\epsilon(\hat x)^{-1}).$$ 
	
	If $\cP_0(x_0)$ is disjoint from $\Omega$, then by definition we have
	$P_n|_{U'} \equiv 1$. In particular, we have
	\ary \label{eq verifyingpnink1}
    \partial_z [\log P_n \circ \Phi^\epsilon_{\hat x}](z) = 0, \quad \forall z \in (\Phi^\epsilon_{\hat x})^{-1}(U').
	\eary
	
	If $\cP_0(x_0)$ intersects $\Omega$, then without loss of generality we may assume that $x_0 \in \Omega$. Then we let $x \in \cP_0(x_0)$ be given by the text above \eqref{eq choose x}.  In particular, \eqref{eq choose x} and \eqref{eq lambdaxlambdax0} hold, and the values of $P_n$ on $\sigma^{-n_1}(\cP_0(x_0))$ is given by \eqref{eq defnofpn1} if we are in case (1); and is given by \eqref{eq defnofpn2} if we are in case (2).   
	
	By \eqref{eq choose x}, \eqref{eq lambdaxlambdax0} and \eqref{eq local metric comparable},
	we have $\cP_0(x_0) \subset W^u_g(x_0, 4C_1\Lambda^\epsilon(x)^{-1})$. Then, by letting $n_1$ be sufficiently large depending on $g, C_1$, we have $v(\cP_0(x_0)) \subset W^u_g(v(x_0), \Lambda^\epsilon(x)^{-1})$. 
	
	We now show that 
	\ary\label{eq verifyingpnink2}
    |\partial_z [\log P_n \circ \Phi^\epsilon_{\hat x}](z)| \leq 1, \quad \forall z \in (\Phi^\epsilon_{\hat x})^{-1}(U').
    \eary
	We can assume that $\Lambda^\epsilon(\hat x) < \Lambda^\epsilon(x)$, for otherwise \eqref{eq verifyingpnink2} is clear by \eqref{eq defnofpn1}, \eqref{eq defnofpn2}, and by letting $\kappa_5$ be sufficiently small depending on $n_1$ and $C_1$.
	Then we have $\hat x \in W^u_g(v(x_0), 4\Lambda^\epsilon(\hat x)^{-1})$.
	By \eqref{eq lambdaxlambdax0} and by the hypothesis that $\Omega$ is $(n_1,  C_5)$-adapted to $\Lambda^\epsilon$, we have 
	\aryst
	\Lambda^\epsilon(x) < 2C_1\Lambda^\epsilon(x_0)  < 2C_1C_5 \Lambda^\epsilon(\hat x).
	\earyst
	Then by \eqref{eq defnofpn1}, \eqref{eq defnofpn2}, and by letting $\kappa_5$ be smaller if necessary, but depending only on $n_1, C_1$ and $C_5$, we obtain \eqref{eq verifyingpnink2}.
	
	Finally, notice that $P_n \in \cK_{\Lambda^\epsilon}$ follow from \eqref{eq verifyingpnink1} and \eqref{eq verifyingpnink2}.
\end{proof}

\section{Proof of Lemma \ref{prop exp rec in vol exp}}\label{sec proof of lemma prop exp rec in vol exp}
By Definition \ref{def recu}, there exist $C_4 > 1$ and $\eta_3 > 0$ such that
$\Omega$, and hence $\widetilde \Omega$ as well, are $(n_1, C_4, \eta_3)$-recurrent with respect to $\nu_U$.

For each $x \in U$, we denote the first return time and the first return map to $\widetilde \Omega$ by
\aryst
\tau_{\widetilde \Omega}(x) = \inf\{ n \geq 1 \mid \sigma^{nn_1}(x) \in \widetilde \Omega \}, \quad 
\sigma_{\widetilde \Omega} = \sigma^{\tau_{\widetilde \Omega}(x)}(x).
\earyst
We define
\aryst
R_0(x) = 0,  \quad R_k(x) = R_{k-1}(\sigma_{\widetilde \Omega}(x)) + \tau_{\widetilde \Omega}(x), \quad \forall k \geq 1, x \in U.
\earyst

Recall that $\cP_0, \cP_1$ are $\nu_U$-measurable finite partitions, and we have 
\aryst
\cP_1 =
\sigma^{-n_1}(\cP_0) \succeq \cP_0.
\earyst
We define $\cQ_0 = \cP_0|_{\widetilde \Omega}$ and
\aryst
\cQ_m = \sigma_{\widetilde \Omega}^{- m }(\cQ_0), \quad \forall m \geq 1.
\earyst
Clearly we have 
\aryst
\cQ_{n} \succeq \cQ_{n-1}, \quad \forall n \geq 1.
\earyst

We begin with a rather general lemma which only uses the fact that $\widetilde \Omega$ is a  $\cP_0$-measurable non-empty subset. 	Let $\{ \Omega_n \}_{n \geq 1}$ be given by the lemma. Then by definition,
for each $n \geq 1$, $\Omega_n$ is $\cQ_1$-measurable and $\cQ_0$-dense in $\widetilde \Omega$.
For any $m > 0$, any integer $l \geq 0$ and any $\eta > 0$, we set
\ary \label{eq def c n ell eta}
C_{m, l, \eta} &=& \{ x \in \widetilde \Omega \mid  |\{ 1 \leq j \leq m \mid \sigma_{\widetilde \Omega}^{j}(x) \in \Omega_{l + R_{j}(x)} \}| < \eta m \}.
\eary

\begin{lemma}  \label{lem upp bd for C}
	There is $\eta_4 > 0$ depending only on $g, F, n_1$ and $\Omega$ such that for any integer $\ell \geq 0$
	\aryst
	\nu_U( C_{m, l, \eta_4} ) < e^{-m \eta_4}, \quad \forall m > C'.
	\earyst
\end{lemma}
\begin{proof}
	We adopt the probabilistic notations introduced in the proof of Proposition \ref{prop L2 iteration}. 
	We define for each integer $l, k \geq 0$ that
	\aryst
	 \Omega_{l, k } = \{ x \in \widetilde \Omega \mid \sigma_{\widetilde \Omega}^{k }(x) \in \Omega_{l + R_{k }(x)} \}.
	\earyst
	\begin{lemma} \label{lem conditionmeasurelowerbound}
		There is $\eta_5 > 0$ depending only on $g, F$ and $n_1$  such that for any integers $l, k \geq 0$, 
		$\Omega_{l,k}$ is $\cQ_{k+1}$-measurable, and
		we have 
		\aryst
		\mathbb E( 1_{\Omega_{l, k }} \mid \cQ_k) \geq \eta_5.
		\earyst
	\end{lemma}
	\begin{proof}
		First notice that for any $x \in \widetilde \Omega$ and any $y \in \cQ_k(x)$, we have
		\aryst
		\tau_{\widetilde \Omega}( \sigma_{\widetilde \Omega}^{j}(x)) = \tau_{\widetilde \Omega}(\sigma_{ \widetilde \Omega}^{j}(y)), \quad 0 \leq j \leq k -1.
		\earyst
		Hence $R_{k }(x) = R_{k }(y)$.	
		
		Let $x \in \widetilde \Omega$ and let $y \in \cQ_{k+1}(x)$. By definition, we have
		\aryst
		\sigma_{ \widetilde \Omega}^{j}(y) \in \cQ_0(\sigma_{ \widetilde \Omega}^{j}(x)), \quad  0 \leq j \leq k+1. 
		\earyst
		In particular, we have
		\aryst
		\sigma_{ \widetilde \Omega}^{k}(y) \in  \cQ_1(\sigma_{ \widetilde \Omega}^{k}(x)).
		\earyst
		Since $\Omega_{l + R_{k}(x)}$ is $\cQ_1$-measurable, thus either both $\sigma^{k}_{\widetilde \Omega}(x)$ and $\sigma^{k}_{\widetilde \Omega}(y)$ belong to $\Omega_{l + R_{k}(x)} = \Omega_{l + R_{k}(y)}$, or  both $\sigma^{k}_{ \widetilde \Omega}(x)$ and $\sigma^{k}_{ \widetilde \Omega}(y)$ belong to $\Omega_{l + R_{k}(x)}^c = \Omega_{l + R_{k}(y)}^c$. This proves the $\cQ_{k+1}$-measurability of $\Omega_{l,k}$.
		
		To show the last inequality, notice that for a $\nu_U$-typical $x \in \widetilde \Omega$, we have
		\aryst
		\mathbb E( \Omega_{l, k } \mid \cQ_k)(x)  = \mathbb P (\cQ_k(x))^{-1} \mathbb P( \Omega_{l, k } \cap \cQ_k(x)).
		\earyst
		Since $\sigma_{ \widetilde \Omega}^{k}$ is invertible restricted to $\cQ_k(x)$,  by the Gibbs property of $\nu_U$ we have
		\ary 
		\mathbb P (\cQ_k(x))^{-1} \mathbb P( \Omega_{l, k } \cap \cQ_k(x))
		&\geq&	 (C')^{-1} \mathbb P (\sigma^{k}_{ \widetilde \Omega}\cQ_k(x))^{-1} \mathbb P(\sigma^{k}_{ \widetilde \Omega}( \Omega_{l, k } \cap \cQ_k(x))) \nonumber \\
		&=&  (C')^{-1} \mathbb P(\cQ_0(x'))^{-1} \mathbb P(\Omega_{n'} \cap \cQ_0(x')) \label{eq Omega n' over x'}
		\eary 
		where 
		\aryst
		x' = \sigma^{k}_{\widetilde \Omega}(x), \quad n' = l + R_{k}(x).
		\earyst
		By hypothesis, $\Omega_{n'}$ is $\cP_1$-measurable and $\cQ_0$-dense. Then by the  Gibbs property of $\nu_U$, we have $\eqref{eq Omega n' over x'} \geq \eta_5$ for some $\eta_5 > 0$ depending only on $g, F$ and $n_1$.
	\end{proof} 
	Fix an integer $l \geq 0$. We define the process $\{ X_n(x) \}_{n \geq 0}$ by
	\aryst
	X_j(x) = 1_{\sigma^{j}_{\widetilde \Omega}(x) \in \Omega_{l + R_{j}(x)}}(x), \quad \forall j \geq 0
	\earyst 
	where $x$ has distribution $\nu_U$.  By Lemma \ref{lem conditionmeasurelowerbound}, we see that $\{ X_n(x) \}_{n \geq 0}$ stochastically dominate an i.i.d. coin-flipping process with rate $\eta_5$. The proof then follows from the Large Deviation Principle for such process.
\end{proof}

\begin{proof}[Proof of Lemma \ref{prop exp rec in vol exp}]
	For any integer $L \geq 0$ and $\eta \in (0,1)$, we set
	\aryst
	A_{L, \eta} &=& \{ x \in U \mid |\{ 1 \leq j \leq L \mid \sigma^{jn_1}(x) \in \Omega_j \}| < \eta L \},  \\
	B_{L, \eta} &=& \{  x \in U \mid  |\{ 1 \leq j \leq L \mid \sigma^{jn_1}(x) \in \widetilde \Omega \}| < \eta L  \}.
	\earyst
	We set
	\ary \label{eq def eta1kappa3}
	\kappa_3 = \min( \eta_3, \eta_4\eta_3)/2.
	\eary	
	We now show that, by enlarging $C_4$ if necessary, for all sufficiently small $\epsilon > 0$ we have that
	\aryst
	\nu_U(A_{L, \kappa_3}) < e^{-L\kappa_3}, \quad \forall L \geq C_4.
	\earyst
	
	By the hypothesis that $\Omega$ is  $( n_1, C_4, \eta_3)$-recurrent with respect to $\nu_U$, we have
	\aryst
	\nu_U(B_{L, \eta_3}) < e^{-L\eta_3}, \quad \forall L > C_4.
	\earyst	
	For any $x \in B_{L, \eta_3}^c \cap A_{L, \eta_3\eta_4}$, there is an integer $0 \leq j \leq L-1$ such that
	\aryst
	|\{ 0 \leq i \leq \eta_3 L-1   \mid \sigma_{\widetilde \Omega}^{i} \sigma^{jn_1}(x) \in \Omega_{j + R_{i}(\sigma^{jn_1}(x)) } \}| < L \eta_3\eta_4.
	\earyst
	Consequently, by \eqref{eq def eta1kappa3} we have
	\ary \label{eq b c cap a}
	B_{L, \eta_3}^c \cap A_{L, \eta_3\eta_4} \subset \bigcup_{0 \leq j \leq L-1} \sigma^{-jn_1}(C_{\eta_3L, j, \eta_4}).
	\eary
	
	By \eqref{eq b c cap a}, Lemma \ref{lem upp bd for C} and by enlarging $C_4$ if necessary, we have for all $L > C_4$ that
	\aryst
	\nu_U(A_{L, \kappa_3}) &\leq& \nu_U(A_{L, \eta_3\eta_4}) \leq  \nu_U(B_{L, \eta_3}) + \nu_U(B_{L, \eta_3}^c \cap A_{L, \eta_3\eta_4} ) \\
	&<& 	e^{-L \eta_3} + \sum_{j = 0}^{L-1}  \nu_U( C_{\eta_3 L, j,  \eta_4	 }) \\
	&<&  e^{-L \eta_3} + L e^{-  \eta_4 \eta_3 L} < e^{-L\kappa_3}.
	\earyst	
\end{proof}

\appendix

\section{} \label{ap a}

\begin{proof}[Proof of Lemma \ref{lem subexponential growth on average}]
	
	Let $\gamma_0 > 0$ be a small constant to be determined later. 
	Let $\cP = \{U_\a\}_{\a \in \cA}$ be the natural partition of $U$.
	For each integer $n \geq 1$, we set 
	\aryst
	\cP_n = \sigma^{-n}(\cP),
	\earyst
	and let $S_n: U \to \R$ be a function defined $\nu_U$-almost everywhere by
	\aryst
	S_n(x) =\sup_{y \in W^s_{g,loc}(w)} \sup_{w \in \cP_n(x)} \det (Dg^{\tau_n(y)}(y))^{\gamma_0}.
	\earyst
	Define
	\aryst
	M_n = \sup_{\a \in \cA} \nu_U(U_\a)^{-1}\int_{U_\a}  S_n(x)   d\nu_U(x).
	\earyst	
	Clearly, we have
	\ary \label{eq mnasupperbound}
	\int \det(D g^{\tau_n})^{\gamma_0} d\nu_{U}  \leq \int_U S_n(x) d\nu_U(x) \leq M_n.
	\eary
	By \eqref{eq wsgx1 is long}, we have for all sufficiently large $n$, and for every $x \in U$ that
	\aryst
	\hat \sigma(W^s_{g, loc}(x)) \subset W^s_{g, loc}(\sigma(x)).
	\earyst
	Then by definition, we have for any sufficiently large integer $n \geq 1$, any integer $k \geq 1$, and $\nu_U$-a.e. $x \in U$
	\aryst
	S_{kn}(x) \leq S_{(k-1)n}(x) S_n(\sigma^{(k-1)n}(x)).
	\earyst
	Then
	\aryst
	M_{kn} &\leq& C' \int S_{(k-1)n}(x) S_n(\sigma^{(k-1)n}(x)) d\nu_U(x) \\
	&=& C' \int Q_{(k-1)n}(y) S_n(y) d\nu_U(y)
	\earyst
	where we denote for any integer $m \geq 1$ that
	\aryst
	Q_m(x) = \sideset{}{^*}\sum_{v \in \sigma^{-m}} e^{\hat f_m(v(x))}S_m(v(x)), \quad \forall  x \in U.
	\earyst
    The notation $\sideset{}{^*}\sum$ here is defined in Section \ref{sec: bounds for smoothness}.
	
	Notice that for every $v \in \sigma^{-m}$, $S_m \circ v$ is constant on each $U_\a$.
	Then by distortion estimate, we see that for every $\a \in \cA$ and any $x, y \in U_\a$
	\aryst
	Q_m(x) \sim Q_m(y).
	\earyst
	Thus for all sufficiently large integer $n \geq 1$ we have
	\aryst
    M_{kn} &\leq& C' \sum_{\a \in U} Q_{(k-1)n}(x_\a) \int_{U_\a} S_n(y) d\nu_U(y) \\
    &\leq& C' \sum_{\a \in U} Q_{(k-1)n}(x_\a) \nu_U(U_\a) M_n \\
    &\leq& C' \sum_{\a \in U} \int_{U_\a}Q_{(k-1)n}(y) d\nu_U(y) M_n \\
    &\leq& C' M_{(k-1)n}M_n.
	\earyst
	Thus for all sufficiently large integer $n \geq 1$, and any integer $ k \geq 1$, we have
	\ary \label{eq submultiplicativemn}
	M_{kn} \leq (C' M_n)^k.
	\eary	
	We also have
	\aryst
	\log \det (Dg^{\tau_n(x)}(x)) = \int_0^{\tau_n(x)} \rdiv V_g(g^t(x)) dt.
	\earyst
	By distortion estimate, we have for any integer $n \geq 1$ and for every $x \in U$ that
	\ary \label{eq distortion for int of div}
	|\int_0^{\tau_n(x)} \rdiv V_g(g^t(x)) dt - \log S_n(x)| < C'
	\eary
	Thus for all $n \geq 1$, for all $\gamma_0 \in (0, n^{-1}]$,  we have
	\aryst
	M_n &\leq& \int C' \exp( \gamma_0 \int_0^{\tau_n(x)} \rdiv V_g(g^t(x))  dt ) d\nu_U(x) \\
	&\leq& C' (1 + \gamma_0 \int \int_0^{\tau_n(x)} \rdiv V_g(g^t(x)) dt d\nu_U(x) + O( \gamma_0^2n^2)).
	\earyst
	By \eqref{eq distortion for int of div},  we have
	\aryst
	|\int \int_0^{\tau_n(x)} \rdiv V_g(g^t(x)) dt d\nu_U(x) - \int \int_0^{\tau_n(w)} \rdiv V_g(g^t(w)) dt d\nu_\Pi(w)| < C'.
	\earyst
	On the other hand, we have
	\aryst
	\int \int_0^{\tau_n(w)} \rdiv V_g(g^t(w)) dt d\nu_\Pi(w) = n \int \rdiv V_g(x) d\nu_{g,F}(x) \leq 0.
	\earyst
	Consequently, for any integer $n \geq 1$, we have for any $\gamma_0 \in (0, n^{-1}]$ that
	\aryst
	M_n \leq C'(1 + \gamma_0 C' + O( \gamma_0^2 n^2)).
	\earyst
	By letting $m_0$ be a large integer, and by setting $\gamma_0 = \frac{1}{m_0}$, we obtain
	\aryst
	M_{n_1} \leq C'.
	\earyst
	Then by letting $m_0$ be sufficiently large depending on $\kappa_0$, by letting $\gamma_0 = \frac{1}{m_0}$, by letting $n_0$ be sufficiently large depending on $g, n_1, \kappa_0$, and by \eqref{eq submultiplicativemn}, we have for any $n > n_0$ that
	\aryst
	M_n \leq e^{n\kappa_0\gamma_0}.
	\earyst
	This concludes the proof by \eqref{eq mnasupperbound}.
\end{proof}

\section{} \label{appendix chart normalization}

\begin{proof}[Proof of Lemma \ref{lem normalcoordinatesystem}]

Given any $x \in M$ and $m \in \Z$ we will denote for simplicity that
\aryst
x_m = g^m(x), \quad \check g_m = \check g_{x_m},    \quad \mu_m = \mu_{x_m}, \quad \lambda_m = \lambda_{x_m}.
\earyst
We define $\check \psi_m$, $\check f_m$, $\check f_{m,1}$,$\check f_{m,2}$ analogously.
We fix an arbitrary $\rhone > 0$, and denote 
\aryst
\rhone' = D\rhone
\earyst
for some large constant $D > 1$ to be determined in due course depending only on $g$. We will let $\varrho_0$ be sufficiently large depending on $\varrho'_1$; and
let $\upsilon_*$ in Section \ref{sec basicpropertiesofanosovflow} be sufficiently small, depending on $g, \varrho_0, D$, so that the range of any chart $\iota_x$ is contained in a small ball of $M$.
\begin{lemma} \label{lemma chart normalization}
	For any integer $K  > \frac{\chi_*}{\chi_0}$ the following is true. 
    For any $x \in M$,
    there are  functions $\tilde \rho_x$, $\tilde \xi_x \in C^r(-\rhone', \rhone')$ with $\norm{\tilde \rho_x}_{C^r}, \norm{\tilde \xi_x}_{C^r} < C'$ such that by setting $\tilde g_x = \tilde h_{g^1(x)}^{-1} \check g_x \tilde h_x$ where $\tilde h_x: (-\rhone', \rhone')^3 \to (-C'\rhone', C'\rhone')^3$ is a $C^r$ embedding of form
	\ary 
	 \tilde h_x(z,y,t) = (z, e^{\tilde\rho_x(z)}y, t + \tilde\xi_x(z)y),
	\eary   
	and write 
	\aryst
	\tilde g_x(z,y,t) = (\tilde f_{x,1}(z,y), \tilde f_{x,2}(z,y), t + \tilde \psi_{x}(z,y)),
	\earyst
	we have
	\aryst
	\partial_y\tilde f_{x,2}(\cdot,0) \equiv \mu_x^{-1} \quad \mbox{and} \quad \partial_y\tilde\psi_x(\cdot,0) \in Poly^K. 
	\earyst
\end{lemma}
\begin{proof}[Proof of Lemma \ref{lemma chart normalization}]
Given $x \in M$, we denote for every $z \in (-\rhone',\rhone')$ that
\aryst
z_m &=& \lambda_{m} \cdots \lambda_{-1}z,  \quad \forall m < 0.
\earyst
We set
\ary
\ \ \tilde\rho_x(z) &=& \sum_{n \geq 1}(\log \partial_y \check f_{-n, 2}(z_{-n}, 0) - \log  \mu_{-n}^{-1} ). \label{eq def tilderho}
\eary
By definition, we see that $\tilde \rho_x$ is uniformly $C^r$, and $\tilde\rho_x(0) = 0$.
For any $x \in M$ and $z \in (-\varrho'_1, \varrho'_1)$ we set
\ary
\ \ \ \ 
\tilde\xi_x(z) &=& \int_0^z \frac{(z-w)^K}{K!} \sum_{n \geq 1}   (\prod_{i=1}^{n}\mu_{-i}\lambda_{-i}^K) \partial^K_z( e^{-\tilde \rho_{x_{-n}}(w_{-n})} \partial_y \check\psi_{-n}(w_{-n}, 0) )  dw. \label{eq def tildexix} 
\eary
Then it is direct to verify that
\aryst
\norm{\tilde \rho_x}_{C^r}, \norm{\tilde \xi_x}_{C^r} < C', \quad \forall x \in M.
\earyst

As $\tilde g_x = \tilde h_{g^1(x)}^{-1} \check g_x \tilde h_x$, we have
\ary
\tilde f_{x,1}(z,y) &=& \check f_{x,1}(z, e^{\tilde \rho_x(z)}y),  \\
\tilde f_{x,2}(z,y) &=& e^{-\tilde \rho_{g^1(x)}(\check  f_{x,1}(z, e^{\tilde \rho_x(z)}y))}\check  f_{x,2}(z, e^{\tilde \rho_x(z)}y),  \label{eq tildefx2formula} \\
\tilde \psi_x(z,y) &=& \check  \psi_x(z, e^{\tilde \rho_x(z)}y) + \tilde \xi_x(z) y - \tilde \xi_{g^1(x)}(\check  f_{x,1}(z, e^{\tilde \rho_x(z)}y))  \label{eq tildepsiformula} \\
&&\cdot e^{-\tilde \rho_{g^1(x)}(\check  f_{x,1}(z, e^{\tilde \rho_x(z)}y))}\check  f_{x,2}(z, e^{\tilde \rho_x(z)}y). \nonumber
\eary
Differentiate \eqref{eq tildefx2formula} with respect to $y$ and evaluate at $y = 0$, we obtain
\ary \label{eq normalize derivatives}
\partial_y \tilde f_{x,2}(z,0) = \partial_y \check  f_{x,2}(z,0) e^{\tilde \rho_x(z) - \tilde \rho_{g^1(x)}(\lambda_x^{-1}z)}.
\eary
By \eqref{eq def tilderho}, the right hand side of \eqref{eq normalize derivatives} equals $\mu_x^{-1}$. This shows the first equality of the lemma.

Differentiate \eqref{eq tildepsiformula} with respect to $y$ and evaluate at $y = 0$, we obtain
\aryst
\partial_y \tilde \psi_x(z,0) &=& e^{\tilde\rho_x(z)} \partial_y \check  \psi_x(z,0) + \tilde\xi_x(z) - \tilde \xi_{g^1(x)}(\lambda_x^{-1}z) e^{\tilde\rho_x(z) - \tilde\rho_{g^1(x)}(\lambda_x^{-1}z)} \partial_y \check  f_{x,2}(z,0) \\
&=&e^{\tilde\rho_x(z)} \partial_y \check  \psi_x(z,0) + \tilde\xi_x(z) - \mu_x^{-1}\tilde \xi_{g^1(x)}(\lambda_x^{-1}z).
\earyst
Differentiate the above equality $K$ times with respect to $z$, replace $x$ by $g^{-1}(x)$, then replace $z$ by $\lambda_{g^{-1}(z)}z$, we obtain
\aryst
\partial^K_z\partial_y \tilde \psi_x(z,0) = \partial^K_z(e^{\tilde\rho_x(z)} \partial_y \check  \psi_x(z,0)) + \partial^K_z \tilde\xi_x(z) - \lambda_x^{-K}\mu_x^{-1}\partial^K_z\tilde \xi_{g^1(x)}(\lambda_x^{-1}z).
\earyst
Substitute \eqref{eq def tildexix} into the above equality, we obtain
\aryst
\partial^K_z \partial_y \tilde \psi_x(z,0) = 0, \quad \forall z \in (-\varrho'_1, \varrho'_1).
\earyst
In another words, $\partial_y \tilde\psi_x(\cdot, 0) \in Poly^K$.
\end{proof}

%We apply Lemma \ref{lemma chart normalization} to obtain $\tilde h_x$ for each $x \in M$. We replace $\iota_x$ and $g_x$ by $\iota_x \tilde h_x$ and $\tilde g_x = \tilde h_{g(x)}^{-1} g_x \tilde h_x$ respectively for each $x \in M$.

Similar to Lemma \ref{lemma chart normalization},  by switching $z,y$ coordinates, we have the following statement for $\tilde g_x$ obtained in Lemma \ref{lemma chart normalization}. We omit its proof.
\begin{lemma}  \label{lemma chart normalization 2}
	For any integer $K > \frac{\chi_*}{\chi_0}$ the following is true. 
	For any $x \in M$, there are  functions $\hat \rho_x, \hat \xi_x \in C^r(-\rhone, \rhone)$ with $\norm{\hat \rho_x}_{C^r}, \norm{\hat \xi_x}_{C^r} < C'$ such that by setting $\bar g_x = \hat h_{g^1(x)}^{-1} \tilde g_x \hat h_x $ where $\hat h_x: (-2\rhone, 2\rhone)^3 \to (-\rhone', \rhone')^3$ is a $C^r$ embedding of form
	\aryst
    \hat h_x(z,y,t) &=& (e^{\hat \rho_x(y)}z, y, t + \hat \xi_x(y)z),
	\earyst   
	and write 
	\aryst
	 \bar g_x(z,y,t) = ( \bar f_{x,1}(z,y),  \bar f_{x,2}(z,y), t + \bar \psi_{x}(z,y)),
	\earyst
	we have
	\aryst
	\partial_z \bar f_{x,1}(0,\cdot) = \lambda_x^{-1} \quad \mbox{and} \quad \partial_z \bar \psi_x(0,\cdot) \in Poly^K.
	\earyst
\end{lemma}
To finish the proof of Lemma \ref{lem normalcoordinatesystem}, it suffices to take $\varrho_0 > C'\rhone'$, and take
\aryst
\check  h_x = \tilde h_x \circ \hat h_x, \quad \forall x \in M.
\earyst

\end{proof}

\section{} \label{ap c}

\begin{proof}[Proof of Lemma \ref{lem def of gammax}]
Denote by $\pi_{1,2}$, resp. $\pi_1$, the projection of $(-\rOut, \rOut)^3$ to the first two coordinates, resp. the first coordinate.
For each $y \in (-2\varrho_1, 2\varrho_1)$, we set
\aryst
\Gamma_{x,y}  = \pi_{1,2}(\tWu_x(0, y, 0)).
\earyst
The parametrization $z \mapsto \Phi^u_{\iota_x(0,y,0)}(z)$ gives rise to the vector field  $D\Phi^u_{\iota_x(0,y,0)}(\partial_z)$ on $W^u_{x}(\iota_x(0,y,0))$.
We define a vector field $\Omega_{x,y}$ on each $\Gamma_{x,y}$ by
\aryst
\Omega_{x,y} = D(\pi_{1,2}\iota_x^{-1})D\Phi^u_{\iota_x(0,y,0)}(\partial_z).
\earyst
 Denote by $[\Omega_{x,y}]$ the set of vector fields which are positively proportional to $\Omega_{x,y}$. It is direct to verify that 
 \ary \label{eq dfxomegaxyinvariant}
 [Df_x(\Omega_{x,y})] = [\Omega_{g^1(x), \mu_x^{-1}y}].
 \eary
Let $U_{x, y} \in [\Omega_{x,y}]$ be the unique vector field such that
for each $y \in (-2\varrho_1, 2\varrho_1)$, we have
\ary \label{eq pi1u}
\pi	_1U_{x,y}(0,y) =  1, \quad \forall y \in (-2\varrho_1, 2\varrho_1).
\eary

We define 
\aryst
\gamma_x(z,y) = \phi^{z}_{U_{x,y}}(0,y), \quad \forall y,z \in (-\varrho_1, \varrho_1)
\earyst
where $\phi^t_{U_{x,y}}$ denote the local flow generated by $U_{x,y}$. Clearly, we have item (1).

By \eqref{eq dfxomegaxyinvariant}, \eqref{eq pi1u} and the equality
\aryst
\partial_z f_{x,1}(0,y) = \lambda_x^{-1}, \quad \forall y \in (-2\varrho_1, 2\varrho_1),
\earyst 
we can verify that for any $z \in (-2|\lambda_x|\varrho_1, 2|\lambda_x|\varrho_1)$ and $y \in (-2\varrho_1, 2\varrho_1)$
\aryst
	f_x(\gamma_x(z, y)) &=& \gamma_{g^1(x)}( \lambda_x^{-1}z ,\mu_x^{-1}y).
\earyst
In particular we have 
\ary \label{eq cygxandcyx}
\cY_{g^1(x)}(\lambda_x^{-1}z, \mu_x^{-1}y) = f_{x,2}(\cX_x(z,y), \cY_x(z,y)).
\eary

We can verify \eqref{eq straighten 4} in item (3) by the H\"older continuous dependence of the non-stationary parametrizations $\Phi^u_x$ on $x$ (here we tacitly identify two parametrisations which only differ in orientations).

To prove item (3),
we first show that $\partial_y\cY_x(z,0)$ exists for all $z \in (-\varrho_1, \varrho_1)$.
Notice that
\ary \label{eq derivative y 1}
f_{x,2}(z,0) = 0, \quad \forall z \in (-2\varrho_1, 2\varrho_1). 
\eary
Since $W^{cu}$ is a $C^{1+\theta}$-foliation with $C^r$ leaves,
we deduce that $\cF_x = \{  \Gamma_{x,y}  \}_{y \in (-\varrho_1, \varrho_1)}$ is also a $C^{1+\theta}$-foliation of a subset $\Gamma_x \subset (-2\varrho_1, 2\varrho_1)^2$ with $C^r$ leaves. 
Moreover, by letting $\nu_*$ in Section \ref{sec basicpropertiesofanosovflow} be sufficiently small, we can see that for every $(z,y) \in (-\varrho_1, \varrho_1)^2$, there is a well-defined map $H_{x,z}: (-\varrho_1, \varrho_1) \to (-2\varrho_1, 2\varrho_1)$ given by equation
\aryst
(z, H_{x,z}(y)) = \{z\} \times (-2\varrho_1, 2\varrho_1) \cap \Gamma_{x,y}.
\earyst
By definition, we have
\ary \label{eq yequalsh}
\cY_x(z,y) = H_{x, \cX_x(z,y)}(y).
\eary
As the holonomy map along $\cF_x$ from $\{0 \} \times (-\varrho_1, \varrho_1)$ to $\{z\} \times (-2\varrho_1,2\varrho_1)$ is $C^{1+\theta}$ with  $C^{1+\theta}$ norm uniformly bounded in $x,y,z$, we have the following:
\enmt
\item for all $z \in (-\varrho_1, \varrho_1)$, $y \in (-\varrho_1/C', \varrho_1/C')$, we have
\ary \label{eq derivative y 2}
|y|/C' < |\cY_x(z,y)| < C' |y|;
\eary
\item for any $\kappa > 0$, there exists $\eta > 0$ such that 
\ary \label{eq derivative y 3}
(1-\kappa)|y| < |\cY_x(z,y)| < (1+\kappa) |y|, \quad \forall z \in (-\eta, \eta);
\eary
\item for any $z \in (-\varrho_1, \varrho_1)$, we have
\ary \label{eq derivative y 4}
\norm{H_{x,z}}_{C^{1+\theta}} < C'.
\eary
\eenmt
 Then by \eqref{eq cygxandcyx}, \eqref{eq derivative y 1}, \eqref{eq derivative y 2}, \eqref{eq straighten 4} and the Taylor expansion of $f_{x,2}$ at $(z,0)$,
we obtain
\aryst
\cY_{g^1(x)}(\lambda_x^{-1}z, \mu_x^{-1}y) = \partial_y f_{x,2}(z,0) \cY_x(z,y) + O(|y|^{1+\delta}).
\earyst
Then by \eqref{eq straignten 1}, we have
\ary
 (\mu_x^{-1}y)^{-1}\cY_{g^1(x)}(\lambda_x^{-1}z, \mu_x^{-1}y) 
&=&  \mu_x \partial_y f_{x,2}(z,0) y^{-1}\cY_x(z,y) + O(|y|^{\delta}) \nonumber \\
&=& y^{-1}\cY_x(z,y)  + O(|y|^{\delta}). \label{eq fractionrelation}
\eary
For each $x \in M$, we define a function $w_x: (-\varrho_1, \varrho_1) \to \R$  by
\aryst
w_x(z) = \liminf_{y \to 0}  y^{-1}\cY_x(z,y), \quad \forall z \in (-\varrho_1, \varrho_1).
\earyst
By \eqref{eq fractionrelation}, for all $x \in M$ we have
\aryst
w_{g^1(x)}(z) = w_x(\lambda_x z), \quad \forall z \in (-\varrho_1, \varrho_1).
\earyst
Thus we have
\aryst
w_x(z) = \lim_{n \to \infty} w_{x_{-n}}( \lambda_{-n} \cdots \lambda_{-1} z), \quad \forall  z \in (-\varrho_1, \varrho_1).
\earyst
On the other hand, by \eqref{eq derivative y 3} and $\cY_x(0,y) = y$ we have
\aryst
\lim_{z \to 0} w_x(z) = w_x(0) = 1,
\earyst
and the convergence of the above limit is uniform for all $x \in M$.
Thus 
\aryst
w_x(z) = 1, \quad \forall z \in (-\varrho_1, \varrho_1).
\earyst
Similarly, we can show that
\aryst
\limsup_{y \to 0}  y^{-1}\cY_x(z,y) = 1, \quad \forall z \in (-\varrho_1, \varrho_1).
\earyst
Consequently,
\aryst
\lim_{y \to 0}  y^{-1}\cY_x(z,y) = 1, \quad \forall z \in (-\varrho_1, \varrho_1).
\earyst
Now it is straightforward to show that 
\aryst
DH_{x, z}(0) = 1, \quad \forall z \in (-\varrho_1, \varrho_1).
\earyst
Lemma \ref{lem def of gammax}(3) follows from \eqref{eq yequalsh} and \eqref{eq derivative y 4}.

 It remains to show Lemma \ref{lem def of gammax}(2). 
 %Let charts $\{  \iota_x \}_{x \in M}$ be given by Section \ref{sec constructionofcharts}. 
 Let $C > 1$ be a large constant to be determined later. Given $y \in (-\varrho_1, \varrho_1)$. Assume that function $w: (-\varrho_1, \varrho_1) \to \R$ satisfies
\aryst
w(z) = O(C|y|)z +  \cdots + O(C|y|)z^{r-1} + O( C z^{r}).  
\earyst
Let $n \geq 1$ be an integer to be determined.  We abbreviate $f_{g^m(x)}$, $f_{g^m(x), i}$ as $f_{m}$, $f_{m,i}$ respectively. We denote  for any $0 \leq m \leq n$ that
\aryst
F_m = (F_{m,1}, F_{m,2}) &=& f_{m-1} \cdots f_0, \\
\mu_{0, m} = \mu_0 \cdots \mu_{m-1}, &\quad& \lambda_{0, m} = \lambda_0 \cdots \lambda_{m-1}.
\earyst
We let
\aryst
Z(z) &=& F_{n,1}(z, y + w(z)), \\
Y(z) &=& F_{n,2}(z, y + w(z)).
\earyst

Notice that we have 
\aryst
F_{n,2}(z,0) = 0, \quad \forall z \in (-\varrho_1, \varrho_1).
\earyst
Then we have
\aryst
\partial^i_z F_{n,2}(0,0) = 0, \quad \forall 1 \leq i \leq r.
\earyst
Thus
\aryst
\partial^i_z F_{n,2}(0,y) = O(C_n|y|), \quad \forall 1 \leq i \leq r
\earyst
where we use $C_n$ to denote a constant depending only on $n$ and $r$ (we suppress the dependence on $r$ in notation) in this section.
It is direct to show by induction that
\aryst
\partial_y F_{n,2}(z,0) = \mu_{0,n}^{-1}.
\earyst
Hence by a similar argument as above we obtain
\aryst
\partial_z^i\partial_y F_{n,2}(0,y) = O(C_n |y|), \quad \forall 1 \leq i \leq r-2.
\earyst
By Taylor expansion of $F_{n,2}$ at $(0,y)$, we obtain
\aryst
Y(z) &=& \mu_{0,n}^{-1} y + \sum_{i = 1}^{r-1} \partial_z^i F_{n,2}(0,y) z^i +  \partial_y F_{n,2}(0,y)  w(z) +
 \sum_{i = 1}^{r-2} \partial_z^i \partial_y F_{n,2}(0,y) z^i w(z) \\
 && +
 \sum_{i=1}^{r-3} \sum_{j=2}^{r-i-1} \partial_z^i \partial^j_y F_{n,2}(0,y) z^i w(z)^j  +  O(C_n |z|^r). 
\earyst

Hence we have
\aryst
Y(z) &=& \mu_{0,n}^{-1} y + \sum_{i=1}^{r-1} O( C_{n} |y|) z^i  + 
 \sum_{i=1}^{r-1} O(|\mu_{0,n}^{-1}|C|y|) z^i \\
 && + \sum_{i=2}^{r-1} O(CC_n|y|^2) z^i + O((C_n + |\mu_{0,n}^{-1}| C + CC_n|y|)|z|^r) \\
 &=& \mu_{0,n}^{-1} y+  O(C_n|y| + |\mu_{0,n}^{-1}|C|y|) z  + \sum_{i=2}^{r-1} O( |\mu_{0,n}^{-1}|C|y| + C_{n} |y|+CC_n|y|^2) z^i  \\
 &&  + O((C_n + |\mu_{0,n}^{-1}| C + CC_n|y|)|z|^r).
\earyst
By letting $C$ be sufficiently large depending on $g,n$, and by letting $|y|$ be sufficently small depending on $g,n$, we have
\ary \label{eq Y to z}
Y(z) = \mu_{0,n}^{-1} y  + \sum_{i=1}^{r-1} O( |\mu_{0,n}^{-1}| C |y| ) z^i    + O(|\mu_{0,n}^{-1}| C |z|^r).
\eary
It is straightforward to see that
\ary \label{eq z to Z}
z = O(|\lambda_{0,n}|)Z + \cdots + O(|\lambda_{0,n}|)Z^{r-1} + O(C_n |Z|^r).
\eary
By  substituting \eqref{eq z to Z} into \eqref{eq Y to z}, and by letting $|y|$ be sufficiently small,  we obtain
\aryst
Y(z) = \mu_{0,n}^{-1} y  + \sum_{i=1}^{r-1} O( |\mu_{0,n}^{-1}y|C \lambda_{0,n} ) Z^i    + O((|\mu_{0,n}^{-1}| C + C_n) |Z|^r).
\earyst 
By letting $n$ be sufficiently large, and 
by letting $C$ be sufficiently large depending on $g, n$, we have
\aryst
Y = \mu_{0,n}^{-1}y + O(|\mu_{0,n}^{-1}y| C) Z + \cdots +O(|\mu_{0,n}^{-1}y| C) Z^{r-1} +O(C |Z|^r).
\earyst
We then conclude the proof of Lemma \ref{lem def of gammax} following the proof of \cite[Prop 2.1]{Ha} using the above estimate instead of \cite[Lemma 2.2]{Ha}.
\end{proof}

\section{} \label{ap d}

\begin{proof}[Proof of Lemma \ref{lem tpltperp}]
	The existence of $\xi^{u, \perp}_x$ follows from Lemma \ref{lemma chart normalization}. Indeed, for each $x \in M$, we define $\widetilde \tplt^\perp_x \in \Tplt^{u,r}_0(x)$ by \eqref{eq tildexiperpxtilde0x}.
	More explicitly, we have for any $z \in (-\varrho_1, \varrho_1)$ that
	\aryst
	 \langle \widetilde \xi^{\perp}_x(\Phi^u_x(z)), (D\iota_x)_{(z,0,0)}(a,b,c)  \rangle = b.
	\earyst
	Then it is direct to verify that $\{ \widetilde \xi^{\perp}_x \}_{x \in M}$ satisfies the required properties.
	
	To show the uniqueness, we notice that any $\{  \xi^{u, \perp}_x \}_{x \in M}$ in the lemma can be expressed as
	\aryst
	\xi^{u, \perp}_x(\Phi^u_x(z)) = 
	\widetilde\varphi_x(z) \widetilde \xi^{\perp}_x(\Phi^u_x(z)),
	\earyst
	where $\{ \widetilde\varphi_x \}_{x \in M}$ is a bounded subset of $C^r(-\varrho_1, \varrho_1)$. Moreover, we have the following properties:
	\enmt
	\item for any $x \in M$, we have
	\aryst
	\widetilde \varphi_x(0) = \pm 1;
	\earyst
	\item for any $x \in M$, any $n > 0$, we have
	\aryst
	\widetilde \varphi_x(z) = \widetilde \varphi_{g^{-n}(x)}( \lambda_{g^{-n}(z)} \cdots \lambda_{g^{-1}(z)} z), \quad \forall z \in (- \varrho_1, \varrho_1).
	\earyst
	\eenmt
	Then it is easy to conclude that for any $x \in M$, we have 
	\aryst
	\widetilde \varphi_x \equiv \pm 1.
	\earyst
\end{proof}

\begin{proof}[Proof of Lemma \ref{lem choose straight sections}]
For each $x \in M$, we define $\widetilde \xi^0_x \in \Tplt^{u,r}_1(x)$ by \eqref{eq tildexiperpxtilde0x2}.
More explicitly, we have for any $z \in (-\varrho_1, \varrho_1)$ that
	\aryst
	\langle \widetilde\xi^{0}_x(\Phi^u_x(z)), (D\iota_x)_{(z,0,0)}(a,b,c)  \rangle = c.
	\earyst
	
Take an arbitrary $x \in M$.	Let $\varphi \in C^r(-\varrho_1, \varrho_1)$ and set
	\aryst
	\xi(\Phi^u_x(s)) = \widetilde\xi^{0}_x(\Phi^u_x(s)) + \varphi(s)\xi^{u,\perp}_x(\Phi^u_x(s)), \quad \forall s \in (-\varrho_1, \varrho_1). 
	\earyst
	
	Let $\{ \xi_n \in \Tplt^{u,r}(g^{-n}(x)) \}_{n \geq 1}$ be a sequence of sections such that for all $n \geq 1$, we have $L^n_{g^{-n}(x)} \tplt_n = \tplt$. We set
	\aryst
	\varphi_n = \varphi^u_{g^{-n}(x), \tplt_n}, \quad \forall n \geq 1.
	\earyst
	Then by \eqref{eq straignten 1} we have
	\aryst
	\varphi \in  \pm \clb \mu_{-n} \cdots \mu_{-1} \varphi_n(\lambda_{-1} \cdots \lambda_{-n}\cdot ) + Poly^K, \quad \forall n \geq 1.
	\earyst
	Thus
	\ary \label{eq dkformula}
	D^K\varphi =  \pm \clb \mu_{-n} \cdots \mu_{-1}(\lambda_{-1} \cdots \lambda_{-n})^K D^K\varphi_n(\lambda_{-1} \cdots \lambda_{-n}\cdot ).
	\eary
	
	If $\varphi \in Poly^K$, then by \eqref{eq dkformula}, for any $n \geq 1$, the restriction of $\varphi_n$ to interval $(-|\lambda_{-1} \cdots \lambda_{-n}|\varrho_1, |\lambda_{-1} \cdots \lambda_{-n}|\varrho_1)$ belongs to $Poly^K$. 
	
	Conversely, if $\norm{D^K\varphi_{n}}_{(-\varrho_1, \varrho_1)}$ is uniformly bounded in $n$, then by \eqref{eq dkformula} and by $K > \frac{\chi_*}{\chi_0}$, we see that $D^K \varphi = 0$. 
\end{proof}

\end{document}